\documentclass{amsart}
\usepackage{amssymb,amsmath, amsthm,latexsym}
\usepackage{graphics}
\usepackage{amscd}
\usepackage{graphics}
\usepackage{here}
\newcommand{\cal}[1]{\mathcal{#1}}
\theoremstyle{plain}

\newtheorem{lemma}{Lemma}[section]

\newtheorem{proposition}[lemma]{Proposition}
\newtheorem{corollary}[lemma]{Corollary}
\parskip=\bigskipamount

\let\egthree=\phi
\let\phi=\varphi
\let\varphi=\egthree

\begin{document}
\title{Geometry of the mapping class groups III: Quasi-isometric
rigidity}
\author{Ursula Hamenst\"adt}
\thanks
{
Research
partially supported by Sonderforschungsbereich 611}
\date{January 28, 2007}

\begin{abstract}
Let $S$ be an oriented surface of genus
$g\geq 0$ with $m\geq 0$ punctures and $3g-3+m\geq 2$.
We show that for every finitely generated group $\Gamma$
which is quasi-isometric to the mapping class group
${\cal M}(S)$ of $S$ there is a finite index subgroup
$\Gamma^\prime$ of $\Gamma$ and
a homomorphism
$\rho:\Gamma^\prime\to {\cal M}(S)$ with finite kernel
and finite index image. We also give a new proof of
the following result of Behrstock and Minsky:
The geometric rank of ${\cal M}(S)$ as well as the
homological dimension of the asymptotic cone of
${\cal M}(S)$ equal $3g-3+m$.
\end{abstract}

\maketitle

\tableofcontents

\section{Introduction}

Let $S$ be an oriented surface of finite type, i.e. $S$ is a
closed surface of genus $g\geq 0$ from which $m\geq 0$
points, so-called \emph{punctures},
have been deleted. We assume that $3g-3+m\geq 2$,
i.e. that $S$ is not a sphere with at most $4$
punctures or a torus with at most $1$ puncture.
We then call the surface $S$ \emph{non-exceptional}.
The \emph{mapping
class group} ${\cal M}(S)$ of all isotopy classes of
orientation preserving self-homeomorphisms of $S$ is finitely
presented \cite{I02}, indeed it acts as a group of automorphisms
on a contractible cell complex
with finite cell stabilizers and compact
quotient. In particular, it
is finitely generated. We refer to the survey of
Ivanov \cite{I02} for more about the mapping class
group and for references.

For a number $L\geq 1$, an \emph{$L$-quasi-isometric embedding}
of a metric space $(X,d)$ into a metric space
$(Y,d)$ is a map $F:X\to Y$ which satisfies
\[d(x,y)/L-L\leq d(Fx,Fy)\leq Ld(x,y)+L\quad
\text{for all}\quad x,y\in X.\]
The map $F$ is called an  \emph{$L$-quasi-isometry}
if moreover
for every $y\in Y$ there is some $x\in X$ with
$d(Fx,y)\leq L$. The spaces $(X,d),(Y,d)$ are then called
\emph{quasi-isometric}.

Every finitely generated group $G$ admits a natural
family of metrics which are pairwise quasi-isometric.
Namely,
choose a finite symmetric set ${\cal G}$ of generators for
$G$. Then every element $g\in G$ can
be represented as a word in the alphabet ${\cal G}$. The minimal
length $\vert g\vert $ of such a word defines the
\emph{word norm} of $g$. This word norm induces
a metric on $G$ which is invariant under left
translation
by defining $d(g,h)=\vert g^{-1}h\vert$.
The word norm $\vert \,\vert^\prime$ defined by
a different set of generators is equivalent to $\vert\,\vert$
and hence the induced metrics $d,d^\prime$ are
quasi-isometric. In particular,
we can talk
about quasi-isometric finitely generated groups.
Note that a finitely generated group is
quasi-isometric to each of its finite index subgroups
and quasi-isometric to its image under a homomorphism
with finite kernel.
The main purpose of this note is to show.

\bigskip

\noindent
{\bf Theorem A:} {\it Let $\Gamma$ be a finitely
generated group which is quasi-isometric to ${\cal M}(S)$.
Then there is a finite index subgroup $\Gamma^\prime$
of $\Gamma$ and a
homomorphism $\rho:\Gamma^\prime\to {\cal M}(S)$
with finite kernel and finite index image.}

\bigskip

For surfaces with precisely one puncture (i.e. in the
case $m=1$), our theorem was earlier shown by
Mosher and Whyte (see \cite{M03b} for more details).
Kida \cite{K06} recently showed an analogous rigidity
result in the context of \emph{measure equivalence}.
Namely, call two countable groups $\Gamma,\Lambda$
\emph{measure equivalent} if $\Gamma,\Lambda$
admit commuting measure preserving actions on a
standard Borel
space $X$ equipped with a Radon measure $\mu$ and
with finite measure fundamental domains.
Motivated by deep results of Zimmer and Furman,
Kida showed that for a countable group $\Gamma$ which
is measure equivalent to the mapping class group
the conclusion of our theorem holds true.
Note however that measure equivalence for finitely
generated groups is neither implied
by nor implies quasi-isometry.

A choice of a word norm for the mapping class group and of
a non-principal
ultrafilter on $\mathbb{N}$ determines
an \emph{asymptotic cone} of ${\cal M}(S)$. The homological
dimension of this cone, i.e. the maximal number
$n\geq 0$ such that there are two open subsets
$V\subset U$ with $H_n(U,U-V)\not=0$,
is independent of the choices. We show the following
version of a
result of Behrstock and Minsky \cite{BM05}.

\bigskip

\noindent
{\bf Theorem B \cite{BM05}:} {\it The homological dimension of
an asymptotic cone of ${\cal M}(S)$
equals $3g-3+m$.}

\bigskip

The \emph{geometric rank} of a metric space $X$ is defined to be
the maximal number $k\geq 0$ such that there is a quasi-isometric
embedding $\mathbb{R}^k\to X$; it is not bigger than the
homological dimension of an asymptotic cone for ${\cal M}(S)$.
Farb, Lubotzky and Minsky \cite{FLM01} showed that the geometric
rank of ${\cal M}(S)$ is at least $3g-3+m$; thus as an
immediate corollary of our theorem we obtain \cite{BM05}.

\bigskip

\noindent
{\bf Corollary \cite{BM05}:} {\it The geometric rank of
${\cal M}(S)$ equals $3g-3+m$.}

\bigskip

The organization of this paper is as follows.
In Section 2 we summarize those of the properties
of the \emph{train track complex}
${\cal T\cal T}$ introduced in \cite{H06a}
which are needed for our purpose.
Section 3 discusses train tracks which hit efficiently.
Building on results
from \cite{H06b}, we analyze in Section 4 in more
detail the distance in the train track complex.
This is used in Section 5 to single out
a collection of infinite
subsets of ${\cal M}(S)$ whose asymptotic cones
are homeomorphic to euclidean cones of dimension
at most $3g-3+m$. In Section 6 we establish a
fairly precise description of the asymptotic cone of
${\cal M}(S)$ which leads to the proof of Theorem B.
Section 7 then contains the proof of Theorem A.

\section{The complex of train tracks}

In this section we summarize some results
and constructions from
\cite{PH92,H06a,H06b} which will be used throughout the paper
(compare also \cite{M03}).

Let $S$ be an
oriented surface of
genus $g\geq 0$ with $m\geq 0$ punctures and where $3g-3+m\geq 2$.
A \emph{train track} on $S$ is an embedded
1-complex $\tau\subset S$ whose edges
(called \emph{branches}) are smooth arcs with
well-defined tangent vectors at the endpoints. At any vertex
(called a \emph{switch}) the incident edges are mutually tangent.
Through each switch there is a path of class $C^1$
which is embedded
in $\tau$ and contains the switch in its interior. In
particular, the branches which are incident
on a fixed switch are divided into
``incoming'' and ``outgoing'' branches according to their inward
pointing tangent at the switch. Each closed curve component of
$\tau$ has a unique bivalent switch, and all other switches are at
least trivalent.
The complementary regions of the
train track have negative Euler characteristic, which means
that they are different from discs with $0,1$ or
$2$ cusps at the boundary and different from
annuli and once-punctured discs
with no cusps at the boundary.
We always identify train
tracks which are isotopic. A train track is called
\emph{maximal} if its complementary components
are all trigons or once punctured monogons.
A train track $\tau$ is called \emph{large} if its complementary
components are all topological discs and once
punctured topological discs.

A \emph{trainpath} on a train track $\tau$ is a $C^1$-immersion
$\rho:[m,n]\to \tau\subset S$ which maps each interval $[k,k+1]$
$(m\leq k\leq n-1)$ onto a branch of $\tau$. The integer $n-m$ is
then called the \emph{length} of $\rho$. We sometimes identify a
trainpath on $S$ with its image in $\tau$. Each complementary
region of $\tau$ is bounded by a finite number of trainpaths which
either are simple closed curves or terminate at the cusps of the
region.
A \emph{subtrack} of a train track $\tau$ is a subset $\sigma$ of
$\tau$ which itself is a train track. Thus every switch of
$\sigma$ is also a switch of $\tau$, and every branch of $\sigma$
is an embedded trainpath of $\tau$. We write $\sigma<\tau$ if
$\sigma$ is a subtrack of $\tau$.

A train track is called \emph{generic} if all switches are
at most trivalent.
The train track $\tau$ is called \emph{transversely recurrent} if
every branch $b$ of $\tau$ is intersected by an embedded simple
closed curve $c=c(b)\subset S$ which intersects $\tau$
transversely and is such that $S-\tau-c$ does not contain an
embedded \emph{bigon}, i.e. a disc with two corners at the
boundary.

A \emph{transverse measure} on a train track $\tau$ is a
nonnegative weight function $\mu$ on the branches of $\tau$
satisfying the \emph{switch condition}:
For every switch $s$ of $\tau$, the sum of the weights
over all incoming branches at $s$
is required to coincide with the sum of
the weights over all outgoing branches at $s$.
The train track is called
\emph{recurrent} if it admits a transverse measure which is
positive on every branch. We call such a transverse measure $\mu$
\emph{positive}, and we write $\mu>0$.
If $\mu$ is any transverse measure on a train track
$\tau$ then the subset of $\tau$ consisting of all
branches with positive $\mu$-mass is a recurrent
subtrack of $\tau$.
A train track $\tau$ is called \emph{birecurrent} if
$\tau$ is recurrent and transversely recurrent.

A \emph{geodesic lamination} for a complete
hyperbolic structure on $S$ of finite volume is
a \emph{compact} subset of $S$ which is foliated into simple
geodesics.
A geodesic lamination $\lambda$ is called \emph{minimal}
if each of its half-leaves is dense in $\lambda$. Thus a simple
closed geodesic is a minimal geodesic lamination. A minimal
geodesic lamination with more than one leaf has uncountably
many leaves and is called \emph{minimal arational}.
Every geodesic lamination $\lambda$ consists of a disjoint union of
finitely many minimal components and a finite number of isolated
leaves. Each of the isolated leaves of $\lambda$ either is an
isolated closed geodesic and hence a minimal component, or it
\emph{spirals} about one or two minimal components
\cite{CEG87,O96}.

A geodesic
lamination is \emph{finite} if it contains only finitely many
leaves, and this is the case if and only if each minimal component
is a closed geodesic. A geodesic lamination is \emph{maximal}
if its complementary regions are all ideal triangles
or once punctured discs with one cusp at the boundary.
The space of all geodesic laminations on $S$
equipped with the \emph{Hausdorff topology} is
a compact metrizable space. A geodesic lamination $\lambda$
is called \emph{complete} if $\lambda$ is maximal and
can be approximated in the Hausdorff topology by
simple closed geodesics. The space ${\cal C\cal L}$
of all complete geodesic laminations equipped with
the Hausdorff topology is compact.
Every geodesic lamination $\lambda$
which is a disjoint union of finitely many minimal components
is a \emph{sublamination} of
a complete geodesic lamination, i.e. there
is a complete geodesic lamination which contains
$\lambda$ as a closed subset \cite{H06a}.

A train track or a geodesic lamination $\sigma$ is
\emph{carried} by a transversely recurrent train track $\tau$ if
there is a map $F:S\to S$ of class $C^1$ which is isotopic to the
identity and maps $\sigma$ into $\tau$ in
such a way that the restriction
of the differential of $F$ to the tangent space of $\sigma$
vanishes nowhere; note that this makes sense since a train track
has a tangent line everywhere. We call the restriction of $F$ to
$\sigma$ a \emph{carrying map} for $\sigma$. Write $\sigma\prec
\tau$ if the train track or the geodesic
lamination $\sigma$ is carried by the train track
$\tau$. If $\sigma$ is a train track which is carried
by a train track $\tau$ then
every geodesic lamination $\lambda$ which is carried
by $\sigma$ is also carried by $\tau$.

A train track
$\tau$ is called \emph{complete} if it is generic and
transversely recurrent and if
it carries
a complete geodesic lamination. A complete train track
is maximal and birecurrent.
The space
of all complete geodesic laminations which are
carried by a fixed complete
train track $\tau$ is open and closed in ${\cal C\cal L}$.
In particular, the space ${\cal C\cal L}$
is totally disconnected \cite{H06a}.

A half-branch $\hat b$ in a generic train track $\tau$ incident on
a switch $v$ of $\tau$ is called
\emph{large} if every trainpath containing $v$ in its interior
passes through $\hat b$. A half-branch which is not large
is called \emph{small}.
A branch
$b$ in a generic train track
$\tau$ is called
\emph{large} if each of its two half-branches is
large; in this case $b$ is necessarily incident on two distinct
switches, and it is large at both of them. A branch is called
\emph{small} if each of its two half-branches is small. A branch
is called \emph{mixed} if one of its half-branches is large and
the other half-branch is small (for all this, see \cite{PH92} p.118).

There are two simple ways to modify a train track $\tau$
to another train track. First, we can \emph{shift}
$\tau$ along a mixed branch to a train track $\tau^\prime$ as shown
in Figure A below. If $\tau$ is complete then the same is true for
$\tau^\prime$. Moreover, a train track or a
lamination is carried by $\tau$ if and
only if it is carried by $\tau^\prime$ (see \cite{PH92} p.119).
In particular, the shift $\tau^\prime$ of $\tau$ is
carried by $\tau$. Note that there is a natural
bijection of the set of branches of $\tau$ onto
the set of branches of $\tau^\prime$.

\begin{figure}[ht]
\includegraphics{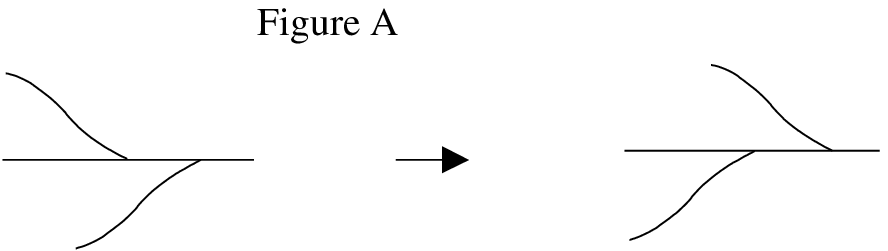}
\end{figure}

Second, if $e$ is a large branch of $\tau$ then we can perform a
right or left \emph{split} of $\tau$ at $e$ as shown in Figure B.
A right split at $e$ is uniquely
determined by the orientation of $S$ and does not
depend on the orientation of $e$.
Using the labels in the figure, in the case of a right
split we call the branches $a$ and $c$ \emph{winners} of the
split, and the branches $b,d$ are \emph{losers} of the split. If
we perform a left split, then the branches $b,d$ are winners of
the split, and the branches $a,c$ are losers of the split.
The split $\tau^\prime$ of a train track $\tau$ is carried
by $\tau$, and there is a natural choice of a carrying map which
maps the switches of $\tau^\prime$ to the switches of $\tau$. The
image of a branch of $\tau^\prime$ is then a trainpath on $\tau$
whose length either equals one or two. There is
a natural bijection of the set of branches
of $\tau$ onto the set of branches of $\tau^\prime$ which
maps the branch $e$ to the diagonal $e^\prime$ of the split.
The split of a maximal
transversely recurrent generic train track is maximal,
transversely recurrent and generic. If $\tau$ is complete and if
$\lambda\in {\cal C\cal L}$ is carried by $\tau$, then there is a
unique choice of a right or left split of $\tau$ at $e$ with the
property that the split track $\tau^\prime$ carries $\lambda$.
We call such a split a \emph{$\lambda$-split}.
The train track $\tau^\prime$ is recurrent and hence
complete. In particular, a complete train track $\tau$ can always
be split at any large branch $e$ to a complete train track
$\tau^\prime$; however there may be a choice of a right or left
split at $e$ such that the resulting
train track is not recurrent any
more (compare p.120 in \cite{PH92}). The reverse of a split is called a
\emph{collapse}. Define moreover a \emph{collision}
of a train track $\tau$ at a large branch $e$ to
be a (right or left) split of $\tau$ at $e$ followed
by the removal of the diagonal of the split. The train track
obtained from $\tau$ by a collision at $e$ is
carried by both train tracks obtained from $\tau$ by a split
at $e$. The number of its branches equals the number
of branches of $\tau$ minus one. If $\tau$ is complete
then the collision of $\tau$ at $e$ is recurrent if and
only if both train tracks obtained from $\tau$ by a split
at $e$ are complete (Lemma 2.1.3 of \cite{PH92}).

\begin{figure}[ht]
\includegraphics{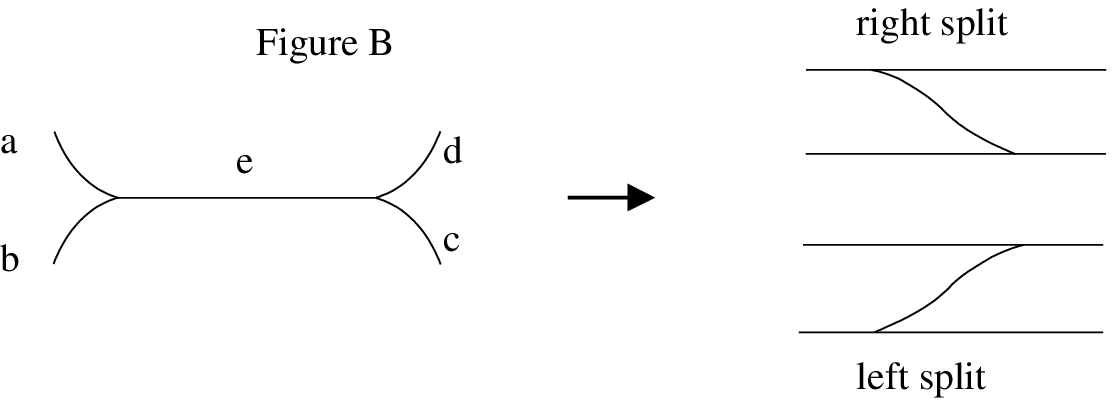}
\end{figure}

Denote by ${\cal T\cal T}$ the directed graph whose
vertices are the isotopy classes of complete train tracks on $S$
and whose edges are determined as follows. The train track
$\tau\in {\cal T\cal T}$
is connected to the train track $\tau^\prime$
by a directed edge if and only $\tau^\prime$ can be obtained
from $\tau$ by a single split.
The graph ${\cal T\cal T}$ is connected \cite{H06a}.
As a consequence, if we identify each edge in ${\cal T\cal T}$
with the unit interval $[0,1]$ then this provides
${\cal T\cal T}$ with the structure of a connected locally finite metric
graph. Thus ${\cal T\cal T}$ is a locally compact complete geodesic
metric space. In the sequel we always assume that
${\cal T\cal T}$
is equipped with this metric without further comment.
The mapping class group ${\cal M}(S)$ of $S$ acts properly
and cocompactly on ${\cal T\cal T}$ as a group of
isometries. In particular, ${\cal T\cal T}$
is ${\cal M}(S)$-equivariantly quasi-isometric to
${\cal M}(S)$ equipped with any word metric \cite{H06a}.

In the sequel we write $\tau\in {\cal V}({\cal T\cal T})$ if
$\tau$ is a \emph{vertex} of the graph
${\cal T\cal T}$, i.e. if $\tau$ is a complete train track on $S$.
Define a \emph{splitting sequence} in ${\cal T\cal T}$
to be a sequence
$\{\alpha(i)\}_{0\leq i\leq m}\subset{\cal V}({\cal T\cal T})$ with
the property that for every $i\geq 0$ the train track
$\alpha(i+1)$ can be obtained from $\alpha(i)$ by
a single split. We view such a splitting sequence
as a simplicial path in the graph ${\cal T\cal T}$
which maps the interval $[i,i+1]$ onto
the edge in ${\cal T\cal T}$ connecting
$\alpha(i)$ to $\alpha(i+1)$.

Recall from the introduction the definition of an $L$-quasi-isometric
embedding of a metric space $(X,d)$
into a metric space $(Y,d)$.
A \emph{$c$-quasi-geodesic} in a metric space
$(X,d)$ is a $c$-quasi-isometric embedding of a closed connected
subset of $\mathbb{R}$ into $X$.
The following two results from
\cite{H06b} will be important in the sequel.

\begin{proposition}\label{quasigeodesic}
There is a number $c>0$
such that every splitting sequence in ${\cal T\cal T}$
is a $c$-quasi-geodesic.
\end{proposition}

\begin{proposition}\label{density}
There is a number $d>0$
with the following property. For arbitrary train tracks
$\tau,\sigma\in {\cal V}({\cal T\cal T})$ there is a train track
$\tau^\prime$ contained in the $d$-neighborhood of
$\tau$ which is splittable to a train track $\sigma^\prime$
contained in the $d$-neighborhood of $\sigma$.
\end{proposition}

\section{Train tracks hitting efficiently}

In this section we construct complete train tracks with some
specific properties which are used to obtain a geometric control
on the train track complex ${\cal T\cal T}$. First,
define a \emph{bigon track} on $S$ to be an embedded 1-complex on
$S$ which satisfies all the requirements of a train track except
that we allow the existence of complementary bigons. Such a bigon
track is called \emph{maximal} if all complementary components are
either bigons or trigons or once punctured monogons. Recurrence,
transverse recurrence, birecurrence and carrying for bigon tracks
are defined in the same way as they are defined for train tracks.
Any complete train track is a maximal birecurrent bigon track in
this sense. A \emph{tangential measure} for a maximal bigon track
$\zeta$ assigns to each branch $b$ of $\zeta$ a nonnegative weight
$\nu(b)\in [0,\infty)$ with the following properties. Each side of
a complementary component of $\zeta$ can be parametrized as a
trainpath $\rho$ on $\zeta$. Denote by $\nu(\rho)$ the sum of the
weights of the branches contained in $\rho$ counted with
multiplicities. If $\rho_1,\rho_2$ are the two distinct sides of a
complementary bigon then we require that
$\nu(\rho_1)=\nu(\rho_2)$, and if $\rho_1,\rho_2,\rho_3$ are the
three distinct sides of a complementary trigon then we require
that $\nu(\rho_i)\leq \nu(\rho_{i+1})+\nu(\rho_{i+2})$ where
indices are taken modulo 3. A bigon track is transversely
recurrent if and only if it admits a tangential measure which is
positive on every branch \cite{PH92}.

A bigon track is called \emph{generic} if all switches are
at most trivalent. A bigon track $\tau$ which is not generic can
be \emph{combed} to a generic bigon track
by successively modifying $\tau$ as shown in Figure C.
By Proposition 1.4.1 of \cite{PH92} (whose proof is also valid
for bigon tracks),
the combing of a recurrent bigon track is recurrent.
However, the combing of a transversely recurrent bigon track
need not be transversely recurrent (see the discussion on
p.41 of \cite{PH92}).

\begin{figure}[ht]
\includegraphics{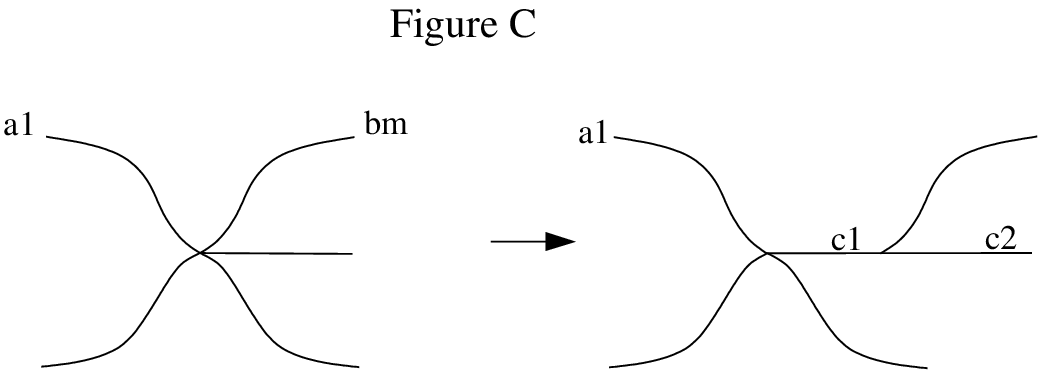}
\end{figure}

The next Lemma gives a criterion for a non-generic
maximal transversely recurrent bigon track to be
combable to a generic maximal transversely recurrent bigon track.
For its formulation, we say that a positive
tangential measure
$\nu$ on a maximal bigon track
$\sigma$ satisfies the \emph{strict triangle inequality
for complementary trigons} if
for every
complementary trigon of $\sigma$ with sides
$e_1,e_2,e_3$ we have $ \nu(e_i)<\nu(e_{i+1})+
\nu(e_{i+2})$. By
Theorem 1.4.3 of \cite{PH92}, a \emph{generic} maximal
train track
is transversely recurrent if and only if it admits a
positive tangential measure satisfying the strict triangle
inequality for complementary trigons.
We have.

\begin{lemma}\label{combing} Let $\zeta$ be a maximal bigon track
which admits a positive
tangential measure satisfying the strict triangle
inequality for complementary trigons.
Then $\zeta$ can be combed to a
generic transversely recurrent bigon track.
\end{lemma}

\begin{proof} Let $\sigma$ be an arbitrary maximal
bigon track.
Then $\sigma$
does not have any bivalent switches.
For a switch $s$ of $\sigma$ denote the valence of $s$ by
$V(s)$ and define the excessive
total valence ${\cal V}(\sigma)$ of
$\sigma$ to be
$\sum_s (V(s)-3)$ where
the sum is taken over all switches
$s$ of $\sigma$; then ${\cal V}(\sigma)=0$
if and only if $\sigma$ is generic. By induction it is
enough to show that a maximal non-generic
bigon track $\sigma$ which admits
a positive tangential measure $\nu$ satisfying the
strict triangle inequality for complementary trigons can be
combed to a bigon track $\sigma^\prime$ which admits
a positive tangential measure $\nu^\prime$ satisfying the
strict triangle inequality for complementary trigons and such that
${\cal V}(\sigma^\prime)<{\cal V}(\sigma)$.

For this let $\sigma$ be such a non-generic maximal bigon track
with tangential measure $\nu$ satisfying the strict triangle
inequality for complementary trigons and let $s$ be a switch of
$\sigma$ of valence at least $4$. Assume that $s$ has $\ell$
incoming and $m$ outgoing branches where $1\leq \ell\leq m$ and
$\ell +m\geq 4$. We number the incoming branches in
counter-clockwise order $a_1,\dots,a_\ell$ (for the given
orientation of $S$) and do the same for the outgoing branches
$b_1,\dots,b_m$. Then the branches $b_m$ and $a_1$ are contained
in the same side of a complementary component of $\sigma$, and the
branches $b_{m-1}, b_m$ are contained in adjacent (not necessarily
distinct) sides $e_1,e_2$ of a complementary component $T$ of
$\sigma$. Assume first that $T$ is a complementary trigon. Denote
by $e_3$ the third side of $T$; by assumption, the total weight
$\nu(e_3)$ is strictly smaller than $\nu(e_1)+\nu(e_2)$ and
therefore there is a number $q\in (0,\min\{\nu(b_i)\mid 1\leq
i\leq m\})$ such that $\nu(e_3)<\nu(e_1)+\nu(e_2)-2q$. Move the
endpoint of the branch $b_m$ to a point in the interior of
$b_{m-1}$ as shown in Figure C; we obtain a bigon track
$\sigma^\prime$ with ${\cal V}(\sigma^\prime)< {\cal V}(\sigma)$.

The branch $b_{m-1}$ decomposes in $\sigma^\prime$ into the union
of two branches $c_1,c_2$ where $c_1$ is incident on $s$ and on
an endpoint of the image $b_m^\prime$ of $b_m$ under our move.
Assign the weight $q$ to the branch $c_1$, the weight $\nu(b_m)-q$
to the branch $b_m^\prime$ and the weight $\nu(b_{m-1})-q$ to the
branch $c_2$. The remaining branches of $\sigma^\prime$ inherit
their weight from the tangential measure $\nu$ on $\sigma$. This
defines a positive weight function on the branches of
$\sigma^\prime$ which by the choice of $q$
is tangential measure satisfies the strict triangle inequality for
complementary trigons.

Similarly, if the complementary component $T$ containing
$b_m$ and $b_{m-1}$ in its boundary is
a bigon or a once punctured monogon, then we can shift $b_m$ along
$b_{m-1}$ as before and modify our tangential measure to a
positive tangential measure on the combed track with the desired
properties. \end{proof}

Following \cite{PH92} we say that a train track $\tau$ on our
surface $S$ \emph{hits efficiently} a train track or a geodesic
lamination $\sigma$ if $\tau$ can be isotoped to a train track
$\tau^\prime$ which intersects $\sigma$ transversely in such a way
that $S-\tau^\prime-\sigma$ does not contain any embedded bigon.
As in \cite{H06a} we define a \emph{splitting and shifting
sequence} to be a sequence
$\{\tau_i\}\subset {\cal V}({\cal T\cal T})$
such that for every $i$ the train track $\tau_{i+1}$ can
be obtained from $\tau_i$ by a sequence of shifts and
a single split. Denote by $d$ the distance on ${\cal T\cal T}$.
We have.

\begin{proposition}\label{backwards}
There is a number $q>0$ and for every
$\tau\in {\cal V}({\cal T\cal T})$ and every complete geodesic lamination
$\lambda$ which hits $\tau$ efficiently there is a complete train
track $\tau^{*}\in {\cal V}({\cal T\cal T})$
with the following properties.
\begin{enumerate}
\item $d(\tau,\tau^*)\leq q$.
\item $\tau^{*}$ carries $\lambda$.
\item
Let $\sigma\in {\cal V}({\cal T\cal T})$
be a train track which hits $\tau$ efficiently and
carries $\lambda$; then $\tau^{*}$
carries a train track $\sigma^\prime$ which carries
$\lambda$ and can be obtained from $\sigma$ by
a splitting and shifting sequence
of length at most $q$.
\end{enumerate}
\end{proposition}

\begin{proof} By Lemma 3.4.4 and
Proposition 3.4.5 of \cite{PH92}, for every
complete train track $\tau$ there is a maximal birecurrent
\emph{dual bigon track} $\tau_b^*$ with the following property. A
geodesic lamination or a train track $\sigma$ hits $\tau$
efficiently if and only if $\sigma$ is carried by $\tau_b^*$. We
construct the train track $\tau^{*}$ with the properties stated
in the lemma from
this dual bigon track and a complete geodesic
lamination $\lambda\in {\cal
C\cal L}$ which hits $\tau$ efficiently and hence is carried by
$\tau_b^*$.

For this we recall from p.194 of \cite{PH92} the precise
construction of the dual bigon track $\tau_b^*$ of a complete
train track $\tau$. Namely, for each branch $b$ of $\tau$ choose a
short arc $b^*$ meeting $\tau$ transversely in a single point in
the interior of $b$ and such that all these arcs are pairwise
disjoint. Let $T\subset S-\tau$ be a complementary trigon of
$\tau$ and let $E$ be a side of $T$ which is composed of the
branches $b_1,\dots,b_\ell$. Choose a point $p\in T$ and extend
all the arcs $b_1^*,\dots,b_\ell^*$ within $T$ in such a way that
they end at $p$, with the same inward pointing tangents at $p$. In
the case $\ell\geq 2$ we then add an arc which connects $p$ within
$T$ to a point $p^\prime\in T$ and whose inward pointing tangent at $p$
equals the outward pointing tangent at $p$ of the arcs
$b_1^*,\dots,b_\ell^*$. We do this in such a way that the
different configurations from the different sides of $T$ are
disjoint. If $q^\prime\in T$ is the point in $T$ arising in this
way from a second side, then we connect $p^\prime$ (or $p$ if
$\ell=1$) and $q^\prime$ by a smooth arc whose outward pointing
tangent at $p^\prime,q^\prime$ coincides with the inward pointing
tangents of the arcs constructed before which end at
$p^\prime,q^\prime$. In a similar way we construct the
intersection of $\tau_b^*$ with a complementary once punctured
monogon of $\tau$. Note that the resulting graph $\tau_b^*$ is in
general not generic, but its only vertices which are not trivalent
arise from the sides of the complementary components of $\tau$.
Figure D shows the intersection of the dual bigon track $\tau_b^*$
with a neighborhood in $S$ of a complementary trigon of $\tau$ and
with a neighborhood in $S$ of a complementary once punctured
monogon.

\begin{figure}[ht]
\includegraphics{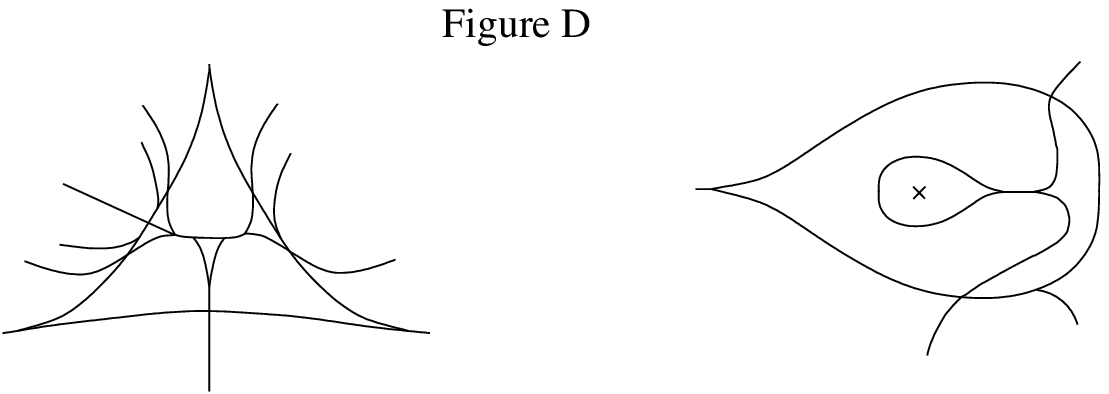}
\end{figure}

Following \cite{PH92}, $\tau_b^*$ is a maximal birecurrent bigon
track, and the number of its branches is bounded from above by a
constant only depending on the topological type of $S$.
Each complementary trigon of
$\tau$ contains exactly one complementary trigon of $\tau_b^*$ in
its interior, and  these are the only complementary trigons. Each
complementary once punctured monogon of $\tau$ contains exactly
one complementary once punctured monogon of $\tau_b^*$ in its
interior. All other complementary components of $\tau_b^*$ are
bigons. The number of these bigons is uniformly bounded.

Now let $\mu$ be a positive integral transverse measure on $\tau$
with the additional property that the $\mu$-weight of every branch
of $\tau$ is at least $4$.
This weight then defines a
\emph{simple multicurve} $c$ carried by $\tau$
in such a way that $\mu$
is just the counting measure for $c$ (see \cite{PH92}).
Here a simple multicurve consists of a disjoint
union of essential simple closed curves which
can be realized disjointly; we allow that some of the curves
are freely homotopic. For every
side $\rho$ of a complementary component of $\tau$ there are at
least $4$ connected subarcs of $c$ which are mapped by the natural
carrying map $c\to \tau$ \emph{onto} $\rho$. Namely, the number of
such arcs is just the minimum of the $\mu$-weights of a branch
contained in $\rho$.

Assign to a branch $b^*$ of $\tau_b^*$ which is dual to the branch
$b$ of $\tau$ the weight $\nu(b^*)=\mu(b)$, and to a branch of
$\tau_b^*$ which is contained in the interior of a complementary
region of $\tau$ assign the weight $0$. The resulting weight
function $\nu$ is a tangential measure for $\tau_b^*$, but it is
not positive (this relation between transverse measures
on $\tau$ and tangential measures on $\tau_b^*$ is
discussed in detail in Section 3.4 of \cite{PH92}).
However by construction, every branch of vanishing
$\nu$-mass is contained in the interior of a complementary trigon
or once punctured monogon of $\tau$, and positive mass can be
pushed onto these branches by ``sneaking up'' as described on p.39
and p.200 of \cite{PH92}. Namely, the closed multicurve $c$ defined by
the positive integral transverse measure $\mu$ on $\tau$ hits the
bigon track
$\tau_b^*$ efficiently. For every branch $b$ of $\tau$
the weight $\nu^*(b)$ equals the number of
intersections between $b^*$ and $c$.
For each side of a complementary component
$T$ of $\tau$ there are at least $4$ arcs from $c$ which are
mapped by the carrying map onto this side. If the side consists of
more than one branch then we pull two of these arcs into $T$ as
shown in Figure E. If the side consists of a single branch then we
pull a single arc into $T$ in the same way.
\begin{figure}[ht]
\includegraphics{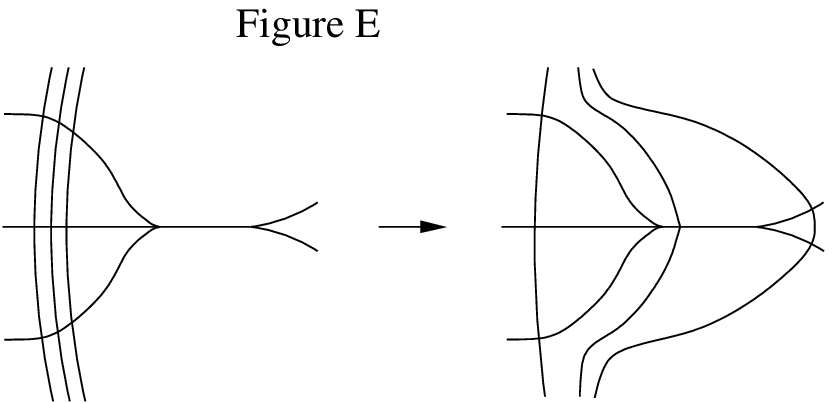}
\end{figure}

For a branch
$e$ of $\tau_b^*$ define $\mu^*(e)$ to be the
number of intersections between $e$ and the deformed
multicurve. The resulting weight function
$\mu^*$ is
a positive integral tangential measure for $\tau_b^*$.
Note that the weight of each
side of a complementary trigon in $\tau_b^*$ is exactly 2 by
construction, and the weight of a side of a once punctured
monogon is 2 as well. In particular, the tangential measure
$\mu^*$ satisfies the strict triangle inequality for complementary
trigons: If $T$ is any complementary trigon with sides
$e_1,e_2,e_3$ then $\mu^*(e_i)<\mu^*(e_{i+1})+\mu^*(e_{i+2})$
(compare the proof of Lemma 3.1).

We now modify our dual bigon track $\tau_b^*$ in a uniformly
bounded number of steps to a complete
train track $\tau^*$ as required in the lemma with
a (non-deterministic) algorithm as follows.

The set of input data
for our algorithm is the set ${\cal B}$ of quadruples
$(\eta,\lambda,\nu,B)$
which consist of a
maximal birecurrent bigon track
$\eta$, a complete geodesic lamination $\lambda$ carried by
$\eta$, a positive tangential measure $\nu$ on $\eta$
which satisfies the strict triangle inequality
for complementary trigons and
a complementary bigon $B$ of $\eta$. If $\eta$ does not
have any complementary bigons, i.e. if
$\eta$ is a train track, then we put $B=\emptyset$.
The algorithm
modifies the quadruple $(\eta,\lambda,\nu,B)\in {\cal B}$ to a
quadruple $(\eta^\prime,\lambda,\nu^\prime,B^\prime)\in {\cal B}$
with $B^\prime=\emptyset$
as follows.

{\sl Step 1:}

Let $(\eta,\lambda,\nu,B)\in {\cal B}$ be an input
quadruple. If $\eta$ does not contain a bigon, i.e.
if $B=\emptyset$, then the algorithm stops.
Otherwise proceed to Step 2.

{\sl Step 2:}

Let $(\eta,\lambda,\nu,B)\in {\cal B}$ be an input
quadruple with
$B\not=\emptyset$. Let $E,F$ be the
sides of the complementary bigon $B$ in the maximal
birecurrent bigon track $\eta$.
Check whether the boundary $\partial B$ of $B$ is embedded in $\eta$.
If this is not the case then go to Step 3.
Otherwise we construct
from $\eta$ and $\nu$ a maximal birecurrent bigon track $\tilde \eta$
which carries $\eta$ and hence $\lambda$ as
follows. Assume that the side $E$ of the bigon $B$ consists of the
ordered sequence of branches $e_1,\dots,e_\ell$ and that the
second side $F$ of $B$ consists of the branches $f_1,\dots,f_k$.
Assume also that the branches $e_1,f_1$ begin at a common cusp of
the bigon $B$. We collapse the bigon $B$ to a single arc in $S$
with a map $\Psi$ which identifies $E$ and $F$ as follows. If for
some $p\geq 1$, $q\geq 1$  we have $\sum_{j=1}^{q-1}
\nu(f_j)< \sum_{i=1}^p \nu(e_i)<\sum_{j=1}^q
\nu(f_j)$ then $\Psi$ maps the subarc $e_1\cup\dots\cup e_p$ of
$E$ homeomorphically onto a subarc of $F$ which contains
$f_1,\dots,f_{q-1}$ and has its endpoint in the interior of the
edge $f_q$. If $\sum_{i=1}^p\nu(e_i)=\sum_{j=1}^q
\nu(f_j)$ then we map $e_1\cup \dots \cup e_{p}$ onto
$f_1\cup\dots \cup f_q$, i.e. an endpoint of $e_p$ is mapped to an
endpoint of $f_q$. The resulting bigon track $\tilde \eta$ carries
$\eta$ and it is maximal. The natural carrying map $\Phi:\eta \to
\tilde \eta$ maps each complementary trigon of $\eta$ to a
complementary trigon of $\tilde \eta$. By construction, the
positive tangential measure $\nu$ on $\eta$ induces a positive
weight function $\tilde \nu$ on the branches of $\tilde \eta$.
Note that the total weight of $\tilde \nu$ is
strictly smaller than the total weight of $\nu$ and that the
$\nu$-weight of a side $\rho$ of a complementary component $T\not=
B$ in $\eta$ coincides with the $\tilde \nu$-weight of the side
$\Phi(\rho)$ of the complementary component
$\Phi(T)$ in $\tilde \eta$.
In particular, the weight function $\tilde \nu$
is a tangential measure $\tilde \nu$ on $\tilde \eta$
which satisfies the strict triangle inequality for
complementary trigons.
The number of complementary bigons in $\tilde \eta$ is strictly
smaller than the number of complementary bigons in $\eta$.
Namely, there is a one-to-one correspondence
between the complementary bigons of $\tilde \eta$
and the complementary bigons of $\eta$ distinct from $B$.
The image of the bigon $B$ under the map $\Phi$ is an
embedded arc in $\tilde \eta$. The number of branches of
$\tilde \eta$ does not exceed the number of branches of $\eta$.
Every complete geodesic lamination which is carried
by $\eta$ is also carried by $\tilde\eta$.
Choose an input quadruple of the form
$(\tilde \eta,\lambda,\tilde\nu,\tilde B)\in {\cal B}$
for a complementary bigon $\tilde B$ of $\tilde \eta$
(or $\tilde B=\emptyset$ if $\tilde \eta$ is a train track)
and continue with Step 1 above for this input quadruple.

{\sl Step 3:}

Let $(\eta,\lambda,\nu,B)\in {\cal B}$ be an input quadruple
such that $B\not=\emptyset$ and that the
boundary $\partial B$ of $B$ is
\emph{not}
embedded. Then the sides $E,F$ of $B$
are immersed arcs of class $C^1$ on $S$ which
intersect or have self-intersections.
Check whether the cusps of $B$ coincide. If this
is not the case, continue with Step 4.

Otherwise the two cusps of $B$ are a common
switch $s$ of $\eta$ which is necessarily at least
4-valent. By Lemma 3.1 and its proof, we can modify $\eta$ with a
sequence of combings near
$s$ to a maximal birecurrent bigon track $\tilde\eta$
in such a way that the two cusps of the
complementary bigon $\tilde B$ in $\tilde \eta$ corresponding
to $B$ under the combing are distinct and such that
the tangential measure $\nu$ on $\eta$
induces a tangential measure $\tilde \nu$
on $\tilde\eta$ which satisfies the strict
triangle inequality for complementary trigons.
Note that $\tilde\eta$ carries the complete
geodesic lamination $\lambda$.
Now continue with Step 2 above with the input
quadruple $(\tilde \eta,\lambda,\tilde \nu,\tilde B)\in {\cal B}$.

{\sl Step 4:}

Let $(\eta,\lambda,\nu,B)\in {\cal B}$ be our input quadruple
where $B$ is a bigon in $\eta$ with sides $E,F$ and
distinct cusps.
Check whether the boundary $\partial B$ of $B$ contains
any isolated self-intersection points.
Such a point is a switch $s$ contained in
the interior of at least two distinct embedded subarcs
$\rho_1,\rho_2$ of
$\partial B$ of class $C^1$ with the additional property
that $\rho_1-\{s\}\cap \rho_2-\{s\}=\emptyset$.
If $\partial B$ does not contain such an isolated
self-intersection point
then continue
with Step 5 below with the input
quadruple $(\eta,\lambda,\nu,B)\in {\cal B}$.

Otherwise any such isolated self-intersection
point $s$ is a switch of $\eta$ which is at
least 4-valent. Thus we can
modify $\eta$ with a sequence combing
to a complete birecurrent bigon track
$\tilde \eta$ with the property that
all self-intersection points of the boundary
of the bigon $\tilde B$ in $\tilde \eta$
corresponding to $B$
are non-isolated, i.e. they
are contained in a self-intersection branch, and that
the tangential measure $\nu$ on $\eta$
induces a tangential
measure $\tilde \mu$ on $\tilde \eta$ satisfying
the strict triangle inequality for
complementary trigons.
Continue with Step 5 with the input quadruple
$(\tilde \eta,\lambda,\tilde \nu,\tilde B)\in {\cal B}$.

{\sl Step 5:}

Let $(\eta,\lambda,\nu,B)\in {\cal B}$ be an input quadruple
where $B$ is a complementary bigon for $\eta$
whose boundary $\partial B$ does not contain any
isolated self-intersection points, whose
cusps are distinct and such that
the self-intersection of $\partial B$ is not empty.

Check whether there is a
branch $e$ of $\eta$ contained in the self-intersection
of $\partial B$ and
which is not incident on any one of the two
cusps $s_1,s_2$ of $B$. Since the
interior of the bigon $B$ is an embedded topological disc in $S$,
such a branch $e$ is necessarily a large branch. Now $\eta$
carries $\lambda$ and therefore
there is a bigon track $\tilde \eta$ which is the
image of $\eta$ under a split at $e$ and which carries $\lambda$.
To each complementary region of
$\tilde \eta$ naturally corresponds a complementary
region of $\eta$ of the same
topological type (compare the discussion in Section 3 of
\cite{H06b} and in Section 4 of this paper).
In particular, the number of complementary bigons in
$\tilde\eta$ and $\eta$ coincide, and the bigon
$B$ in $\eta$ corresponds to a bigon $\tilde B$ in
$\tilde \eta$.
The bigon track $\tilde \eta$
is recurrent, and it admits
a positive tangential measure $\tilde \nu$ induced
from the measure $\nu$ on $\eta$
which satisfies
the strict triangle inequality for complementary trigons.
The number of branches contained in the self-intersection locus
of the boundary $\partial \tilde B$ of
the bigon $\tilde B$
is strictly smaller than the number of branches in
the self-intersection locus of $\partial B$.

After a number of splits of this kind which is bounded from
above by the number of branches of the bigon track $\eta$
we obtain from $\eta$ a
bigon track $\eta_1$ which is maximal and birecurrent.
There is a natural bijection from the
collection of complementary bigons of $\eta$
onto the collection of complementary bigons of $\eta_1$.
If $B_1$ is the bigon of $\eta_1$ corresponding to $B$
then
the self-intersection locus of the boundary
$\partial B_1$ of $B_1$ is a union of branches
which are incident on one of the two
cusps $s_1\not=s_2$ of $B_1$. As before, $\eta_1$ admits
a positive tangential measure $\nu_1$ which satisfies
the strict triangle inequality for complementary trigons
and is induced from $\nu$. Moreover, $\eta_1$ carries $\lambda$.

If the boundary $\partial B_1$ of $B_1$ is embedded
then we proceed
with Step 2 above for the input
quadruple $(\eta_1,\lambda,\nu_1,B_1)$.
Otherwise there is a self-intersection
branch $b$ of $\partial B_1$ which is incident on
a cusp $s_1$ of $B_1$. Note that the branch
$b$ can \emph{not} be large, so it is either small
or mixed.

For a small branch $b$, there are again two possibilities
which are shown in Figure F. A small branch $b$ as shown
on the left hand side of Figure F can
be collapsed to a large branch.
\begin{figure}[ht]
\includegraphics{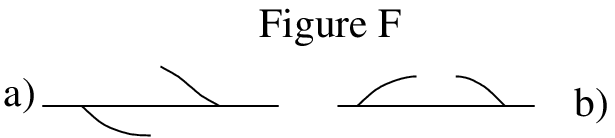}
\end{figure}
Using the fact that the boundary of an embedded
bigon on a bigon track admits a natural orientation
induced from the orientation of $S$,
the small branch $b$ is contained in the intersection
of the two distinct sides of $B_1$.
Since the tangential measure $\nu_1$ on $\eta_1$ is
positive by assumption, the construction in
Step 2 above can be used to collapse the bigon
$B_1$ in $\eta_1$ to a single simple closed curve.
As in Step 2 above, we obtain a maximal birecurrent bigon
track $\tilde \eta$ which carries $\lambda$ and
admits a positive tangential measure $\tilde \nu$ satisfying
the strict triangle inequality for complementary trigons.
The number of bigons of $\tilde \eta$ is strictly
smaller than the number of bigons of $\eta$.
Choose an arbitrary complementary bigon $\tilde B$
in $\tilde \eta$ or put $\tilde B=\emptyset$ if there
is no such bigon and
continue with Step 1 above and the input
quadruple $(\tilde \eta,\lambda,\tilde \nu,\tilde B)\in {\cal B}$.

A small branch $b$ as shown on the right hand side of Figure F
can not be collapsed. In this case the
branch $b$ coincides with a side $E_1$
of the bigon $B_1$, and the second side $F_1$ of
$B_1$ contains $b$ as a proper subarc.
Since the tangential measure
$\nu_1$ on $\eta_1$ is \emph{positive} by
assumption, this is impossible.

If the branch $b$ is mixed then
$b$ and the cusp $s_1$ of $B_1$ are
contained in the interior
of a side $E_1$ of $B_1$. The bigon track
$\eta_1$ can be modified with a single shift to
a maximal birecurrent bigon track $\eta_2$ such
that the switch $s_1$ is not contained any more
in the interior of a side of the bigon $B_2$
corresponding to $B_1$ in $\eta_1$.
The tangential measure $\nu_1$ on $\eta_1$
naturally induces a positive tangential
measure $\nu_2$ on $\eta_2$ which satisfies the
strict triangle inequality for complementary trigons.
We now proceed with Step 2 for the input
quadruple $(\eta_2,\lambda,\nu_2,B_2)$.
This completes the description of our algorithm.

We now apply our algorithm to the
bigon track $\tau_b^*$, the tangential measure
$\mu^*$ for $\tau_b^*$ constructed from a suitably
chosen transverse measure $\mu$ on $\tau$
and a complete geodesic lamination $\lambda$ which
hits $\tau$ efficiently and hence is carried
by $\tau_b^*$. Since the number of branches of
$\tau_b^*$ is bounded from above by a constant
only depending on the topological type of
$S$, there is a universal upper bound
$p>0$ for the number of modifications of $\tau_b^*$
needed in our above algorithm
to construct from these data a (possibly non-generic)
birecurrent train track $\chi$ which carries $\lambda$
and admits a positive tangential measure
satisfying the strict triangle inequality
for complementary trigons. By Lemma 3.1, this
train track can be combed to a maximal birecurrent
generic train track $\tau^*$ which satisfies
property 2) stated in the proposition.
The train track $\tau^*$ is not unique, and it depends
on $\lambda$ and $\mu$. However,
since our algorithm stops after
a uniformly bounded number of steps and each step involves
only a uniformly bounded number of choices,
the number of such train tracks which can be obtained
from $\tau_b^*$ by this procedure is uniformly bounded.
Moroever, our algorithm is equivariant with respect to
the action of the mapping class group ${\cal M}(S)$
and therefore by invariance under
the action of ${\cal M}(S)$, the distance between
$\tau$ and $\tau^*$ is uniformly bounded.
In other words, $\tau^*$ has property 1) stated in the
proposition as well.

To show property 3),
let $\sigma$ be a complete train track on $S$ which hits $\tau$
efficiently and carries a complete geodesic lamination $\lambda\in
{\cal C\cal L}$. Then $\sigma$ is carried by $\tau_b^*$
and hence it is carried by every bigon track which
can be obtained from $\tau_b^*$ by a sequence
of combings, shifts and collapses.
On the other hand, if $\eta_i$ $(0\leq i\leq k)$
are the successive bigon tracks obtained from
an application of our algorithm to
$\eta_0=\tau_b^*$ and if $\eta_i$ is obtained from
$\eta_{i-1}$ by a split at a large branch $e$,
then this split is a $\lambda$-split.
By Lemma 4.5
of \cite{H06a} and its proof (which is also
valid for bigon tracks)
there is a universal number $p>0$ such that
if $\sigma$ is carried by $\eta_{i-1}$ and carries
$\lambda$ then there is
a train track $\tilde \sigma$ which
carries $\lambda$, which is carried by $\eta_i$ and which
can be obtained from $\sigma$ by a splitting and shifting
sequence of length at most $p$.
Since the number of splits occuring during our
modification of $\tau_b^*$ to $\tau^*$
is uniformly bounded, this means that
$\tau^*$ also satisfied the third requirement
in our proposition.
This completes the proof of our proposition.
\end{proof}

{\bf Remark:} We call the train
track $\tau^*$ constructed in the
proof of Proposition 3.2 from a complete
train track $\tau$ and a complete geodesic lamination $\lambda$
which hits $\tau$ efficiently
a \emph{$\lambda$-collapse} of $\tau_b^*$. Note
that a $\lambda$-collapse is not unique, and that
in general it is neither carried by the dual bigon
track of $\tau$ nor carries this bigon track.
Thus in general a $\lambda$-collapse of
$\tau$ does not hit $\tau$ efficiently.
However, the number of different train
tracks which can be obtained from our construction is
bounded by a constant only depending on the topological
type of $S$.

\section{Flat strip projection}

In this section we use the results from Sections 3 and
from \cite{H06b} to obtain
a control on distances in the train track complex
${\cal T\cal T}$.

Note first that
a complementary
component $C$ of a
train track $\sigma$ on $S$ is bounded by a finite
number of arcs of class $C^1$, called \emph{sides}.
Each side either is a closed curve of class $C^1$ (i.e. the
side does not contain any cusp)
or an arc of class $C^1$ with endpoints at two not
necessarily distinct
cusps of the component. We call a side of $C$ which does not
contain cusps a \emph{smooth side} of $C$.
If $C$ is a complementary component of $\sigma$ whose
boundary contains precisely $k\geq 0$ cusps, then the \emph{Euler
characteristic} $\chi(C)$ is defined by $\chi(C)=\chi_0(C)-k/2$
where $\chi_0(C)$ is the usual Euler characteristic of
$C$ viewed as a topological surface with boundary. Note that the
sum of the Euler characteristics of the complementary components
of $\sigma$ is just the Euler characteristic of $S$.
If $T$ is a smooth side
of a complementary component $C$ of $\sigma$
then we mark a point on $T$ which is contained in the interior
of a branch of $\sigma$.
If $T$ is a common smooth side of two distinct complementary
components $C_1,C_2$ of $\sigma$ then we assume that
the marked points on $T$ defined by $C_1,C_2$ coincide.

A \emph{complete extension} of a train track
$\sigma$ is a complete train track $\tau$
containing $\sigma$ as a subtrack and whose switches
are distinct from the marked points in $\sigma$.
Such a complete extension $\tau$ intersects
each complementary component $C$ of $\sigma$ in an
embedded graph. The closure $\tau_C$
of $\tau\cap C$ in $S$
is a graph whose univalent
vertices are contained in the complement of the
cusps and the marked points of the boundary $\partial C$
of $C$. We call two such graphs
$\tau_C,\tau^\prime_C$
\emph{equivalent}
if there is an isotopy of $C$ of class $C^1$ which fixes
a neighborhood of the
cusps and the marked points in $\partial C$ and which
maps $\tau_C$ onto $\tau^\prime_C$.
The complete extensions
$\tau,\tau^\prime$ of $\sigma$ are called
\emph{$\sigma$-equivalent} if for each complementary component
$C$ of $\sigma$ the graphs $\tau_C$ and
$\tau^\prime_C$ are equivalent in this sense.

For two complete extensions $\tau,\tau^\prime$
of $\sigma$ and a complementary component $C$ of $\sigma$
define the \emph{$C$-intersection number}
$i_C(\tau,\tau^\prime)$ between $\tau$ and $\tau^\prime$
to be the minimal number
of intersection points contained in $C$ between
any two complete extensions $\eta,\eta^\prime$ of $\sigma$
which are $\sigma$-equivalent to $\tau,\tau^\prime$
and with the following additional properties.
\begin{enumerate}
\item[a)] A switch $v$ of $\eta$
is also a switch of $\eta^\prime$ if and only
if $v$ is a switch of $\sigma$.
\item[b)] For every intersection point $x\in \eta\cap \eta^\prime\cap C$
there is a branch $b$ of $\eta$ containing $x$ in its interior
and a branch $b^\prime$ of $\eta^\prime$ containing $x$ in
its interior. Moreover, the branches $b,b^\prime$ intersect
transversely at $x$.
\end{enumerate}
We define $i_\sigma(\tau,\tau^\prime)=\sum_Ci_C(\tau,\tau^\prime)$
to be the sum of the $C$-intersection numbers between $\tau$
and $\tau^\prime$ where
$C$ runs through the complementary components of $\sigma$.
For every number
$m>0$ there is a constant $q(m)>0$ not depending on $\sigma$
so that for every complete extension $\tau$ of $\sigma$
the number of $\sigma$-equivalence classes
of complete extensions $\tau^\prime$ of $\sigma$
with $i_\sigma(\tau,\tau^\prime)\leq m$ is bounded
from above by $q(m)$ (see \cite{H06b}).

Let again $\tau$ be a complete extension
of a train track $\sigma$.
To simplify our notation we do not
distinguish between $\sigma$ as a subset of $\tau$ (and hence
containing switches of $\tau$ which are bivalent in $\sigma$) and
$\sigma$ viewed as a subtrack of $\tau$, i.e. the graph from which
the bivalent switches not contained in simple closed curve
components have been removed. Define the \emph{$\sigma$-complexity
$\chi(\tau,\sigma)$ of $\tau$} to be the number of
branches of $\tau$ contained in $\sigma$. Note that
this complexity is not smaller than the
number of branches of $\sigma$, and it coincides with
this number if and only if $\sigma$ is itself a complete
train track. Similarly,
if $b$ is a branch in $\sigma$ then we define the
\emph{$b$-complexity} $\chi(\tau,b)$ of $\tau$
to be the number of branches of $\tau$ contained in $b$.
Note that $\chi(\tau,b)=1$ if and only if $b$ is a branch of $\tau$.

Let $b$ be any branch of $\sigma$ which is not a branch
of $\tau$. Then $b$ defines a trainpath $\rho:[0,m]\to \tau$
of length $m\geq 2$. For each $i<m$ the branch $\rho[i,i+1]$
of $\tau$ is a proper subset of $b$ and will be
called a \emph{proper subbranch} of $b$ in $\tau$.
A proper subbranch $h$ of $b$ in $\tau$ is incident on at least
one switch $v$ in $\tau$ which
is bivalent in $\sigma$, and there is a unique branch
$a\in \tau-\sigma$ which is incident on $v$. We call
$a$ a \emph{neighbor} of $b$ at $v$ or simply a
neighbor of $\sigma$ at $v$. A neighbor of
$\sigma$ at an interior point $v$ of a branch of $\sigma$
is small at $v$. A \emph{proper subbranch of $\sigma$}
is a proper subbranch of some branch of $\sigma$.
If $e$ is a large proper subbranch of $\sigma$
then $\tau$ can be split at $e$ to a train track
which contains $\sigma$ as a subtrack. We call such a
(not necessarily unique) split a \emph{$\sigma$-split}
of $\tau$ at $e$.

If $\sigma<\tau$ is a recurrent subtrack of a complete
train track $\tau$ then
there is a \emph{simple geodesic multi-curve} $\nu$ on $S$, i.e. a
collection of pairwise disjoint simple closed geodesics on $S$,
which is carried by $\sigma$ and such that
a carrying map $\nu\to\sigma$ is surjective
\cite{PH92}. We call such a multi-curve \emph{filling} for
$\sigma$. More generally, we call a geodesic lamination
$\nu$ carried by $\sigma$ \emph{filling} if the restriction to
the minimal components of $\nu$ of a carrying map $\nu\to \sigma$
is surjective. It follows from the results of
Section 2 in \cite{H06a} that for
every $\sigma$-filling lamination $\nu$
which is a disjoint union of minimal components
there is a complete geodesic
lamination $\lambda$ which is carried by $\tau$ and contains $\nu$
as a sublamination. We call $\lambda$ a \emph{complete
$\tau$-extension} of $\nu$. Whenever the existence
of such a lamination $\lambda$ is needed in the sequel, this
lamination is given explicitly so we refrain
from a more detailed discussion.
We call the train track $\tau$
\emph{tight} at a large branch $e$ of $\sigma$ if
$e$ is a large branch in $\tau$, i.e. if $\chi(\tau,e)=1$.
In Lemma 4.3 of \cite{H06b} we showed.

\begin{lemma}\label{tight}
There is a number $p>0$ with the
following properties. Let $\sigma$ be a
recurrent train track, let $e$ be
a large branch of $\sigma$ and let
$\nu$ be a $\sigma$-filling geodesic
lamination. Then there is
an algorithm which associates to every complete extension $\tau$
of $\sigma$ and every
complete $\tau$-extension $\lambda$ of $\nu$ a
complete train track $\tau^\prime$ with the following properties.
\begin{enumerate}
\item $\tau^\prime$ can be obtained from $\tau$ by a sequence of
$\lambda$-splits of length at most $p$ at proper subbranches of $e$,
and it contains $\sigma$ as a subtrack.
\item $\tau^\prime$ is tight at $e$.
\item If no marked point of $\sigma$
is contained in the branch $e$ then we have
$i_\sigma(\tau^\prime,\eta)\leq i_\sigma(\tau,\eta)+
p(\chi(\tau,\sigma)-\chi(\tau^\prime,\sigma))$ for
every complete extension $\eta$ of $\sigma$.
\item If there is a marked point of $\sigma$ contained in $e$
then we have $i_\sigma(\tau^\prime,\eta)\leq
i_\sigma(\tau,\eta)+p$ for every complete extension $\eta$ of
$\sigma$.
\end{enumerate}
\end{lemma}

For train tracks $\sigma<\tau$ as in Lemma 4.1 we call
$\tau^\prime$ the \emph{$\sigma$-modification} of $\tau$ at
the large branch $e$ with respect to the complete
$\tau$-extension $\lambda$ of $\nu$.

Let $\sigma$ be a recurrent train track on $S$
without closed curve components. Let
$\tau\in {\cal V}({\cal T\cal T})$ be a
complete extension of $\sigma$ and
let $\{\sigma(i)\}_{0\leq i\leq \ell}$ be a splitting
sequence issuing from $\sigma(0)=\sigma$.
We call this splitting sequence
\emph{recurrent} if each of the train tracks
in the sequence is recurrent.
Call a splitting
sequence $\{\eta(j)\}_{0\leq j\leq m}
\subset {\cal V}({\cal T\cal T})$
beginning at $\eta(0)=\tau$
\emph{induced} by the sequence $\{\sigma(i)\}_{0\leq
i\leq\ell}$ if there is an injective strictly increasing
map $q:\{0,\dots,\ell\}\to \{0,\dots,
m\}$ with the following properties.
\begin{enumerate}
\item[{a)}]
$q(0)=0,q(\ell)\leq m$ and for $q(i)\leq
j<q(i+1)$ the train track $\eta(j)$ contains a subtrack isotopic
to $\sigma(i)$.
\item[{b)}] Let $i<\ell$ and
assume that the split $\sigma(i)\to \sigma(i+1)$
is a right (or left) split at a large branch $e_i$. Then
there is a complete geodesic lamination $\lambda\in {\cal C\cal L}$
which is a complete $\eta(q(i+1)-1)$-extension of
a filling lamination for $\sigma(i)$ and
such that $\eta(q(i+1)-1)$ is the $\sigma(i)$-modification of
$\eta(q(i))$ at $e_i$ with respect to $\lambda$.
The train track $\eta(q(i+1)-1)$ is tight at $e_i$ and the
split $\eta(q(i+1)-1)\to \eta(q(i+1))$
is a right (or left) split at $e_i$.
\item[{c)}] For $q(\ell)\leq j<m$ the train track
$\eta(j+1)$ is obtained from $\eta(j)$ by a split
at a large proper subbranch of $\sigma(\ell)$.
\end{enumerate}

In the next lemma we compare
distances
between
complete train tracks which are obtained from splitting sequences
induced by a splitting sequence of a common subtrack.

\begin{lemma}\label{inducing}
There is a number $R_0>0$ with the following
property. Let $\sigma(0)$ be a train track on $S$
without closed curve components and let
$\{\sigma(i)\}_{0\leq i\leq \ell}$ be a finite
recurrent splitting sequence
issuing from $\sigma(0)$. Let $\tau(0),\eta(0)\in
{\cal V}({\cal T\cal T})$
be two complete extensions of $\sigma(0)$ and
let $\{\tau(j)\}_{0\leq j\leq m}\subset{\cal V}({\cal T\cal T}),
\{\eta(p)\}_{0\leq p\leq n}\subset {\cal V}({\cal T\cal T})$
be splitting sequences beginning at $\tau(0),\eta(0)$ which are
induced by $\{\sigma(i)\}_{0\leq i\leq \ell}$; then $d(\tau(m),
\eta(n))\leq d(\tau(0),\eta(0))+ R_0.$
\end{lemma}

\begin{proof}
Let $\sigma$ by a recurrent train track on $S$
without closed
curve components. Then
$\sigma$ decomposes $S$ into a finite
number of complementary components $C^1,\dots,C^u$. Among these
complementary components there are
components $C^{1},\dots,C^s$ which contain essential simple closed
curves not homotopic into a puncture, and there
are components $C^{s+1},\dots,C^u$ which are
topological discs or once
punctured topological discs. As in the beginning of this
section, we mark a point on each smooth boundary component
of the sets $C^i$.

For $k\leq u$ the boundary $\partial C^k$ of $C^k$ consists of a
finite number of connected components. Each of these
components consists of finitely many sides.
Such a side
either is a simple closed curve or
an arc terminating at two not necessarily distinct
cusps of the component, and it
can be represented as a not necessarily
embedded trainpaths on $\sigma$. The union
of these trainpaths determine a
closed curve on $S$ which is freely homotopic to a simple closed
curve, and this simple closed curve is essential if
and only if $k\leq s$, i.e. if the component $C^k$
is different from a disc or a once punctured disc.
In particular, for every $k\leq s$ there is a bordered oriented
surface $S^k\subset C^k$ with smooth boundary
$\partial S^k$ consisting of
a finite number of embedded circles and such that $C^k-S^k$
is a finite union of essential open
annuli, one annulus for each boundary component
of $S^k$. The Euler characteristic of $S^k$ is negative.
Define $S_0=\cup_{k=1}^s S^k$; then $S_0$ is an embedded
bordered subsurface of $S$.

Let ${\cal M}_0(S_0)$ be the subgroup of finite index of the
mapping class group of the bordered
surface $S_0$ consisting of all mapping classes which can
be represented by diffeomorphisms preserving each
component of $S_0$ and each of the
boundary components of $S_0$. By convention, the group
${\cal M}_0(S_0)$ contains in its center the free abelian group
generated by the Dehn twists about the boundary components of
$S_0$.
Then ${\cal M}_0(S_0)$ can naturally be
viewed as a subgroup of ${\cal M}(S)$ consisting
of mapping classes which can be represented by
diffeomorphisms fixing the complement of a small neighborhood of
$S_0$ pointwise. Since $\sigma$ is
contained in $S-S_0$, the group ${\cal M}_0(S_0)$ acts on
the set ${\cal E}(\sigma)\subset {\cal V}({\cal T\cal T})$
of complete extensions
of $\sigma$. Since each of the complementary components
$C^k$ $(s+1\leq k\leq u)$ of
$\sigma$ which is contained in $S-S_0$ is a topological disc
or a once punctured topological disc,
for $s+1\leq k\leq u$ the boundary $\partial C^k$ of $C^k$
is connected and contains at least one cusp.
Namely, otherwise this boundary
defines a closed trainpath on $\sigma$
which necessarily is an essential curve on $S$.
There are only finitely many equivalence classes of graphs
$\tau_{C^k}$ in the sense defined in the beginning of this
section
which occur as the closures of the intersection
of a complete train track $\tau$ with $C^k$. As a consequence,
there is a number
$r>0$ only depending on the topological type
of $S$ such that there are at most $r$ orbits in ${\cal E}(\sigma)$
under the action of ${\cal M}_0(S_0)$. Moreover, there
is a number $q>0$ with the following property.
For any $\tau,\eta\in {\cal E}(\sigma)$ there is
an element $\Theta(\tau,\eta)\in {\cal M}_0(S_0)$
(see \cite{H06a}) such that
$i_\sigma(\tau,\Theta(\tau,\eta)\eta)\leq q$.

Let $\{\sigma(i)\}_{0\leq i\leq \ell}$ be a
recurrent sequence of splits of
$\sigma=\sigma(0)$. Since we only allow right or left splits, for
each $i$ there is a natural
diffeomorphism $\phi_i$ of $S-\sigma$ onto $S-\sigma(i)$
which can be chosen to be the identity on $S_0$. The
restriction $\phi_i^k$ of this diffeomorphism to
a complementary component $C^k$ of $\sigma$
maps $C^k$ diffeomorphically
onto a complementary
component $C_i^k$ of $\sigma(i)$ and extends continuously
to a bijection of the cusps of $\partial C^k$ onto the
cusps of $\partial C^k_i$.
Namely, if we split $\sigma(0)$ at a large branch $e$
to the train track
$\sigma(1)$, then for each complementary component $C^k$ of
$\sigma(0)$ there is a unique complementary component $C^k_{1}$
of $\sigma(1)$ which is diffeomorphic
to $C^k$ with a diffeomorphism $\phi_1^k$ which
fixes pointwise the complement of a neighborhood of $e$.
For $k\leq s$ the diffeomorphism $\phi_1^k$
can be chosen to be the identity on
the subsurface $S^k\subset C^k$ and to induce
a natural (oriented) bijection between the sides
of $C^k$ and the sides of $C^k_1$ which
is the identity outside a neighborhood of $e$.
For every smooth boundary component $T$ of $C^k$ there is
a marked point $x$ contained in the interior of a branch $b\subset T$.
If $b\not=e$ then we require that
the image of $b$ under the natural bijection of
the branches of $\sigma$ onto the branches of $\sigma(1)$
contains a marked point in its interior, and if $b=e$ then
we place a
marked point in the interior of the unique
losing branch of the split which is
contained in the smooth boundary component of $C_1^k$
corresponding to $T$.
For $i\geq 2$ the maps $\phi_i^k$ are constructed inductively
from the maps $\phi_{i-1}^k$ in this way.

Let $\tau,\eta$ be any complete extensions of $\sigma$ and let
$\Theta=\Theta(\tau,\eta)
\in {\cal M}_0(S_0)$ be such that $i_\sigma(\tau,\Theta(\eta))\leq q$.
Since for every $i\leq \ell$ the map
$\phi_i$ restricts to the identity on $S_0$, for every splitting
sequence $\{\eta(j)\}_{0\leq j\leq n}$ induced by
a splitting sequence
$\{\sigma(i)\}_{0\leq i\leq \ell}$ of $\sigma=\sigma(0)$ and
issuing from $\eta=\eta(0)$, the
sequence $\{\Theta \eta(j)\}_{0\leq j\leq n}$
is a splitting sequence issuing from $\Theta(\eta)$ and induced
by the sequence $\{\sigma(i)\}_{0\leq i\leq \ell}$. For every
complete extension $\xi$ of $\sigma(\ell)$ we have
$d(\xi,\eta(n))\leq d(\xi,\Theta\eta(n))+
d(\Theta\eta(n),\eta(n))$. Moreover,
since the action of ${\cal M}_0(S_0)$ commutes with
the action of the subgroup of
the mapping class group of $S-S_0$ consisting of all mapping
classes which preserve every boundary component of $S-S_0$,
by invariance we have
$d(\Theta\eta(n),\eta(n))\leq
d(\Theta\eta(0),\eta(0))+c$ where $c>0$ is a universal constant.
Thus for the proof of the
lemma, it is enough to show the existence of a number
$b>0$ only depending on the
topological type of $S$ with the following properties.
Let $\tau,\eta$ be any complete extensions of $\sigma$
with $i_\sigma(\tau,\eta)\leq q$ and let
$\{\tau(j)\}_{0\leq j\leq m}\subset{\cal V}({\cal T\cal T}),
\{\eta(k)\}_{0\leq k\leq n}\subset {\cal V}({\cal T\cal T})$
be any two
splitting sequences issuing from the train tracks $\tau(0)=\tau,
\eta(0)=\eta$ which are
induced by the splitting sequence
$\{\sigma(i)\}_{0\leq i\leq \ell }$; then $i_{\sigma(\ell)}
(\tau(m),\eta(n))\leq b$.

To show that this is indeed the case we use
Lemma 4.1 inductively. Namely, let $\nu$ be
a geodesic lamination which is carried by
$\sigma(\ell)$ and fills $\sigma(\ell)$.
Such a lamination exists since the
splitting sequence $\{\sigma(i)\}_{0\leq i\leq \ell}$
is recurrent by assumption, and it is carried
by $\sigma=\sigma(0)$ and fills $\sigma$.
Let $\lambda$ be a complete $\tau(m)$-extension of
$\nu$ and let $\mu$ be a complete $\eta(n)$-extension of
$\nu$. Then $\lambda,\mu$ are complete
$\tau,\eta$-extensions of $\nu$.
Assume that $\sigma(1)$ is obtained from
$\sigma$ by a split at a large branch $e$.
Let $\tau^\prime,\eta^\prime\in {\cal V}
({\cal T\cal T})$ be the $\sigma$-modification of
$\tau,\eta$ at $e$ determined by $\lambda,\mu$.
By the definition of an induced splitting sequence,
there are numbers $m(1)\geq 1,
n(1)\geq 1$ such that $\tau(m(1)),\eta(n(1))$ can be
obtained from $\tau^\prime,\eta^\prime$ by a single
split at $e$, and $\tau(m(1)),\eta(n(1))$ contain
$\sigma(1)$ as a subtrack. We call the modification
which transforms $\sigma$ to $\sigma(1)$ and
$\tau,\eta$ to $\tau(m(1)),\eta(n(1))$ a \emph{move}
and we distinguish two types of move.

{\sl Type 1:} $\sigma$ does not contain a marked point in
the interior of $e$.

By Lemma 4.1 and the above discussion, with $p>0$ as in
Lemma 4.1 we have
\begin{align}
i_{\sigma(1)}(\tau(m(1)),\eta(n(1)))=i_\sigma(\tau^\prime,
\eta^\prime)\leq
i_\sigma(\tau^\prime,\eta)+
p(\chi(\eta,\sigma)-\chi(\eta^\prime,\sigma)) \notag\\
\leq
i_\sigma(\tau,\eta)+
p\bigl[\chi(\tau,\sigma)-\chi(\tau(m(1)),\sigma(1))+
\chi(\eta,\sigma)-\chi(\eta(n(1)),\sigma(1))\bigr]\notag
\end{align} and
therefore $i_{\sigma(1)}(\tau(m(1)),\eta(n(1)))\leq
i_\sigma(\tau,\eta)+2r$ for a universal constant $r>0$.
Moreover, we have equality $i_{\sigma(1)}(\tau(m(1)),\eta(n(1)))=
i_\sigma(\tau,\eta)$
if $\chi(\tau(m(1)),\sigma(1))=\chi(\tau,\sigma)$ and
$\chi(\eta(n(1)),\sigma(1))=\chi(\eta,\sigma)$.

{\sl Type 2:} There is a marked point of $\sigma$
contained in the interior
of the branch $e$.

Using again Lemma 4.1 and the
above discussion,
we conclude as before that $i_{\sigma(1)}(\tau(m(1)),\eta(n(1)))
\leq i_\sigma(\tau,\eta)+2p$.

Let $m(\ell)>0,n(\ell)>0$ be the smallest numbers such that the
train tracks $\tau(m(\ell)),\eta(n(\ell))$ contain $\sigma(\ell)$
as a subtrack. By the definition of a splitting sequence induced
from the sequence $\{\sigma(i)\}_{0\leq i\leq \ell}$, the train
tracks $\tau(m(\ell)),\eta(n(\ell))$ can be obtained from
$\tau,\eta$ by $\ell$ moves. The train tracks $\tau(m),\eta(n)$
contain $\sigma(\ell)$ as a subtrack and are obtained from
$\tau(m(\ell))$ by a sequence of splits at large proper
subbranches of $\sigma(\ell)$. Since $\sigma(\ell)$ does not have
any closed curve components, the length of a splitting sequence
connecting $\tau(m(\ell)),\eta(n(\ell))$
to $\tau(m),\eta(n)$ is bounded from above by a
universal constant (compare the detailed discussion in Section 4
of \cite{H06b}). In particular, the distance between
$\tau(m)$ and $\eta(n)$ is bounded from above
by $d(\tau(m(\ell)),\eta(n(\ell))+\tilde c$ where
$\tilde c>0$ is a universal constant.
Since moreover
$\chi(\tau,\sigma)+\chi(\eta,\sigma)$ is bounded from above by a
universal constant, the lemma follows if we can show that the
number of times a move of type 2 occurs in our splitting sequences
is bounded from above by a universal constant.

However, if for some $i<\ell$ the split $\sigma(i)\to
\sigma(i+1)$ is a split at a large branch $e_i$ containing
a marked point of $\sigma$ in its interior then
$e_i$ is contained in a smooth component $A$ of the boundary of
a complementary component $C_i^k$ of $\sigma(i)$. In particular,
$A$ defines a closed trainpath $\rho$ on $\sigma(i)$ which
is freely homotopic to a simple closed curve
defining a boundary component of $C_i^k$.
Thus the trainpath $\rho$ passes through any
branch of $\sigma(i)$ at most twice and hence
its length is uniformly bounded.
The train track
$\sigma(i+1)$ obtained from $\sigma(i)$
by a single split at $e_i$
contains $A$ as an embedded trainpath $\rho^\prime$
whose length is strictly smaller than the length of $\rho$.
As a consequence, the number of moves of type 2 is bounded from
above by a universal constant.
This completes the proof of the lemma.
\end{proof}

Next we
estimate distances
in ${\cal T\cal T}$ between train tracks
which do not carry a common geodesic lamination.
We show.

\begin{lemma}\label{reverse}
For every $R>0$ there is a number
$\beta_0=\beta_0(R)
>0$ with the following property. Let $\tau,\eta\in
{\cal V}({\cal T\cal T})$ with $d(\tau,\eta)\leq R$ and let
$\tau^\prime,\eta^\prime$ be complete train tracks which can be
obtained from $\tau,\eta$ by any splitting sequence.
If $\tau,\eta$ do not carry any common geodesic
lamination then there is a point
$\tau^{\prime\prime}\in {\cal V}({\cal T\cal T})$ in the
$\beta_0(R)$-neighborhood of $\tau^\prime$ which can be connected
to a point $\eta^{\prime\prime}\in {\cal V}({\cal T\cal T})$
contained in the $\beta_0(R)$-neighborhood of $\eta^\prime$
by a splitting sequence which passes through the
$\beta_0(R)$-neighborhood of $\tau$.
Moreover, for any complete
geodesic lamination $\nu$ which is carried by
$\eta^\prime$, we can assume that $\eta^{\prime\prime}$
carries $\nu$.
\end{lemma}

\begin{proof} For a fixed number $R>0$ there
are only finitely many orbits under the
action of the mapping class group of
pairs $(\tau,\eta)\in {\cal V}({\cal T\cal T})\times
{\cal V}({\cal T\cal T})$ where $d(\tau,\eta)\leq R$ and such
that $\tau,\eta$ do not carry any common
geodesic lamination. Thus by invariance
under the mapping class group it is enough
to show the lemma for two fixed train tracks
$\tau,\eta\in {\cal V}({\cal T\cal T})$ which do not
carry any common geodesic lamination and with
a constant $\beta_0>0$ depending on $\tau,\eta$.

For a complete train track $\xi$ denote
by ${\cal C\cal L}(\xi)$ the set of all complete
geodesic laminations which are carried by $\xi$.
By Lemma 2.4 of \cite{H06a}, the set ${\cal C\cal L}(\xi)$
is open and closed in ${\cal C\cal L}$.
We first show that there are finitely many
complete train tracks $\tau_1,\dots,\tau_\ell$ and
$\eta_1,\dots,\eta_m$ with the following
properties.
\begin{enumerate}
\item For each
$i\leq \ell$ the train track $\tau_i$ is carried
by $\tau$ and $\cup_i{\cal C\cal L}(\tau_i)=
{\cal C\cal L}(\tau)$.
\item For each
$j\leq m$ the train track $\eta_j$ is carried
by $\eta$ and $\cup_j{\cal C\cal L}(\eta)=
{\cal C\cal L}(\eta)$.
\item For all $i\leq \ell,j\leq m$ the train tracks
$\tau_i,\eta_j$ hit efficiently.
\end{enumerate}

For this observe that since $\tau,\eta$ do not
carry any common geodesic lamination,
every lamination $\lambda\in {\cal C\cal L}(\tau)$
intersects every lamination $\mu\in {\cal C\cal L}(\eta)$
transversely. Namely, a complete geodesic lamination
decomposes the surface $S$ into ideal triangles and
once punctured monogons. Thus if $\ell$ is any
simple geodesic on $S$ whose closure in $S$ is compact and if
$\ell$ does \emph{not} intersect
the complete geodesic lamination $\mu$ transversely
then $\ell$ is contained in $\mu$ and hence the
closure of $\ell$ is a sublamination
of $\mu$.
Now if $\ell$ is a leaf of the complete geodesic lamination
$\lambda$ then the closure of $\ell$
is a sublamination of $\lambda$ as well.
Since $\lambda$ is carried by $\tau$ and $\mu$ is
carried by $\sigma$, this violates our assumption
that $\tau,\sigma$ do not carry a common geodesic
lamination.

As a consequence, for every
lamination $\lambda\in {\cal C\cal L}(\tau)$,
every train track $\xi$ which is sufficiently
close to $\lambda$ in the Hausdorff topology
hits every lamination $\nu\in {\cal C\cal L}(\eta)$
efficiently. Thus by Lemma 2.2 and Lemma 2.3 of
\cite{H06a}, for every complete geodesic lamination
$\lambda\in {\cal C\cal L}(\tau)$ there is a
train track $\tau(\lambda)\in {\cal V}({\cal T\cal T})$
which carries $\lambda$, which
is carried by $\tau$ and which hits every
lamination $\nu\in {\cal C\cal L}(\eta)$ efficiently.
The set ${\cal C\cal L}(\tau(\lambda))$ is
an open subset of ${\cal C\cal L}(\tau)$ and
therefore by compactness of ${\cal C\cal L}(\tau)$ there
are finitely many laminations
$\lambda_1,\dots,\lambda_\ell\in {\cal C\cal L}(\tau)$ such that
${\cal C\cal L}(\tau)=\cup_i{\cal C\cal L}(\tau(\lambda_i))$.
Write $\tau_i=\tau(\lambda_i)$.

Every lamination $\nu\in {\cal C\cal L}(\eta)$
hits each of the train tracks $\tau_i$ efficiently.
As a consequence, if $\xi$ is a train track
which is sufficiently close to $\nu$ in the
Hausdorff topology then $\xi$ hits each of the
train tracks $\tau_i$ $(i\leq \ell)$ efficiently.
As before, this implies that we can find a finite
family $\eta_1,\dots,\eta_m\in {\cal V}({\cal T\cal T})$
of train tracks
which are carried by $\eta$, which hit each of the
train tracks $\tau_i$ $(i\leq \ell)$ efficiently and such that
$\cup_j{\cal C\cal L}(\eta_j)={\cal C\cal L}(\eta)$.
This shows our above claim.

Let $k=\max\{d(\tau,\tau_i),d(\eta,\eta_j)\mid
i\leq \ell, j\leq m\}$. Let
$\tau^\prime,\eta^\prime$ be obtained
from $\tau,\eta$ by a splitting sequence and let
$\lambda\in {\cal C\cal L}(\tau^\prime)$ be
a complete geodesic lamination carried by $\tau^\prime$.
Then $\lambda\in {\cal C\cal L}(\tau)$ and hence
by our above construction, there is some $i\leq \ell$
such that $\tau_i$ carries $\lambda$.
By Corollary 4.6 of \cite{H06a} there is a complete
train track
$\xi$ which carries $\lambda$, is
carried by both $\tau^\prime$ and $\tau_i$ and
whose distance to $\tau^\prime$ is bounded from above
by a universal constant $q>0$ only depending on $k$.
Similarly, for a complete train track $\eta^\prime$ which
can be obtained from $\eta$ by a splitting sequence
and for some $\nu\in {\cal C\cal L}(\eta^\prime)$
there is some $j\leq m$ and
a train track $\zeta$
which is carried by both $\eta^\prime$ and $\eta_j$,
which carries $\nu$ and whose distance to $\eta^\prime$
does not exceed $q$.

Since $\tau_i,\eta_j$
hit efficiently and $\tau_i$ carries $\xi$,
the train tracks $\xi,\eta_j$ hit efficiently.
In particular, the geodesic lamination
$\nu$ hits $\xi$ efficiently. Therefore
by Proposition 3.2, a $\nu$-collapse
$\xi^*$ of the dual bigon track $\xi_b^*$
of $\xi$ carries a train track $\beta$
which carries $\nu$ and can be obtained from
$\eta_j$ by a splitting and shifting sequence of length
at most $q$. In particular, the distance
between $\eta_j$ and $\beta$ and hence
the distance between $\tau$ and $\beta$ is
uniformly bounded.
Since $\eta_j$ carries $\zeta$ and
$\zeta$ carries $\nu$, we conclude from Corollary 4.6
of \cite{H06a}
that $\beta$ carries a train track $\sigma$
contained in a uniformly bounded neighborhood of
$\zeta$ and hence in a uniformly bounded 
neighborhood of $\eta^\prime$.

The distance between
$\tau^\prime$ and $\xi$ is uniformly bounded and therefore
the distance between $\xi^*$ and
$\tau^\prime$ is uniformly bounded as well
(compare Proposition 3.2 and the following remark).
On the other hand, since $\beta$
is carried by $\xi^*$ the train track
$\xi^*$ can be connected to $\beta$ by a splitting
and shifting sequence (Theorem 2.4.1 of \cite{PH92}).
We deduce from Lemma 5.4 of
\cite{H06a} that $\xi^*$ is splittable to train tracks
$\beta^\prime,\sigma^\prime$ which carry $\nu$ and
are contained in a uniformly bounded neighborhood of
$\beta,\sigma$. It then follows from Lemma 5.1 of \cite{H06a}
that $\beta^\prime$ is splittable to a train track
$\sigma^{\prime\prime}$ which carries $\nu$ and is
contained in a uniformly bounded neighborhood of
$\sigma^\prime$. Since the distance between $\sigma^\prime$ and
$\eta^\prime$ is uniformly bounded,
this implies the lemma.
\end{proof}

Our next goal is to obtain a suitable extension
of Lemma 4.3
to train tracks $\tau,\eta\in {\cal V}({\cal T\cal T})$
containing a common subtrack which is a union of simple closed
curves. We first have to overcome a technical difficulty
arising from the fact that the split of a complete train
track need not be recurrent. For this call a large branch $e$ of a
complete train track $\eta$ \emph{rigid} if only one of the two
train tracks obtained from $\eta$ by a split at $e$ is complete.
Note that $e$ is a rigid large branch of $\eta$
if and only if the (unique) complete train track
obtained from $\eta$ by a split at $e$ carries \emph{every}
complete geodesic lamination which is carried by $\eta$. It
follows from the results of \cite{PH92} that a branch $e$ in
$\eta$ is rigid if and only if the train track $\zeta$ obtained
from $\eta$ by a collision, i.e. a split followed by the removal
of the diagonal of the split, is \emph{not} recurrent. We have.

\begin{lemma}\label{rigid} There is a number $a_1>0$ with the following
property. For every $\eta\in {\cal V}({\cal T\cal T})$ there is a
splitting sequence $\{\eta(i)\}_{0\leq i\leq s}$ issuing from
$\eta=\eta(0)$ of length $s\leq a_1$ such that for every $i$,
$\eta(i+1)$ is obtained from $\eta(i)$ by a split at a rigid large
branch and such that $\eta(s)$ does not contain any rigid large
branches.
\end{lemma}

\begin{proof} We show the existence of a constant $a_1>0$ as in the
lemma with an argument by contradiction. Assume that our claim
does not hold. Then there is a sequence of pairs
$(\beta_i,\xi_i)\in {\cal V}( {\cal T\cal T})$ such that $\xi_i$
can be obtained from $\beta_i$ by a splitting sequence of length
at least $i$ consisting of splits at rigid large branches. Every
complete geodesic lamination which is carried by $\beta_i$ is also
carried by $\xi_i$. By invariance under the action of the mapping
class group ${\cal M}(S)$ and the fact that there are only
finitely many ${\cal M}(S)$-orbits on ${\cal V}( {\cal T\cal T})$,
by passing to a subsequence we may assume that there is some
$\eta\in {\cal V}({\cal T\cal T})$ such that $\beta_i=\eta$ for
all $i$. Since $\eta$ has only finitely many large branches, by a
standard diagonal procedure we construct an \emph{infinite}
splitting sequence $\{\eta(i)\}_{i\geq 0}$ issuing from
$\eta=\eta(0)$ such that for every $i$ the train track $\eta(i+1)$
is obtained from $\eta(i)$ by a split at a rigid large branch.
Then for every $i$, every complete geodesic lamination which is
carried by $\eta$ is also carried by $\eta(i)$. Now for every
\emph{projective measured geodesic lamination} whose support $\nu$
is carried by $\eta$ there is a complete geodesic lamination
$\lambda$ carried by $\eta$ which contains $\nu$ as a
sublamination \cite{H06a}. Thus for each $i$, the space ${\cal
P\cal M}(i)$ of projective measured
geodesic laminations carried by
$\eta(i)$ coincides with the space ${\cal P\cal M}(0)$ of
projective measured geodesic laminations carried by $\eta$. On the
other hand, the set ${\cal P\cal M}(0)$ contains an open subset of
the space of all projective measured geodesic laminations since
$\eta$ is complete (see \cite{PH92}). Therefore $\cap_i{\cal P\cal
M}(i)={\cal P\cal M}(0)$ contains an open subset of projective
measured geodesic
lamination space which contradicts Theorem 8.5.1 in
\cite{M03}. This shows our lemma.
\end{proof}

To each complete geodesic lamination $\lambda$ on $S$ and every
simple closed curve component $c$ of $\lambda$ we associate a sign
${\rm sgn}(\lambda,c)\in \{+,-\}$ as follows. Let $S_1$ be the
surface obtained from $S$ by cutting $S$ open along $c$. We view
$S_1$ as a bordered surface with two boundary components
$c_1,c_2$ corresponding to $c$. The orientation of $S$ then
induces a boundary orientation for $c_1,c_2$. For our complete
geodesic lamination $\lambda$ containing $c$ as a minimal
component and for $i=1,2$ there is at least one leaf of $\lambda$
which is contained in $S_1$ and spirals about $c_i$. We associate
to $\lambda$ and $c$ the sign ${\rm sgn}(\lambda,c)=+$ if the
spiraling direction of such a leaf coincides with the boundary
orientation of $c_i$ for $i=1,2$, and we associate to $\lambda$
and $c$ the sign ${\rm sgn}(\lambda,c)=-$ otherwise. Since
$\lambda$ is complete by assumption, if ${\rm sgn}(\lambda,c)=-$
then the spiraling direction of a leaf of $\lambda$ spiraling
about $c_i$ $(i=1,2)$ is opposite to the orientation of $c_i$ as a
component of the boundary of $S_1$ (compare Section 2 of
\cite{H06a}).

Next we look at a complete train track $\tau$ containing a
subtrack $\beta$ which is a union of $k\geq 1$ simple
closed curves embedded in $\tau$. In other words, $\beta$ is a
subtrack of $\tau$ without large branches. Such a train track
$\tau$ carries a complete geodesic lamination $\lambda$ which
contains the $k$ simple closed curve components $c_1,\dots,c_k$ of
$\beta$ as minimal components. By the above, $\lambda$ determines
for each of the curves $c_i$ a sign. Recall the definition of a
\emph{Dehn twist} about an essential simple closed curve $c$ in
$S$. We define the twist to be \emph{positive} if the direction of
the twist coincides with the direction given by the boundary
orientation of $c$ in the surface $S_1$ obtained by cutting $S$
open along $c$. We use these signs to estimate distances in ${\cal
T\cal T}$ obtained by splitting complete train tracks
along simple closed curve
subtracks. We show.

\begin{lemma}\label{Dehn} There is a constant $a_2>0$ with the
following property. Let $\tau\in {\cal V}({\cal T\cal T})$, let
$c<\tau$ be an embedded simple closed curve and let $\phi_c$ be
the positive Dehn twist about $c$. Let $\{\tau_i\}_{0\leq i\leq
m}$ be a splitting sequence issuing from $\tau_0=\tau$ which
consists of $c$-splits at large proper subbranches of $c$. Let
$\lambda\in {\cal C\cal L}$ be a complete geodesic lamination
which is carried by $\tau$ and contains $c$ as a minimal
component. Then there is some $i\geq 0$ such that the distance
between $\tau_m$ and the train track $\phi_c^{{\rm
sgn}(\lambda,c)i}\tau$ is at most $a_2$.
\end{lemma}

\begin{proof}
A \emph{standard twist connector} in a
complete train track $\xi$ is an embedded closed curve $\alpha$
in $\xi$ which consists of a large branch and a small branch,
connected at two switches. If $\xi^\prime$ is obtained from
$\xi$ by an $\alpha$-split at the large branch in $\alpha$,
then $\xi^\prime$ is obtained from $\xi$ by a half-Dehn twist
about $\alpha$ whose sign is determined by the neighbors
of the subtrack $\alpha<\xi$ (compare the discussion in
Section 2 of \cite{H06a}).
In particular, for every complete geodesic lamination
$\lambda$ which is carried by $\xi$ and contains $\alpha$ as a minimal
component, the sign ${\rm sgn}(\lambda,\alpha)$ is determined by the
twist connector and hence does not depend on $\lambda$.

By the considerations in Section 4 of \cite{H06b},
there is a number $a_0>0$ only depending on the topological type
of $S$ such that for every complete train track $\tau$ containing
an embedded simple closed curve $c$ the image $\eta$ of $\tau$
under a sequence of $c$-splits of length at most $a_0$ contains $c$
as a \emph{simple vertex cycle}, i.e. $\eta$ can be shifted to a
train track $\eta^\prime$ which contains $c$ as a standard
twist connector.
Thus any sequence of $c$-splits modifies $\eta$ to a
train track $\eta_1$ which is
contained in uniformly bounded neighborhood of the image of
$\eta$ under a multiple of the Dehn twist along $c$
whose direction is determined by ${\rm sgn}(\lambda,c)$
where $\lambda$ is a complete extension of $c$ carried by $\eta$.
From this the lemma is immediate.
\end{proof}

A \emph{simple multicurve} $c$ on $S$ consists of a finite collection
$c_1,\dots, c_k$ of essential simple closed curves which are not
mutually freely homotopic and which can be realized disjointly. A
simple multicurve $c$ can be viewed as a train track without large
branches by adding a single switch to each component of $c$. A
\emph{Dehn multitwist} of a multicurve $c=\cup_ic_i$ is an element
$\phi\in {\cal M}(S)$ which can be represented in the form
$\phi=\phi_{c_1}^{m_1}\circ\dots\circ \phi_{c_k}^{m_k}$ for some
$m_i\in \mathbb{Z}$ and where as before, $\phi_{c_i}$ is the
positive Dehn twist about $c_i$. The next lemma is an extension of
Lemma 4.3. For its formulation (and later use),
for a train track $\tau\in
{\cal V}({\cal T\cal T})$ which is splittable to
a train track $\tau^\prime\in {\cal V}({\cal T\cal T})$
define the \emph{flat
strip} $E(\tau,\tau^\prime)$ to be the maximal subgraph
of ${\cal T\cal T}$ whose vertices consist of the
collection of all train tracks
which can be obtained from $\tau$ by a splitting sequence and are
splittable to $\tau^\prime$.
We have.

\begin{lemma}\label{multitwist} For every $R>0$ there is a number
$\beta_1(R)>0$ with the following property. Let $\tau,\eta\in
{\cal V}({\cal T\cal T})$ with $d(\tau,\eta)\leq R$.
Assume that
$\tau,\eta$ contain a common subtrack $\sigma$ which is a
multicurve and that the components of this multicurve are
precisely the minimal geodesic laminations which are carried by
both $\tau$ and $\eta$.
Let
$\tau^\prime,\eta^\prime$ be complete train tracks which can be
obtained from $\tau,\eta$ by a splitting sequence.
Then there is a Dehn multitwist $\phi$
about $c$ and a point $\tau^{\prime\prime}$ in the
$\beta_1(R)$-neighborhood of $\tau^\prime$ which can be connected
to a point $\eta^{\prime\prime}$ in the
$\beta_1(R)$-neighborhood of $\eta^\prime$
by a splitting sequence passing
through the $\beta_1(R)$-neighborhood of $\phi \tau$. Moreover,
the distance between $\phi(\tau)$ and
$E(\tau,\tau^\prime),E(\eta,\eta^\prime)$
is at most $\beta_1(R)$.
\end{lemma}

\begin{proof} As in the proof of Lemma 4.3, by invariance under the
mapping class group it suffices to show the lemma for some fixed
train tracks $\tau,\eta\in {\cal V}({\cal T\cal T})$
which satisfy
the assumptions in the lemma. In particular, $\tau,\eta$ contain a
common subtrack $c$ which is a simple
multicurve and such that
every minimal geodesic lamination which is carried by both $\tau$ and
$\eta$ is a component of $c$.
By Lemma 4.3, we may assume that $c$ is not empty.

Let $c_1,\dots,c_k$ be the components of $c$. Then for
$i\not=j$, a $c$-split of $\tau$ at
a large proper subbranch of $c_i$
commutes with a $c$-split of $\tau$ at a large proper
subbranch of $c_j$. Assume that $\tau$ is splittable to a train
track $\tau^\prime$. By Lemma 4.4, via replacing $\tau^\prime$ by
its image under a splitting sequence of uniformly bounded length
we may assume that $\tau^\prime$ does not contain any rigid
large branches.
Define inductively a sequence $\{\tau(i)\}_{0\leq i\leq
k}\subset E(\tau,\tau^\prime)$ consisting of train
tracks $\tau(i)$ which contain $c$ as a subtrack as
follows. Put $\tau(0)=\tau$ and assume that $\tau(i-1)$ has been
defined for some $i\in \{1,\dots,k\}$. Let $\tau(i)\in
E(\tau,\tau^\prime)$ be the train track which
can be obtained from $\tau(i-1)$ by a
sequence of $c_i$-splits of maximal length
at proper large subbranches of $c_i$.
Since splits at large subbranches of $c_i,c_j$
for $i\not=j$ commute,
the train track $\tau(k)$ only depends on $\tau,\tau^\prime,c$
but not on the ordering of the components $c_i$ of
$c$. By Lemma 4.5, there are numbers $b_i\in \mathbb{Z}$ such that
for $\phi_\tau=\phi_{c_1}^{b_1}\circ \dots\circ
\phi_{c_k}^{b_k}\in {\cal M}(S)$
the distance between $\tau(k)$ and the train track
$\phi_\tau(\tau)$ is bounded from above by $a_2$.

Similarly, if $\eta$ is splittable to a train track
$\eta^\prime\in {\cal V}({\cal T\cal T})$
without rigid large branches, then there are numbers
$p_i\in \mathbb{Z}$ such that for $\phi_\eta=
\phi_{c_1}^{p_1}\circ\dots\circ \phi_{c_k}^{p_k}\in
{\cal M}(S)$ the distance between the train track
$\eta(k)$ obtained from $\eta$ and $\eta^\prime$ by the above
procedure and $\phi_\eta(\eta)$
is bounded from above by $a_2$.

For $i\leq k$ define $m(i)=0$ if the signs of $b_i,p_i$ are
distinct, and if the signs of $b_i,p_i$
coincide then define
$m(i)={\rm sgn}(a_i)\min\{\vert a_i\vert, \vert b_i\vert \}$. Write
$\phi=\phi_{c_1}^{m(1)}\circ\dots\circ \phi_{c_k}^{m(k)}$. By the
choice of the multi-twist $\phi$ and by invariance of the distance
function on ${\cal T\cal T}$ under the action of the mapping class
group, there is a number $\chi >0$ only depending on
$d(\tau,\eta)$ and there are
train tracks $\tau_1\in E(\tau,\tau^\prime),
\eta_1\in E(\eta,\eta^\prime)$ contained in the
$\chi$-neighborhood of
$\phi(\tau),\phi(\eta)$.
Moreover, by the choice of $\phi$,
we may assume that $\tau_1,\eta_1$ contain the simple multicurve
$c$ as a subtrack and that for each $i$
one of the following two possibilities is satisfied.
\begin{enumerate}
\item $m(i)=0$ and a splitting sequence connecting $\tau,
\eta$ to $\tau_1,\eta_1$ does not contain any split at a large
proper subbranch of  $c_i$.
\item If $\vert b_i\vert \leq \vert p_i\vert$ then
the flat strip $E(\tau,\tau^\prime)$ does not contain
any train track which can be obtained from $\tau_1$
by a $c_i$-split at a large proper subbranch of $\tau_1$
and similarly for $\eta_1$ in the case that
$\vert p_i\vert \leq \vert b_i\vert$.
\end{enumerate}

After reordering we may assume that
there is some $s\leq k$ such that
$m(i)=0$ for $i\leq s$ and
that $m(i)\not=0$ for $i\geq s+1$. Let $i\geq s+1$ and
assume without loss of generality
that $\vert b_i\vert \leq \vert p_i\vert$. Since
the train tracks $\tau_1,\eta_1$ are complete and contain
the curve $c_i$ as a subtrack, both $\tau_1$ and $\eta_1$ contain
a large proper subbranch of $c_i$.
By the choice of $\tau_1$, for each such large subbranch
$e$ of $c_i$, the train track obtained from $\tau_1$ by
a $c_i$-split at $e$ is \emph{not} contained in the flat
strip $E(\tau,\tau^\prime)$.
Let $\hat\tau_1$ be the train
track obtained from $\tau_1$ by the split at $e$ which is
\emph{not} a $c_i$-split. Then either we have
$\hat\tau_1\in E(\tau,\tau^\prime)$, in particular
$\hat\tau_1$ is complete,
or no train track which can be obtained from
$\tau_1$ by a split at $e$ is splittable
to $\tau^\prime$. In the second case, $e$ is a large
branch of $\tau^\prime$ by uniqueness of splitting
sequences (Lemma 5.1 of \cite{H06a}).
Since $\tau^\prime$ does not
contain any rigid large branch by assumption,
both train tracks which can be obtained from $\tau^\prime$
by a single split at $e$ are complete. As a consequence,
the train track $\hat\tau_1$ is complete and is
splittable to a complete train track which can be obtained from
$\tau^\prime$ by a single split at $e$.
The train track $\hat\tau_1$
does not carry the simple
closed curve $c_i$.

Since splits at large proper subbranches of the distinct
components of $c$ commute, we construct in this way
successively in $k-s$ steps
from the train tracks $\tau_1,\eta_1$
complete train tracks $\tau_2,\eta_2$ which can be
obtained from $\tau_1,\eta_1$ by a splitting sequence
of length at most $k-s$. The train tracks
$\tau_2,\eta_2$ are splittable to
train tracks $\tau_2^{\prime},\eta_2^{\prime}$
which can be obtained from $\tau^\prime,\eta^\prime$
by splitting sequences of length at most $k-s$.
A minimal component of a geodesic
lamination which is carried by both $\tau_2,\eta_2$ coincides with
one of the curves $c_i$ for $i\leq s$. Moreover,
the train tracks
$\tau_2,\eta_2$ contain the simple closed curves $c_1,\dots,c_s$ as
embedded subtracks, and their distance
in ${\cal T\cal T}$ is bounded from above by a universal constant
only depending on $d(\tau,\eta)$ and the topological type
of $S$.

By the considerations in Section 4 of \cite{H06b}, after
reordering of the components $c_i$ and after possibly
replacing $\tau_2,\eta_2$ by their images under a splitting
sequence of uniformly bounded length and which are splittable to
$\tau_2^\prime,\eta_2^\prime$ we may assume that there is a number
$u\leq s$ with the following property.
For each $i>u$, either the train
track $\tau_2$ or the train track $\eta_2$
contains a large proper subbranch $e$ of $c_i$ with the property
that a $c_i$-split of $\tau_2$ (or $\eta_2$) is not splittable to
$\tau_2^\prime$ (or $\eta_2^\prime$). For every $i\leq u$, the
curve $c_i$ is a simple vertex cycle in both $\tau_2,\eta_2$, i.e.
$\tau_2$ and $\eta_2$
can be shifted to train tracks which contain $c_i$
as a twist connector. Apply the above construction to the train
tracks $\tau_2,\eta_2$ and the simple closed curves $c_i$
$(u+1\leq i\leq s)$. We obtain train tracks $\tau_3,\eta_3$ of
uniformly bounded distance which are splittable to train tracks
$\tau_3^\prime,\eta_3^\prime$
obtained from $\tau_2^\prime,\eta_2^\prime$ by a splitting sequence
of uniformly bounded length and such that a minimal component of a
geodesic lamination carried by both $\tau_3,\eta_3$ coincides with
one of the curves $c_1,\dots,c_u$, and these
curves are simple vertex cycles of
$\tau_3,\eta_3$.

Let $\tau_4$ be the train track which can be obtained from
$\tau_3$ by a sequence of $\cup_{i=1}^uc_i$-splits of maximal length and which
is splittable to $\tau_3^\prime$.
For $\psi=\phi_{c_1}^{b_1}\circ \dots\circ \phi_{c_u}^{b_u}$,
the distance between $\tau_4$ and $\psi(\tau_3)$ is
uniformly bounded. Moreover, since the curves $c_i$ are simple
vertex cycles in $\eta_3$ and since for each $i\leq u$ the
signs of $b_i,p_i$ do \emph{not} coincide, the
train track $\psi \eta_3=\eta_4$ is
splittable to $\eta_3$. Replace $\tau_4$ by its image $\tau_5$
under a splitting sequence of uniformly bounded length which is
splittable to the image of $\tau_3^\prime$ under a splitting
sequence of uniformly bounded length and which does not carry any
of the curves $c_i$. Then the distance between $\tau_5$ and
$\eta_4$ is bounded from above by a constant only depending on
$d(\tau,\eta)$, and $\eta_4$ and $\tau_5$ do not carry any common geodesic
lamination. Moreover, $\eta_4$ is splittable to $\eta_3^\prime$ with
a splitting sequence which passes through $\eta_3$.

Let $\nu$ be a complete geodesic lamination which is carried by
$\eta_3^\prime$.
By Lemma 4.3, applied to the train tracks $\tau_5,\eta_4$
which are splittable to the train tracks $\tau_3^\prime,
\eta_3^\prime$ in uniformly bounded neighborhoods
of $\tau^\prime,\eta^\prime$,
there is a splitting sequence
$\{\alpha(i)\}_{0\leq i\leq n}$ which connects a train track
$\alpha(0)$ in a uniformly bounded neighborhood of $\tau^\prime$
to a train track $\alpha(n)$ in a uniformly bounded neighborhood
of $\eta^\prime$ which carries $\nu$, and such that this splitting
sequence passes through a uniformly bounded neighborhood of
$\eta_4$. Now by the considerations in Section 5 of \cite{H06a},
this splitting sequence can be chosen to pass through a uniformly
bounded neighborhood of $\eta_3$ as well.
Thus our lemma follows with the Dehn-multisplit
$\phi$ as above which maps $\tau,\eta$ into a uniformly
bounded neighborhood of $\tau_1,\eta_1$.
\end{proof}

For a train track $\tau\in {\cal V}({\cal T\cal T})$ containing
a subtrack $\sigma$, recall from the beginning of this
section the definition
of a splitting sequence of $\tau$ induced by a
splitting sequence of $\sigma$.
We have.

\begin{lemma}For every train track $\tau\in {\cal V}({\cal
T\cal T})$ which is splittable to a train track $\eta$ and every
recurrent subtrack $\sigma$ of $\tau$ without closed curve
components there is a unique train track $\tau^\prime\in
E(\tau,\eta)$ with the following properties.
\begin{enumerate}
\item There is a recurrent
splitting sequence $\{\sigma(i)\}_{0\leq i\leq p}$
issuing from $\sigma(0)=\sigma$ such that
$\tau^{\prime}$ can be obtained from
$\tau$ by a splitting sequence induced by $\{\sigma(i)\}$.
\item If $\tilde \tau\in E(\tau,\eta)$ can be
obtained from $\tau$ by a recurrent splitting sequence induced
by a sequence of splits of $\sigma$ then $\tilde\tau$
is splittable to $\tau^\prime$.
\end{enumerate}
\end{lemma}

\begin{proof} Let $\tau\in
{\cal V}({\cal T\cal T})$ be a complete train track
which
is splittable to a train track $\eta\in {\cal V}({\cal T\cal T})$ and
let $\sigma$ be
a subtrack of $\tau$. Recall that a \emph{$\sigma$-split}
of $\tau$ is a split of $\tau$ at a large proper
subbranch of $\sigma$ with the property that the
split track contains $\sigma$ as a subtrack.
Similar to the procedure in the proof of Lemma 4.3 of
\cite{H06b} we define
a (non-deterministic) finite algorithm
which takes $\tau$ as input and yields a finite sequence
$\{\zeta(i)\}_i\subset E(\tau,\eta)$
of train tracks with $\tau=\zeta(0)$ and such
that for each $i$, either $\zeta(i+1)$ is obtained from $\zeta(i)$
by a single $\sigma$-split at a large proper subbranch
of $\sigma$ or $\zeta(i+1)$ is obtained from $\zeta(i)$ by
putting a mark on a large proper subbranch of $\sigma$.
The initial train track $\tau=\zeta(0)$ does not
have marked branches.

In the $i$-th step $(i\geq 1)$
the algorithm begins with checking for the existence of an
unmarked large proper subbranch $e$ of $\sigma$ in
the train track
$\zeta(i-1)$. If there is no such branch then the algorithm stops.
Otherwise the algorithm chooses such a large
proper subbranch $e$ of $\sigma$
and proceeds as
follows.

If there is a train track $\tilde \zeta(i)\in E(\tau,\eta)$
which can be obtained from $\zeta(i-1)$ by a single split
at $e$ and which contains $\sigma$ as
a subtrack, then define $\zeta(i)$ to be the train track
$\tilde \zeta(i)$
equipped with the markings
obtained from the markings of the branches of $\zeta(i-1)$
via the
natural identification of the branches of $\zeta(i-1)$ with the
branches of $\tilde \zeta(i)$. Otherwise
define $\zeta(i)$ to be the train track $\zeta(i-1)$ equipped
with an additional mark on the branch $e$ (see the proof
of Lemma 4.3 of \cite{H06b}).

It follows from the discussion in Section 4 of \cite{H06b}
that there is a universal constant $q>0$
only depending on the topological type of $S$ such that
our algorithm stops
after at most $q$ steps (see the proof of Lemma 4.3 in \cite{H06b}).
It produces a train track
$\tilde \tau(1)\in E(\tau,\eta)$ which contains
$\sigma$ as a subtrack.
For each large branch $e$ of $\tilde \tau(1)$ contained in
$\sigma$, either $e$ is a large branch in $\sigma$,
i.e. $\tilde \tau(1)$ is tight at $e$, or $e$ is
marked. Since splits at distinct large branches
of a train track $\tau$ commute, the marked train track
$\tilde \tau(1)$ only depends on $\tau,\eta,\sigma$
but not on any choices made.
Moreover, if $e_1,\dots,e_s$ are the large branches
of $\sigma$ which correspond to the unmarked large
branches of $\tilde \tau(1)$ then there is a complete
geodesic lamination carried by $\tilde \tau(1)$
such that the train track obtained from $\tau$
by a successive $\sigma$-modification at the branches
$e_1,\dots,e_s$ with respect to $\lambda$ is splittable
to $\tilde \tau(1)$ with a splitting sequence
of uniformly bounded length.

After reordering, there is a number $\hat k\in \{0,\dots,s\}$
such that for every $i\leq \hat k$ the flat strip
$E(\tau,\eta)$ contains a train track $\xi_i$
which
can be obtained from $\tilde \tau(1)$ by a single
right or left split at
$e_i$. For each $i$, there is a unique
choice of a right or left split of $\sigma$ at
$e_i$ such that the resulting train track
$\tilde \sigma_i$ is a subtrack of $\xi_i$.
After reordering, we may assume that there is
some $k\leq \hat k$ such that for every
$i\leq k$ the train track $\tilde \sigma_i$ is
recurrent but that this is not the case
for the train tracks $\sigma_j$ for $k<j\leq \hat k$.
Define $\tau(1)$ to be the train track obtained from
$\tilde \tau(1)$ by splitting $\tilde \tau(1)$ with
a single split at each of the branches $e_i$ for
$i\leq k$ in such a way that the resulting train
track is splittable to $\eta$ and by putting
a mark on each of the branches
$e_{k+1},\dots,e_s$. The marked train track
$\tau(1)$ contains a subtrack $\sigma(1)$ which
can be obtained from $\sigma$ by a single split
at each of the branches $e_1,\dots,e_k$.
By our construction, the subtrack $\sigma(1)$ of $\tau(1)$
is recurrent (compare \cite{PH92}). Moreover,
since splits of a complete
train track at distinct large branches commute,
the train track $\tau(1)$ can be obtained from
$\tau$ by a splitting sequence which is induced
from a splitting sequence connecting $\sigma$
to $\sigma(1)$ as defined in the beginning of this
section. Note that $\tau(1),\sigma(1)$ only depend
on $\tau,\eta,\sigma$ but not on any choices made.

Repeat the above procedure with the train track $\tau(1)$
and the subtrack $\sigma(1)$ of $\tau(1)$.
After finitely many steps we
obtain a train track $\tau^\prime\in E(\tau,\eta)$ which
clearly satisfies
the requirements in the lemma.
\end{proof}

For a complete train track $\tau$ and
a complete geodesic lamination $\lambda$ carried by
$\tau$ define
the \emph{flat strip}
$E(\tau,\lambda)$ to be
the maximal subgraph of ${\cal T\cal T}$ whose
vertices consist of the set of all complete train
tracks which can be obtained from $\tau$ by a splitting sequence
and which carry $\lambda$.
For a convenient formulation of the following lemma, we say that
a train track $\eta$ is splittable to a complete geodesic
lamination $\lambda$ if $\eta$ carries $\lambda$.
We have.

\begin{lemma}\label{projection} Let $\tau\in {\cal V}({\cal T\cal T})$
and let ${\cal S}(\tau)\subset {\cal V}({\cal T\cal T})$ be the
set of all complete train tracks which can be obtained from
$\tau$ by a splitting sequence. Let
$E(\tau,\eta)$ be a flat
strip where either $\eta\in {\cal S}(\tau)$
or $\eta$ is a complete geodesic
lamination carried by $\tau$. Then
there is a projection
$\Pi_{E(\tau,\eta)}^1:{\cal S}(\tau)\to E(\tau,\eta)$ such
that for every $\zeta\in {\cal S}(\tau)$ the train track
$\Pi_{E(\tau,\eta)}^1(\zeta)$ is splittable to both
$\zeta,\eta$ and such that there is no $\chi\in
E(\Pi_{E(\tau,\eta)}^1(\zeta),\eta)- \Pi_{E(\tau,\eta)}^1
(\zeta)$ with this property.
\end{lemma}

\begin{proof} Let $\tau\in {\cal V}({\cal T\cal T})$ be a train
track which is splittable to a train track $\eta$. Define ${\cal
S}(\tau)\subset {\cal V}({\cal T\cal T})$ to be the set of all
train tracks $\zeta$ which can be obtained from $\tau$ by a
splitting sequence. For the proof of our lemma, we construct by
induction on the length $m$ of a splitting sequence connecting
$\tau$ to $\eta$ a projection
$\Pi_{E(\tau,\eta)}^1=\Pi_\eta^1:{\cal S}(\tau)\to
E(\tau,\eta)$. If $m=0$, i.e. if $\tau=\eta$, then we define
$\Pi^1_{\eta}(\zeta)=\tau$ for
every $\zeta\in {\cal S}(\tau)$.
By induction, assume that for some
$m\geq 1$ we determined such a projection of ${\cal S}(\tau)$ into
$E(\tau,\eta)$ for each pair $(\tau,\eta)$ with the property that
$\tau$ is splittable to $\eta$ with a splitting sequence of length
at most $m-1$. Let $\{\alpha(i)\}_{0\leq i\leq m }\subset {\cal
V}({\cal T\cal T})$ be a splitting sequence of length $m$
connecting the train track $\tau=\alpha(0)$ to $\eta=\alpha(m)$
and let $\{e_1,\dots,e_\ell\}$ be the collection of all large
branches of
$\tau$ with the property that the splitting sequence
$\{\alpha(i)\}_{0\leq i\leq m}$ contains a split at $e_i$. Note
that $\ell\geq 1$ since $m\geq 1$. For each $i$, the choice of a
right or left split at $e_i$ is determined by the requirement that
the split track carries $\eta$.

Let $\zeta\in {\cal S}(\tau)$
and assume that there is a large branch $e\in
\{e_1,\dots,e_\ell\}$ with the property that the
train track $\tilde \alpha(1)$
obtained from $\tau$ by a split at $e$ and
which is splittable to $\eta$ is also splittable
to $\zeta$.
There is then a splitting sequence
$\{\tilde \alpha(i)\}_{1\leq i\leq m}$
of length $m-1$ connecting $\tilde\alpha(1)$ to
$\tilde\alpha(m)=\alpha(m)=\eta$ (compare
Lemma 5.1 of \cite{H06a}). The flat
strip $E(\tilde \alpha(1),\eta)$
is contained in the flat strip
$E(\tau,\eta)$, and we
have $\zeta\in {\cal S}(\tilde \alpha(1))$.
By induction hypothesis, there is
a unique projection point $\Pi^1_{\tilde \alpha(1)}(\zeta)\in
E(\tilde \alpha(1),\eta)\subset
E(\tau,\eta)$ with the property that
$\Pi_{\tilde \alpha(1)}^1(\zeta)$
is splittable to $\zeta$ but that this is not
the case for any point $\rho\in
E(\Pi_{\tilde \alpha(1)}^1(\zeta),\eta)-
\Pi_{\tilde \alpha(1)}^1(\zeta)$. We define
$\Pi^1_{\eta}(\zeta)=\Pi^1_{\tilde \alpha(1)}(\zeta)$.
Then $\Pi^1_{\eta}(\zeta)$ is splittable to $\zeta$ and
this is not the case for any train track in
$E(\Pi_\eta^1(\zeta),\eta)-\Pi_\eta^1(\zeta)$. On the other
hand, a splitting sequence connecting $\tau$ to $\zeta$
is unique up to order (see Lemma 5.1 of \cite{H06a} for
a detailed discussion of this fact) and therefore
if $\xi\in E(\tau,\eta)$ is such that $\xi$ is
splittable to $\zeta$ and such that a splitting
sequence connecting $\tau$ to $\xi$ does not
contain a split at $e$, then $\xi$ contains
$e$ as a large branch, and there is a train track
$\xi^\prime$ which can be obtained from $\xi$
by a split at $e$ and which is splittable to $\zeta$.
This just means that the point
$\Pi_\eta^1(\zeta)$ does not depend on the above choice of the
large branch $e$.

If none of the train tracks $\xi\in E(\tau,\eta)$
obtained from $\tau$ by a split at one of the branches
$e_1,\dots,e_\ell$ is splittable to $\zeta$, then no train track
$\beta\in E(\tau,\eta)-\tau$ is splittable to $\zeta$ and
we define $\Pi^1_\eta(\zeta)=\tau$. This completes
our inductive construction of the map $\Pi_\eta^1:
{\cal S}(\tau)\to E(\tau,\eta)$.
Note that we have $\Pi_\eta^1(\zeta)=
\Pi_\zeta^1(\eta)$ for all $\zeta,\eta\in {\cal S}(\tau)$. Namely,
$\Pi_\eta^1(\zeta)$ is splittable to both $\zeta,\eta$, but
this is not the case for any train track which can
be obtained from $\Pi_\eta^1(\zeta)$ by a split.
This shows the lemma.
\end{proof}

Let $F$ be a
\emph{framing} for $S$ (or \emph{marking}
in the terminology of \cite{MM99}),
i.e. $F$ consists of a pants decomposition
$P$ for $S$ and a system of $3g-3+m$ \emph{spanning curves}.
The framing determines a family ${\cal P}(F)$ of
train tracks \emph{in standard form for $F$}
(\cite{PH92} and \cite{H06a}).
Let $X\subset {\cal V}({\cal T\cal T})$ be the set of all
train tracks which can be obtained from a train track in
standard form for $F$ by a
splitting sequence. By Proposition 2.2
(in the slightly more precise version which is immediate from
the proof given in \cite{H06b}), there is a
number $q>0$ such that the $q$-neighborhood of $X$ in ${\cal T\cal
T}$ is all of ${\cal T\cal T}$. Thus
if we equip $X$ with the
restriction of the metric on ${\cal T\cal T}$
then the inclusion $X\to {\cal
T\cal T}$ is a quasi-isometry.

For $\eta\in X$
there is a \emph{unique} train track $\tau$ in
standard form for $F$ which is splittable to $\eta$.
Define the
\emph{flat strip} $E(F,\eta)=E(\tau,\eta)$ to be
the maximal subgraph of ${\cal T\cal T}$ whose vertices
are the train tracks which can be obtained from $\tau$
by a splitting sequence and which are
splittable to $\eta$. If $\lambda$ is any complete geodesic
lamination then $\lambda$ is
carried by a unique train track $\tau$ in standard
form for $F$, and we write $E(F,\lambda)=E(\tau,\lambda)$.

For a complete train track $\tau$ and a complete
geodesic lamination $\lambda$ carried
by $\tau$ define a subset $A$ of
$E(\tau,\lambda)$ to be
\emph{combinatorially convex} if $A$ can be written in the
form $A=\cup_iE(\tau,\sigma_i)$ where for each $i$ we have
$\sigma_i\in E(\tau,\sigma_{i+1})$. The next result is the key to
a geometric understanding of the train track complex. For its
formulation, if $\lambda\in {\cal C\cal L}$ is a complete geodesic
lamination and if $\{\eta(i)\}_{0\leq i}$ is an
infinite splitting
sequence then we say that the sequence \emph{connects $\eta(0)$ to
$\lambda$} (or to a point in the $\delta$-neighborhood of
$\lambda$ for some $\delta >0$)
if $\cap {\cal C\cal L}(\eta(i))=\{\lambda\}$
where as before, we denote for $\eta\in {\cal V}({\cal T\cal T})$
by ${\cal C\cal L}(\eta)$ the set of all complete geodesic
laminations on $S$ which are carried by $\eta$ (we refer
to \cite{M03} for a discussion of a related construction).
We have.

\begin{proposition}\label{shortestdistance}
There is a number $\kappa>0$ with the
following property. Let $F$ be a framing for $S$ and
let $X\subset {\cal V}({\cal T\cal T})$ be the
set of all train complete train tracks which
can be obtained from a train track
in standard form for $F$ by a splitting sequence.
Then for every
$\eta\in X\cup {\cal C\cal L}$ there is a map
$\Pi_{E(F,\eta)}:X\to E(F,\eta)$ such that for every $\zeta\in X$
the following is satisfied.
\begin{enumerate}
\item There is a splitting sequence connecting a train
track $\tau^\prime$ in standard form for $F$
to $\zeta$ which passes through the $\kappa$-neighborhood
of $\Pi_{E(F,\eta)}(\zeta)$.
\item There is a splitting sequence
connecting a point in the $\kappa$-neighborhood of
$\zeta$ to a point in the $\kappa$-neighborhood of
$\eta$ which passes through
the $\kappa$-neighborhood of $\Pi_{E(F,\eta)}(\zeta)$.
\item $d(\Pi_{E(F,\eta)}(\zeta),\Pi_{E(F,\zeta)}(\eta))
\leq \kappa$ for all $\eta,\zeta\in X$.
\item For all $\lambda,\nu\in {\cal C\cal L}$ the set
$\Pi_{E(F,\lambda)}E(F,\nu)\subset E(F,\lambda)$ is
combinatorially convex.
\item $d(\Pi_{E(F,\eta)}(\zeta),\Pi_{E(F,\eta)}(\xi))\leq
\kappa d(\zeta,\xi)+\kappa$ for all $\xi,\zeta\in X$, all
$\eta\in X\cup {\cal C\cal L}$.
\item If $\eta,\zeta\in E(\tau,\lambda)$ for a train track
$\tau$ in standard form for $F$ which carries the complete
geodesic lamination
$\lambda\in {\cal C\cal L}$ then
$\Pi_{E(F,\eta)}(\zeta)=\Pi_{E(\tau,\eta)}^1(\zeta)=
\Pi_{E(\tau,\zeta)}^1(\eta)$.
\end{enumerate}
\end{proposition}

\begin{proof} Let $\tau\in {\cal V}({\cal T\cal T})$ be a complete
train track and let ${\cal S}(\tau)$ be the set of all complete
train tracks which can be obtained from $\tau$ by a splitting
sequence. Let
$\eta\in {\cal S}(\tau)$ and
let $\zeta\in {\cal S}(\tau)$. By Lemma 4.4, via
replacing $\zeta,\eta$ by their images under a splitting sequence
of uniformly bounded length we may assume that $\zeta,\eta$ do not
contain any rigid large branch. Recall from Lemma 4.8 the
definition of the ``minimal distance'' projection
$\Pi^1_{E(\tau,\eta)}=\Pi^1_{\eta}:
{\cal S}(\tau)\to E(\tau,\eta)$. The projection
point $\Pi_\eta^1(\zeta)=\zeta_1$ is uniquely determined by the
requirement that $\Pi^1_{\eta}(\zeta)$ is splittable to $\zeta$
but that no $\chi\in
E(\Pi^1_{\eta}(\zeta),\eta)-\Pi^1_{\eta}(\zeta)$ has this
property. The train track $\Pi_\eta^1(\zeta)$ determines
flat strips $E(\Pi_\eta^1(\zeta),\eta),
E(\Pi_\eta^1(\zeta),\zeta)$ which intersect in the unique
point $\Pi_\eta^1(\zeta)$.
Let ${\cal
E}(\zeta), {\cal E}(\eta)$ be the set of large branches $e$ of the
train track $\Pi^1_{\eta}(\zeta)$ with the property that a
splitting sequence connecting $\Pi^1_{\eta}(\zeta)$ to
$\zeta,\eta$ contains a split at $e$. If ${\cal
E}(\zeta)=\emptyset$ then $\zeta=\Pi_\eta^1(\zeta)\in
E(\tau,\eta)$ and we define $\Pi_{\eta}(\zeta)=\zeta$. Similarly,
if ${\cal E}(\eta)=\emptyset$ then $\Pi^1_{\eta}(\zeta)=\eta$, the
train track $\eta$ is splittable to $\zeta$ and we define
$\Pi_{\eta}(\zeta)=\eta$.

Now assume that the
sets ${\cal E}(\zeta), {\cal E}(\eta)$ are both non-empty. If
${\cal E}(\zeta)\cap {\cal E}(\eta)= \emptyset$ then define
$\Pi_\eta(\zeta)=\Pi_\eta^1(\zeta)$; note that this is in
particular the case if there is some $\lambda\in
{\cal C\cal L}$ such that
$\zeta,\eta\in E(\tau,\lambda)$. Otherwise let
$\{e_1,\dots,e_s\}={\cal E}(\zeta)\cap {\cal E}(\eta)$; then
by the definition of the map $\Pi_\eta^1$, for
each $i\leq s$ a splitting sequence connecting $\Pi^1_\eta(\zeta)$
to $\zeta$ contains a right (or left) split at the branch $e_i$
and a splitting sequence connecting $\Pi^1_{\eta}(\zeta)$ to
$\eta$ contains a left (or right) split at $e_i$. Let
$\zeta_2\in E(\Pi_\eta^1(\zeta),\zeta),
\eta_2\in E(\Pi_\eta^1(\zeta),\eta)$
be the train track obtained from
$\Pi^1_{\eta}(\zeta)=\zeta_1=\eta_1$
by a split at each of the large branches
$e_1,\dots,e_s$. Then
$\zeta_2,\eta_2$ contain a common subtrack $\hat \chi$ which is
obtained from $\Pi^1_{\eta}(\zeta)$ by a collision at each
of the large branches $e_1,\dots,e_s$, i.e. a split followed by
the removal of the diagonal of the split. Note that every geodesic
lamination which is carried by both $\zeta,\eta$ is carried by
$\hat\chi$ and that the number of branches of
$\hat\chi$ is strictly
smaller than the number of branches of $\zeta,\eta$. Moreover, by
a successive application of
Lemma 2.3.1 of \cite{PH92}, the train track $\hat\chi$ is recurrent
since the train tracks $\zeta_2,\eta_2$ are both
complete and hence recurrent. Denote by $\chi$ the
subtrack of $\zeta_2,\eta_2$ obtained from $\hat \chi$ by
removing all simple closed curve components of $\tilde \chi$.

By Lemma 4.7, there is a recurrent
splitting sequence $\{\chi(i)\}_{0\leq
i\leq p}\subset E(\tau,\eta)$
issuing from $\chi=\chi(0)$ which induces a splitting
sequence $\{\alpha(i)\}_{0\leq i\leq k}\subset 
E(\eta_2,\eta)\subset E(\tau,\eta)$ of
maximal length issuing from $\eta_2=\alpha(0)$. The train track
$\chi(p)$ is a recurrent subtrack of $\alpha(k)$, and $\chi(p)$ and
$\alpha(k)$ only depend on $\eta_2,\eta,\chi$ but not
on any choices made for the construction
of the splitting sequences. Similarly there is a recurrent
splitting
sequence $\{\tilde\chi(i)\}_{0\leq i\leq q}$ issuing from
$\chi=\tilde\chi(0)$ which induces a
splitting sequence $\{\beta(j)\}_{0\leq j\leq \ell}\subset
E(\zeta_2,\zeta)\subset E(\tau,\zeta)$
of maximal length issuing from $\zeta_2$. The pairs
of train tracks $(\chi,\chi(p))$ and $(\tilde \chi,\tilde
\chi(q))$ define flat strips $E(\chi,\chi(p)),E(\chi,\tilde
\chi(q))$. These flat strips contain all train tracks which can be
obtained from $\chi$ by a splitting sequence and which are
splittable to $\chi(p),\tilde \chi(q)$. Apply Lemma 4.8 to these
flat strips $E(\chi,\chi(p))$ and $E(\chi,\tilde \chi(q))$; this
is possible since the construction in the proof of Lemma 4.8 does
not use the assumption of completeness for our train tracks.
We find a train track $\sigma=\Pi^1_{\chi(p)}\tilde \chi(q)=
\Pi^1_{\tilde \chi(q)}\chi(p) \in E(\chi,\chi(p))\cap
E(\chi,\tilde \chi(q))$ with the property that $\sigma$ is
splittable to both $\chi(p),\tilde \chi(q)$ but that this is not
the case for any train track which can be obtained from $\sigma$
by a split. By Lemma 4.7, a splitting sequence in $E(\chi,\sigma)$
connecting $\chi$ to $\sigma$ induces splitting sequences
$\{\tilde \zeta(i)\}\subset E(\tau,\zeta), \{\tilde
\eta(j)\}\subset E(\tau,\eta)$ of maximal length issuing from
$\tilde \zeta(0)=\zeta_2,\tilde \eta(0)=\eta_2$ and connecting
$\zeta_2,\eta_2$ to train tracks $\zeta_3,\eta_3$ which contain
$\sigma$ as a subtrack and which are splittable to $\zeta,\eta$.
By Lemma 4.2, the distance between $\zeta_3,\eta_3$ in ${\cal T\cal T}$
is uniformly bounded. Moreover, a geodesic lamination which is
carried by both $\zeta_3,\eta_3$ is carried by
the union $\hat\sigma$ of $\sigma$ with the simple
closed curve components of $\hat \chi$.

Repeat this construction with the train track $\sigma$
instead of $\Pi_\eta^1(\zeta)$
and the flat strips
$E(\sigma,\chi(p))$ and $E(\sigma,\tilde\chi(q))$.
By Lemma 4.7 we obtain a recurrent
splitting sequences contained
in $E(\sigma,\chi(p)),E(\sigma,\tilde\chi(q))$
which then induce splitting sequences
in $E(\tau,\eta),E(\tau,\zeta)$.

After a
uniformly bounded number of steps we obtain a pair of train tracks
$\eta^\prime,\zeta^\prime\in {\cal V}({\cal T\cal T})$ with the
following properties.
\begin{enumerate}
\item The distance between $\eta^\prime,
\zeta^\prime$ is bounded from above by a universal constant.
\item $\zeta^\prime,\eta^\prime$ contain a common
recurrent subtrack $\beta$ (which possibly is a union of simple
closed curves) which carries
every geodesic lamination carried by both
$\zeta^\prime,\eta^\prime$.
\item For every large branch $e$ of $\beta$,
one of the following (not mutually exclusive)
possibilities holds.
\begin{enumerate}
\item[a)] One of the two train
tracks $\zeta^\prime$ or $\eta^\prime$ is not tight at $e$
and for every large proper subbranch $e^\prime$ of $e$
in $\zeta^\prime$ (or $\eta^\prime$)
the $\beta$-split of $\zeta^\prime$ (or $\eta^\prime$)
at $e^\prime$
is not splittable to $\zeta$ (or $\eta$).
\item[b)]
One of the train tracks $\zeta^\prime$ (or
$\eta^\prime$) is tight at
$e$ and no train track which can be obtained from
$\zeta^\prime$ (or $\eta^\prime$) by a single
split at $e$ is splittable to $\zeta$ (or $\eta$).
\end{enumerate}
\end{enumerate}

Let $\beta_0\subset\beta$ be the union of the simple closed curve
components of $\beta$. We claim that there is a universal number
$r>0$ with the following properties.
\begin{enumerate}
\item[a)] The train tracks $\zeta^{\prime},
\eta^{\prime}$ can be split with a splitting sequence of length at
most $r$ to train tracks $\hat\zeta, \hat\eta$ which contain
a simple multicurve $c\supset \beta_0$
as a subtrack and such that every minimal geodesic
lamination carried by both $\hat\zeta,\hat\eta$ is one of the
simple closed curve components of $c$.
\item[b)] $\hat\zeta,
\hat\eta$ are splittable to train tracks which can be obtained
from $\zeta,\eta$ by a splitting sequence of length at most $r$.
\end{enumerate}

If $\beta_0=\beta$ then there is nothing to show, so assume that
$\beta-\beta_0=\beta^\prime\not=\emptyset$. Define a
\emph{$\beta^\prime$-fake collision branch} of the train track
$\zeta^{\prime}$ to be a large branch $e$ in $\zeta^{\prime}$
which is a proper subbranch of $\beta^\prime$ and such that every
train track obtained from $\zeta^{\prime}$ by a split at $e$
contains $\beta$ as a subtrack. If $\tilde\zeta$ is obtained from
$\zeta^{\prime}$ by any split at $e$ then the number of branches of
$\tilde\zeta$ contained in $\beta^\prime$ is strictly smaller than
the number of branches of $\zeta^{\prime}$ contained in
$\beta^\prime$. By the definition of $\zeta^{\prime}$, if $e$ is a
$\beta$-fake collision branch of $\zeta^{\prime}$ then $e$ is a
branch of $\zeta$, i.e. no train track
obtained from $\zeta^\prime$ by a split at $e$ is
splittable to $\zeta$.
Namely, otherwise a splitting sequence
connecting $\zeta^{\prime}$ to $\zeta$ contains a split at $e$
which is necessarily a $\beta^\prime$-split. However, this
violates our choice of $\zeta^{\prime}$. Thus via replacing
$\zeta,\eta$ by their images under a splitting sequence whose
length does not exceed the number $q$ of branches of a complete
train track on $S$, we may assume that the train tracks
$\zeta^\prime,\eta^\prime$ do not have any $\beta^\prime$-fake
collision branches.

Let $b$ be a large branch of $\beta^\prime$ and let
$e$ be any large branch of $\zeta^{\prime}$ which is a proper
subbranch of $b$.
Note that if $\zeta^\prime$ is not tight at $b$,
such a branch always exists.
By assumption, $e$ is
not a $\beta$-fake collision branch.
Since $\beta$ is recurrent by assumption,
there is a simple closed multicurve $\nu$ which
is carried by $\beta$ and which fills $\beta$.
The train track $\zeta^\prime$ carries a complete
extension $\lambda$ of $\nu$, and the
$\beta$-split of $\zeta^\prime$ at $e$ is necessarily
a $\lambda$-split. This means that the train track
$\xi$ obtained from
$\zeta^{\prime}$ by a $\beta^\prime$-split at $e$ is
necessarily complete.
However, by the
construction of $\zeta^{\prime}$, the train track $\xi$ is
\emph{not} splittable to $\zeta$. Thus if $\xi$ is obtained from
$\zeta^{\prime}$ by say a right split (for convenience of
notation), then either the train track $\zeta^{\prime}(1)$
obtained from $\zeta^{\prime}$ by a left split at $e$ is
splittable to $\zeta$ or no train track which can be obtained from
$\zeta^{\prime}$ by a split at $e$ is splittable to $\zeta$. In
the first case we define $\zeta(1)=\zeta$. In the second case the
branch $e$ can naturally be viewed as a large branch in $\zeta$.
By assumption on $\zeta$, this branch is not rigid and the train
track $\zeta(1)$ obtained from $\zeta$ by a left split at $e$ is
complete. Now the train track $\zeta^{\prime}(1)$ obtained from
$\zeta^{\prime}$ by a left split at $e$ is splittable to
$\zeta(1)$ and therefore the train track $\zeta^{\prime}(1)$ is
complete as well.

By construction, a geodesic lamination which is carried by both
$\zeta^{\prime}(1), \eta^{\prime}$ is carried by the
largest recurrent subtrack $\beta(1)$ of $\beta$ which does
\emph{not} contain the branch $b$. In other words, the number of
branches of $\beta(1)$ is strictly smaller than the number of
branches of $\beta$. Every large proper subbranch $a$ of a
large branch
of $\beta(1)$ contained in $\zeta^{\prime}(1)$ is a large branch
in $\zeta^{\prime}$ and therefore if $a$ is contained in
$\zeta(1)$ then $a$ is not rigid. Thus we can repeat the
above
construction with the train tracks $\beta(1)$ and $\zeta(1)$. In a
number $s\geq 0$ of steps
which is bounded from above by the number $q$
of branches of a complete train track we obtain in this way a
train track $\zeta^{\prime}(s)$ containing a recurrent subtrack
$\beta(s)$ of $\beta$ as a subtrack with the additional property
that $\zeta^{\prime}(s)$ does not contain any proper subbranches
of large branches of $\beta(s)$. Moreover, a geodesic lamination
which is carried by both $\zeta^{\prime}(s)$ and $\eta$ is carried
by $\beta(s)$, and
$\beta(s)$ contains $\beta_0$ as a subtrack.
The train track $\zeta^{\prime}(s)$ is splittable
to a complete
train track $\zeta(s)$ which can be obtained from $\zeta$ by
a splitting sequence of uniformly bounded length.

If $\beta(s)$ contains components which are not simple closed
curves then $\beta(s)$ contains large branches $e_1,\dots,e_\ell$,
and each such branch is tight in $\zeta^{\prime}(s)$.
There is a number $k\leq \ell$ such that
for each $i\leq k$ a splitting sequence connecting $\zeta^{\prime}(s)$
to $\zeta(s)$ contains a split at $e_i$. Let $\zeta^{\prime}(s+1)$
be the train track which is splittable to $\zeta(s)$ and which can
be obtained from $\zeta^{\prime}(s)$ by a single split at each of
the large branches $e_i$ $(1\leq i\leq k)$. Then
$\zeta^{\prime}(s+1)$ contains a subtrack $\beta(s+1)$ which can
be obtained from $\beta(s)$ by a collision at each of the branches
$e_1,\dots,e_\ell$. Every geodesic lamination which is carried by
$\zeta^{\prime}(s+1)$ and $\eta^{\prime}$ is carried by
$\beta(s+1)$.

Repeat the above construction with the train track
$\eta^{\prime}$ and its subtrack $\beta(s)$.
We obtain a subtrack $\beta(s+1)$ of $\beta(s)$ containing $\beta_0$
and a train track $\eta^\prime(s+1)$ which can
be obtained from $\eta^\prime$ by a splitting sequence
of uniformly bounded length. Every large branch
$e$ of $\beta(s+1)$ is tight in both
$\zeta^\prime(s)$ and $\eta^\prime(s+1)$.
After finitely
many steps we obtain train tracks $\hat\zeta,\hat \eta$
which satisfy the requirements
a),b).

Now the distance between $\hat \zeta,\hat \eta$ is uniformly
bounded, and every minimal geodesic lamination which is carried by
both $\hat\zeta,\hat\eta$ is one of the simple closed curves
which form the components of $c$.
Lemma 4.6 applied to $\hat\zeta,\hat\eta$ then yields a
train track
$\tilde \zeta=\Pi_{E(F,\eta)}(\zeta)$
which satisfy the properties 1),2),3)
stated in the proposition, and property 4) follows immediately from our
construction. If $\eta$ is a complete geodesic lamination
carried by $\tau$ then choose an infinite
splitting sequence $\{\tau(i)\}$ issuing from $\tau(0)=\tau$
with $\cap {\cal C\cal L}(\tau(i))=\{\eta\}$.
By construction, for every $i$ the train track $\Pi_{E(F,\tau(i))}\zeta$
is splittable to $\Pi_{E(F,\tau(i+1))}\zeta$ and there is
a number $i_0>0$ such that $\Pi_{E(F,\tau(i))}\zeta=
\Pi_{E(F,\tau(j))}\zeta=\Pi_{E(F,\eta)}\zeta$ for
all $i,j\geq i_0$. The train track $\Pi_{E(F,\eta)}\zeta$
satisfies properties 1)-4) in the proposition.

Now let $\zeta,\eta\in X$ be arbitrary; then there are unique
train tracks $\tau,\tau^\prime$ in standard form for $F$ so that
$\tau$ is splittable to $\zeta$ and $\tau^\prime$ is splittable to
$\eta$. Let ${\cal M}(\tau),{\cal M}(\tau^\prime)$
be the set of all measured
geodesic laminations carried by $\tau,\tau^\prime$. If
${\cal M}(\tau)\cap {\cal M}(\tau^\prime)=\{0\}$ then we define
$\Pi_{E(F,\zeta)}(\eta)=\tau$ and $\Pi_{E(F,\eta)}(\zeta)=
\tau^\prime$. On the other hand, if ${\cal M}(\tau)\cap
{\cal M}(\tau^\prime)\not= \{0\}$ then $\tau,\tau^\prime$ contain a
common maximal recurrent subtrack $\chi$ which carries the support
of every lamination in ${\cal M}(\tau)\cap
{\cal M}(\tau^\prime)$ (see
Lemma 4.4 of \cite{H06b}). We apply our
above construction to the train tracks $\tau,\tau^\prime$,
which are splittable to $\zeta,\eta$ and the
common subtrack $\chi$ of $\zeta,\eta$
and extend in this way
the maps $\Pi_{E(F,\zeta)},\Pi_{E(F,\eta)}$ to all of $X$ in such
a way that properties 1)-4) stated in the lemma are satisfied.

To show property 5) in the proposition,
let $\sigma,\zeta_1,\zeta_2\in X$ and
let $\tau$ be a train track in standard form for $F$ which is
splittable to $\sigma$. Write $\chi_i=\Pi_{E(F,\sigma)}(\zeta_i)$
$(i=1,2)$ and let $\chi=\Pi_{E(F,\chi_1)}(\chi_2)$. Since
$\chi_1,\chi_2$ are both contained in the same flat strip
$E(\tau,\sigma)$, by our above construction the train track $\chi$
is splittable to both $\chi_1,\chi_2$ and there are disjoint sets
${\cal E}_1,{\cal E}_2$ of large branches of $\chi$ such that a
splitting sequence connecting $\chi$ to $\chi_i$ contains a split
at a large branch $e$ if and only if $e\in {\cal E}_i$ $(i=1,2)$.
Let $\ell_i\geq 0$ be the length of a splitting sequence
connecting $\chi$ to $\chi_1,\chi_2$ $(i=1,2)$ and let
$\ell=\max\{\ell_1,\ell_2\}$. Then the distance between
$\chi_1,\chi_2$ is not bigger than $2\ell$.

Let $\tau_1,\tau_2$ be train tracks in standard form for $F$
such that $\tau_i$ is splittable to $\zeta_i$
$(i=1,2)$. It follows from our above construction that
there is a universal constant $\kappa >0$ and there is a splitting
sequence connecting $\tau_i$ to $\zeta_i$ which passes through the
$\kappa$-neighborhood of $\chi_i$ $(i=1,2)$.
Since by Proposition 2.1 splitting
sequences are uniform quasi-geodesics, we conclude that the
distance between $\zeta_1$ and $\zeta_2$ is bounded from below by
$c \ell$ for a universal constant $c>0$. This finishes the proof
of the lemma. \end{proof}

For every $\tau\in {\cal V}({\cal T\cal T})$
and every complete geodesic
lamination $\lambda$ carried by $\tau$, the flat strip
$E(\tau,\lambda)\subset {\cal T\cal T}$ is connected and hence can
be equipped with the intrinsic metric $d_\lambda$. The following
observation is a consequence of Lemma 4.9.

\begin{corollary}\label{qiembed} There is a number $c>0$ with the
following property. For every $\tau\in {\cal V}({\cal T\cal T})$
and every complete geodesic lamination $\lambda$ carried by
$\tau$, the natural inclusion $(E(\tau,\lambda),d_\lambda)\to
{\cal T\cal T}$ is a $c$-quasi-isometric embedding.
\end{corollary}

\begin{proof} Since splitting sequences are uniform quasi-geodesics
in ${\cal T\cal T}$ which define geodesics in $E(\tau,\lambda)$
(see Lemma 5.1 of \cite{H06a}), we only have to show the
existence of a number $c>0$ with the following property. Let
$\tau\in {\cal V}({\cal T\cal T})$, let $\lambda\in {\cal C\cal
L}$ and let $\sigma_1,\sigma_2\in E(\tau,\lambda)$. Using the
notations from Lemma \ref{projection}, let
$\nu=\Pi^1_{E(\tau,\sigma_1)}(\sigma_2)=
\Pi^1_{E(\tau,\sigma_2)}(\sigma_1) \in E(\tau,\lambda)$ be the
unique train track which is splittable to $\sigma_1,\sigma_2$ and
such that no train track which can be obtained from $\nu$ by a
single split has this property. Let $\ell_1,\ell_2\geq 0$ be the
length of a splitting sequence connecting $\nu$ to
$\sigma_1,\sigma_2$; then $d(\sigma_1,\sigma_2)\geq (\ell_1
+\ell_2)/c-c$ where $d$ is the distance of ${\cal T\cal T}$.

By Proposition \ref{shortestdistance},
there is a splitting sequence connecting a point in
the $\kappa$-neigh\-bor\-hood of $\sigma_1$ to a point contained in
the $\kappa$-neighborhood of $\sigma_2$ which passes through the
$\kappa$-neighborhood of $\nu$. Now by Proposition 2.1,
splitting sequences are
$L$-quasi-geodesics in ${\cal T\cal T}$ for a universal number
$L>1$ and therefore the distance between $\sigma_1,\sigma_2$ is
not smaller than
$d(\sigma_1,\nu)/L-L-2\kappa+d(\sigma_2,\nu)/L-L-2\kappa$. On the
other hand, the distance in $E(\tau,\lambda)$ between $\sigma_1$
and $\sigma_2$ is not bigger than $L
d(\sigma_1,\nu)+Ld(\sigma_2,\nu)+2L$ from which the corollary
follows. \end{proof}

\section{The large-scale geometry of flat strips}

In this section we have a closer look at the geometry of flat
strips in ${\cal T\cal T}$. In particular, we compute the
\emph{asymptotic cone} of such a flat strip. Here a flat
strip $E(\tau,\lambda)$ is determined by a complete train track
$\tau\in {\cal V}({\cal T\cal T})$ and a complete geodesic
lamination $\lambda$ carried by $\tau$, and it
is the maximal subgraph of ${\cal T\cal T}$
whose set of vertices consists of all
train tracks $\sigma\in {\cal V}({\cal T\cal T})$ which carry
$\lambda$ and can be obtained from $\tau$ by a splitting sequence.
The flat strip $E(\tau,\lambda)$
is connected and can be equipped with the intrinsic path-metric
$d_\lambda$. By Corollary 4.10, there is a number $c>1$ not
depending on $\tau,\lambda$ such that the natural inclusion
$(E(\tau,\lambda),d_\lambda)\to {\cal T\cal T}$ is a
$c$-quasi-isometric embedding.

By Lemma 5.1 of
\cite{H06a}, there is an isometry of $(E(\tau,\lambda),d_\lambda)$
onto a connected \emph{cubical graph} in $\mathbb{R}^q$ where
$q>0$ is the number of branches of the complete train track
$\tau$. Such an isometry $\Phi$ is determined by the
choice of a point $\Phi(\tau)\in \mathbb{Z}^q$ and the
choice of a numbering of the branches of $\tau$ and has the
following property. Let $x_1,\dots,x_q$ be the standard
basis of $\mathbb{R}^q$. If $\sigma\in E(\tau,\lambda)$
is a complete train track then the numbering of the branches
of $\tau$ induces a numbering of the branches of $\sigma$.
If the branch $i$ in $\sigma$
is large, then the train track $\sigma^\prime\in E(\tau,\lambda)$
obtained from $\sigma$ by a single split at $i$ is mapped
by $\Phi$ to $\Phi(\sigma)+x_i$. We call such an isometry
$\Phi$ of $E(\tau,\lambda)$ onto the cubical graph
$\Phi(E(\tau,\lambda))\subset
\mathbb{R}^q$ \emph{standard}.

To obtain an understanding of the intrinsic geometry
of the graph $\Phi(E(\tau,\lambda))$,
consider for the moment an arbitrary
connected \emph{cubical complex} $K$ as defined on p.111-112
in \cite{BH99} which is
isometrically embedded in the euclidean space
$\mathbb{R}^q$.
Such a complex $K$ is a closed subset of $\mathbb{R}^q$
which is the union of
an at most countable number of
\emph{standard
cubes}, i.e. subsets of $\mathbb{R}^q$ which are isometric
to a cube $[0,1]^\ell$ for some $\ell\leq q$.
The intersection of any two such cubes is either
empty or is again a standard cube. If the
vertices of $K$ are points in the standard
integer lattice $\mathbb{Z}^q$ then we
call the cubical complex \emph{standard}.

Following Definition II.5.15 of \cite{BH99},
call an abstract \emph{simplicial} complex $L$ with
vertex set $V$ a \emph{flag complex} if every finite
subset $A$ of $V$ with the property that any
two distinct points in $A$ are connected
by an edge spans a simplex.
By Theorem II.5.4 and Theorem II.5.18
of \cite{BH99}, a standard cubical complex
$K$ in $\mathbb{R}^q$
has \emph{non-positive curvature} if
and only if for every vertex $v$ of $K$ the link complex
$L(v)$ of $v$ is a flag complex.
Moreover, $L(v)$ is a flag complex if and only
if $L(v)$ equipped with the path metric
induced from the round metric on the
$(q-1)$-dimensional unit sphere in $\mathbb{R}^{q}$
is a ${\rm Cat}(1)$-space.

Let again $\tau\in {\cal V}({\cal T\cal T})$ and let $\lambda\in
{\cal C\cal L}$ be a complete geodesic lamination carried by
$\tau$. Let $\Phi$ be a standard isometry of the flat strip
$E(\tau,\lambda)\subset {\cal T\cal T}$ equipped with its
intrinsic metric $d_\lambda$ onto an embedded standard cubical
graph in $\mathbb{R}^q$. Define the \emph{maximal
extension} of the graph $E(\tau,\lambda)$ to be the maximal
cubical subcomplex $C(\tau,\lambda)$ of $\mathbb{R}^q$ whose
one-skeleton equals $\Phi(E(\tau,\lambda))$. This complex is
uniquely determined by $E(\tau,\lambda)$ up to permutations of
vectors from the standard basis of $\mathbb{R}^q$ and translation
by a vector in $\mathbb{Z}^q$. In particular,
it is uniquely determined by $E(\tau,\lambda)$ up to cubical
isometry. The two-skeleton
$C^2(\tau,\lambda)$ of the complex $C(\tau,\lambda)$ is determined
as follows. Let $x_1,\dots,x_q$ be the standard basis of
$\mathbb{R}^q$; then for some $v\in \mathbb{Z}^q$ the
two-dimensional cube in $\mathbb{R}^q$ with vertices
$v,v+x_i,v+x_j,v+x_i+x_j$ is a cube in $C^2(\tau,\lambda)$ if and
only if each of its four sides is contained in
$\Phi(E(\tau,\lambda))$. For $k\geq 3$ the $k$-skeleton
$C^k(\tau,\lambda)$ of $C(\tau,\lambda)$ is constructed in the
same way by induction: If $Q$ is any $k$-cube in $\mathbb{R}^q$
all of whose sides are contained in $C^{k-1}(\tau,\lambda)$ then
we require that $Q$ is contained in $C^k(\tau,\lambda)$. We have.

\begin{lemma} The maximal extension
$C(\tau,\lambda)$ of a flat strip
$E(\tau,\lambda)\subset {\cal T\cal T}$ is
a complete ${\rm Cat}(0)$-space.
\end{lemma}

\begin{proof} We show first that the maximal extension
$C(\tau,\lambda)$ of the graph $E(\tau,\lambda)$
is of non-positive curvature. For this we
have to show that the link complex $L(v)$ of
every vertex $v$ of $C(\tau,\lambda)$ is a flag complex.

Let $x_1,\dots,x_q$ be the
standard basis of $\mathbb{R}^q$. By construction of the
map $\Phi$ (see above and compare Lemma 5.1 of \cite{H06a}), if
$v\in \Phi(E(\tau,\lambda))\cap \mathbb{Z}^q$
and if $1\leq i\leq  q$ is such that
$v+x_i\in \Phi (E(\tau,\lambda))$
then the line segment in $\mathbb{R}^q$ connecting
$v$ to $v+x_i$
is contained in $\Phi (E(\tau,\lambda))$
as well.

Now assume that $v\in \Phi(E(\tau,\lambda))\cap \mathbb{Z}^q$
and that $1\leq
i<j\leq q$ are such that $v,v+x_i,v+x_j\in \Phi(E(\tau,\lambda))$.
Let $\sigma\in E(\tau,\lambda)$ be such that $\Phi(\sigma)=v$. By
construction of the map $\Phi$, $\sigma$ is a complete train track
equipped with a numbering of its branches such that the branches
with numbers $i,j$ are large. The train track $\sigma^\prime$
obtained from $\sigma$ by the $\lambda$-split at the branch $i$ is
mapped by $\Phi$ to $v+x_i$, and the train track
$\sigma^{\prime\prime}$ obtained from $\sigma$ by the
$\lambda$-split at the branch $j$ is mapped by $\Phi$
to $v+x_j$. By definition, the line segments connecting $v$ to
$v+x_i,v+x_j$ are contained in $\Phi(E(\tau,\lambda))$. Since
$\lambda$-splits at distinct large branches in $\sigma$ commute,
the branch $j$ in the train track $\sigma^\prime$ (with respect to
the numbering inherited from the numbering of the branches of
$\sigma$) is large and the train track $\tilde \sigma$ obtained
from $\sigma^\prime$ by the $\lambda$-split at $j$ is mapped by
$\Phi$ to $v+x_i+x_j$. The same consideration
also shows that $\tilde
\sigma$ can be obtained from $\sigma^{\prime\prime}$ by the
$\lambda$-split at the branch $i$. This implies that the boundary
of the two-dimensional cube $Q$ in $\mathbb{R}^q$ with vertices
$v,v+x_i,v+x_j,v+x_i+x_j$ is contained in $C(\tau,\lambda)$ and
hence the cube $Q$ is contained in $C(\tau,\lambda)$ as well. In
other words, if $v$ is a vertex in $\Phi (E(\tau,\lambda))$ and if
the points $x_i,x_j$ (viewed as directions in the unit sphere at
$v$) are contained in the link complex $L(v)$ of $v$ then the
spherical edge connecting $x_i$ to $x_j$ is contained in $L(v)$ as
well. The obvious extension of this discussion to more than two
of the standard basis vectors $x_1,\dots,x_q$ shows the following.
If $v$ is a vertex in $E(\tau,\lambda)$, if
$k\geq 1$ and if $1\leq i_1<\dots <i_k\leq q$
are such that $v+x_{i_j}\in
\Phi(E(\tau,\lambda))$ for every $j\leq k$ then the $k$-dimensional
standard cube $Q\subset \mathbb{R}^q$ which
is determined by the vertices
$v,v+x_{i_j}$ $(j\leq k)$ is contained in $C(\tau,\lambda)$. Thus
the vertices in the link complex $L(v)$ of $v$
defined by the directions $x_{i_j}$ are pairwise joined by edges,
and their closed convex hull is a spherical simplex contained in $L(v)$.

Let $i,j\leq q$ be such that $v,v+x_i,v-x_j$ are vertices in
$\Phi(E(\tau,\lambda))$ and let $\sigma\in E(\tau,\lambda)$ be such
that $\Phi(\sigma)=v$. Then the branch $i$ in $\sigma$ is large
and the branch $j$ is small, in particular we have $i\not=j$. The
small branch $j$ can be collapsed in $\sigma$, and the branch with
number $j$ in the train track $\sigma^\prime$ obtained from
$\sigma$ by this collapse is large. There are now two
possibilities. The first possibility is that the branch $i$ in
$\sigma^\prime$ is a large branch which
is equivalent to saying that the branches $i$ and $j$ in
$\sigma$ are not incident on a common switch.
Then the train track
$\sigma^{\prime\prime}$ obtained from $\sigma^\prime$ by a
$\lambda$-split at $i$ is mapped by $\Phi$ to $v+x_i-x_j$ and
hence as above, the vertices $v-x_j,v,v+x_i-x_j, v+x_i$ are
contained in $\Phi (E(\tau,\lambda))$ and span a two-dimensional
cube in $C(\tau,\lambda)$. However, if the branch $i$ in
$\sigma^\prime$ is \emph{not} large then the two-dimensional cube
with vertices $v-x_j,v,v-x_j+x_{i},v+x_i$
is \emph{not} contained in
$C(\tau,\lambda)$ and the vertices
$x_i,-x_j$ in the link complex $L(v)$ are
\emph{not} connected by an edge. As a consequence, if
the points
$v,v-x_j,v+x_{i_1},\dots,v+x_{i_k}$ are vertices in
$\Phi(E(\tau,\lambda))$ and if $\sigma^\prime\in E(\tau,\lambda)$
is mapped by $\Phi$ to $v-x_j$ then the spherical
edge in the space of directions at $v$ connecting
the vertices $-x_j,x_{i_\ell}$ $(\ell \leq k)$
is contained in $L(v)$
if and only if the branch $i_\ell$ in
$\sigma^\prime$ is large.
This shows the following. Assume that the edges
in $\Phi(E(\tau,\eta))$ which connect $v$
to $v-x_j,v+x_{i_\ell}$ $(\ell \leq k)$, viewed as vertices in the
link complex $L(v)$ of $v$, are pairwise connected in $L(v)$ by an
edge. Then the branches $j,i_1,\dots,i_k$ in $\sigma^\prime$ are
all large and the
$k+1$-dimensional standard cube in $\mathbb{R}^q$
determined by these
vertices is contained in $C(\tau,\lambda)$. As a consequence, the
$k$-dimensional spherical simplex in $L(v)$ spanned by the
vertices $-x_j,x_{i_1},\dots,x_{i_k}$ is contained in $L(v)$
if and only if any two of these vertices are connected
by an edge.

Now consider a triple of vertices in
$\Phi( E(\tau,\lambda))$ of the
form $v,v-x_i,v-x_j$ for some $v\in \mathbb{Z}^q$. If
$\sigma,\sigma_i,\sigma_j$ are the preimages of $v,v-x_i,v-x_j$
under the map
$\Phi$ then $\sigma_i,\sigma_j$ is obtained from $\sigma$ by
a collapse of the small branch $i,j$. However, both train tracks
$\sigma_i,\sigma_j$ can be obtained from the same train track
$\tau$ by a splitting sequence. Since a splitting sequence
connecting $\tau$ to $\sigma$ is unique up to the order
of the splits (see the
discussion in the proof of Lemma 5.1 of \cite{H06a}), the branch
$i$ in $\sigma_j$ is a small branch and the train track $\eta$
which can be obtained from $\sigma_j$ by a collapse of the branch
$i$ is splittable to both $\sigma_i,\sigma_j$.
In particular,
as before the two-dimensional cube in $\mathbb{R}^q$
with vertices $v,v-x_i,v-x_j,v-x_i-x_j$ is contained in
$C(\tau,\lambda)$. The obvious extension of this consideration to
more than two of the standard basis vectors $x_1,\dots,x_q$ shows
that if $x_{i_1},\dots,x_{i_k}$ are such that
$v,v-x_{i_1},\dots,v-x_{i_k}$ are contained in $C(\tau,\lambda)$
then the same is true for the $k$-dimensional cube determined by
these points. Together with our above discussion we conclude the
following. Let $v=\Phi(\sigma)$ be a vertex of $C(\tau,\lambda)$
and let $x_{i_1},\dots,x_{i_k},x_{j_1},\dots,x_{j_\ell}$ be such
that for each $p\leq k,q\leq \ell$ the points
$v-x_{i_p},v+x_{j_q}$ are vertices of
$C(\tau,\lambda)$. If any two of the directions $x_{i_p},x_{j_q}$
defined by these vertices are connected in $L(v)$ by an edge
then for every $p\leq k,q\leq \ell$ the branch $j_q$ is large in
$\sigma$ as well as in
the train track obtained from $\sigma$ by a single collapse at
$i_p$. Equivalently, the small branch $i_p$ in $\sigma$ does not
have a switch in common with the large branch $j_q$. However, we
observed above that in this case the train track $\nu$ obtained
from $\sigma$ by a collapse of each of the branches
$i_1,\dots,i_p$ contains the large branches $j_1,\dots,j_q$ and
the cube of dimension $p+q$ determined by the vertices
$v,v-x_{i_p},v+x_{i_q}$ for $p\leq k,q\leq \ell$ is contained in
$C(\tau,\lambda)$. This shows that the link complex
$L(v)$ of $v$ is a flag complex as claimed.

Since the map $\Phi:E(\tau,\lambda)\to \mathbb{R}^q$ is proper by
construction, the cubical complex $C(\tau,\lambda)$ is a complete
geodesic metric space. Therefore to show that $C(\tau,\lambda)$
is indeed a complete ${\rm CAT}(0)$-space it is now enough to
establish that $C(\tau,\lambda)$ is simply connected. This in turn
follows if we can show that every closed edge-path in
$C(\tau,\lambda)$ which begins and ends at $\Phi(\tau)$ is
contractible. Note that via the isometry $\Phi$ such an edge-path
can be identified with a path in the graph $E(\tau,\lambda)$.

To show that this is indeed the case we proceed by induction on
the combinatorial length of the path. If this length vanishes then
the claim is trivial, so assume that the claim holds for all flat
strips $E(\sigma,\zeta)$
where $\sigma\in {\cal V}({\cal T\cal T})$ and
where $\zeta\in {\cal C\cal L}$ is carried by $\sigma$
and all closed edge-paths of
combinatorial length at most $m-1$ for some $m\geq 0$ which begin
and end at $\sigma$. Let $\gamma:[0,m]\to E(\tau,\lambda)$ be a
closed edge-path of combinatorial length $m$ beginning and ending
at $\tau$. Then $\gamma(1)$ is a train track which can be obtained
from $\gamma(0)=\tau$ by a single split at a large branch $e$. The
branch $e_0$ in $\gamma(1)$ corresponding to $e$ is small. Assume
without loss of generality that the standard isometric embedding
$\Phi:E(\tau,\lambda)\to \mathbb{R}^q$ satisfies $\Phi(\tau)=0$
and $\Phi(\gamma(1))=x_1$ where $x_1,\dots,x_q$ is the standard
basis of $\mathbb{R}^q$.
Let $\alpha_1,\dots,\alpha_q$ be the basis of $(\mathbb{R}^q)^*$
which is dual to $x_1,\dots,x_q$. Then $\alpha_1(\gamma(1))=1$ and
if $\alpha_1(\gamma(i))>0$ for some $i\in \{1,\dots,m-1\}$ then
$\gamma(i)\in E(\gamma(1),\lambda)$. In particular, if
$\alpha_1(\gamma(i))>0$ for every $i\in \{1,\dots,m-1\}$ then
$\gamma[1,m-1]$ is a loop in $E(\gamma(1),\lambda)$
beginning and ending at $\gamma(1)$ of combinatorial length
$m-2$. By our induction hypothesis, this loop
is contractible in $C(\gamma(1),\lambda)\subset
C(\tau,\lambda)$ and therefore $\gamma$ is contractible
in $C(\tau,\lambda)$.

Otherwise there is a first number
$i_0\in \{2,\dots,m-1\}$ such that
$\alpha_1(\gamma(i_0))=0$. Then $\gamma(i_0)$ can be obtained from
$\tau$ by some splitting sequence not containing a split at $e$.
In particular, $\gamma(i_0)$ contains the large branch
$e$ and $\gamma(i_0-1)$ is obtained from $\gamma(i_0)$ by a
single $\lambda$-split at $e$. Let $i_1$ be the minimum of all
numbers $i> i_0$ such that $\alpha_1(\gamma(i))>0$; if there is
not such $i$ then define $i_1=m$. Note that
$\gamma(i_1)=\gamma(i_0-1)$ if $i_1=i_0+1$.

If $i_1<m$ then $\gamma(i_1)$ can be obtained from $\gamma(i_1-1)$
by a single $\lambda$-split at $e$. Put $\tilde
\gamma(j)=\gamma(j)$ for $j\leq i_0-1$, $\tilde
\gamma(j)=\gamma(j+2)$ for $j\geq i_1-1$ and for $i_0\leq j\leq
i_1-2$ define $\tilde \gamma(j)$ to be the train track which can
be obtained from $\gamma(j+1)$ by a single $\lambda$-split at $e$.
Then the assignment $j\to \tilde \gamma(j)$ $(i_0-1\leq j\leq
i_1-2)$ determines an edge path contained in
$C(\gamma(1),\lambda)$ connecting $\tilde
\gamma(i_0-1)=\gamma(i_0-1)$ to $\tilde
\gamma(i_1-2)=\gamma(i_1)$. For every $j\in \{i_0,\dots,i_1-2\}$
the vertices $\gamma(j),\gamma(j+1), \tilde \gamma(j-1),\tilde
\gamma(j)$ are the vertices of a 2-dimensional cube embedded in
$C(\tau,\lambda)$. Thus by the definition of the maximal extension
$C(\tau,\lambda)$ of the flat strip $E(\tau,\lambda)$, this edge
path is homotopic with fixed endpoints to the edge path
$\gamma[i_0-1,i_1]$. Then $\tilde \gamma$ is homotopic
to $\gamma$ with fixed endpoints. Since the combinatorial length
of $\tilde \gamma$ equals $m-2$, by induction hypothesis the
edge-path $\tilde\gamma$ is contractible in $C(\tau,\lambda)$ and
hence the same holds true for the edge-path $\gamma$.

If $i_1=m$ then we let $\tilde \gamma(j)=\gamma(j)$
for $j\leq i_0-1$ and for $i_0\leq j\leq m-1$ define
$\tilde \gamma(j)$ to
be the unique train track which can be obtained from
$\gamma(j+1)$ by a single split at $e$. Also put
$\tilde \gamma(m)=\tau$. By the above consideration, the
loop $\tilde \gamma$ is homotopic to $\gamma$. On the other
hand, the
curve $\tilde \gamma[1,m-1]$ is a loop contained in
$E(\gamma(1),\lambda)$ of combinatorial length $m-2<m$ and
hence this loop is contractible in $C(\gamma(1),\lambda)\subset
C(\tau,\lambda)$ by induction hypothesis. But then $\tilde \gamma$
is contractible in $C(\tau,\lambda)$ and hence the same
holds true for $\gamma$. This completes the proof of our
lemma.
\end{proof}

Let
$\tau\in {\cal V}({\cal T\cal T})$ be a complete train track
which is splittable to a train track $\eta\in
{\cal V}({\cal T\cal T})$. Then the
flat strip $E(\tau,\eta)$ is defined. The
proof of Lemma 5.1 can be applied without modification
to $E(\tau,\eta)$ and shows that $E(\tau,\eta)$ admits
a natural ${\rm CAT}(0)$-cubical extension
$C(\tau,\eta)$. For every complete geodesic
lamination $\lambda$ which is carried by $\eta$,
this extension is naturally a subspace of
the extension $C(\tau,\lambda)$ of the flat
strip $E(\tau,\lambda)$. Moreover,
the maximal extension $C(\eta,\lambda)$ of the flat
strip $E(\eta,\lambda)$ is a subspace
of $C(\tau,\lambda)$ as well. We have.

\begin{lemma}\label{convex} Let $\tau\in {\cal V}({\cal T\cal T})$ and
let $\lambda$ be a complete geodesic lamination which is carried
by $\tau$. Let $\Phi:E(\tau,\lambda)\to \mathbb{R}^q$
be a standard isometric embedding and let $d$ be the
intrinsic metric on the cubical complex $C(\tau,\lambda)$.
\begin{enumerate}
\item For every vertex $\eta\in E(\tau,\lambda)$
the maximal extensions $C(\tau,\eta),C(\eta,\lambda)$
are convex subspaces of $C(\tau,\lambda)$, and
$d(C(\eta,\lambda),\Phi(\tau))=d(\Phi(\eta),\Phi(\tau))$.
\item The restriction of every coordinate function
$\alpha^i$ of $\mathbb{R}^q$ to a geodesic
ray $\gamma:[0,\infty)\to C(\tau,\lambda)$ issuing
from $\gamma(0)=\Phi(\tau)$ is non-increasing.
\end{enumerate}
\end{lemma}

\begin{proof} Let $\tau\in {\cal V}({\cal T\cal T})$,
let $\lambda$ be a complete geodesic lamination
carried by $\tau$ and let $\Phi:E(\tau,\lambda)\to
\mathbb{R}^q$ be a standard isometric embedding with
$\Phi(\tau)=0$.

Let $\eta\in {\cal V}({\cal T\cal T})$ be
a vertex in $E(\tau,\lambda)$.
We show first that $C(\eta,\lambda)\subset C(\tau,\lambda)$
is convex. For this note that if
$X$ is a any complete ${\rm Cat}(0)$-space and
if $A\subset X$ is a closed convex subset,
then $A$ is a complete ${\rm Cat}(0)$-space.
Moreover, every closed convex subset $B\subset A$ is
convex in $X$. Using this fact inductively,
we conclude that it is enough to show that
$C(\eta,\lambda)$ is convex subspace of $C(\tau,\lambda)$
for every complete
train track $\eta\in E(\tau,\lambda)$
which can be obtained from $\tau$ by a single
split.

Let as before $\alpha_1,\dots,\alpha_q$ be the
basis of $(\mathbb{R}^q)^*$ which is dual to the
standard basis $x_1,\dots,x_q$
of $\mathbb{R}^q$. Assume without
loss of generality that $\Phi(\eta)=x_1$; then
a point $z\in C(\tau,\lambda)$ is contained in
$C(\eta,\lambda)$ if and only if $\alpha_1(z)\geq 1$.
By the construction of the map $\Phi$,
for every point $(z_1,z_2,\dots,z_q)\in C(\tau,\lambda)$ with
$z_1<1$ the point $(1,z_2,\dots,z_q)$ is contained in
$C(\tau,\lambda)$ as well. Thus
the cubical complex $C(\tau,\lambda)\subset
\mathbb{R}^q$ is invariant under the natural
distance-non-increasing shortest distance projection
$\rho$ of $\mathbb{R}^q$ onto the closed
half-space $\{\alpha_1\geq 1\}$
which maps a point $z=(z_1,\dots,z_q)$ with
$z_1<1$ to $\rho(z)=(1,z_2,\dots,z_q)$.
Since $C(\tau,\lambda)$ is equipped
with the complete path metric
induced from the euclidean metric,
the restriction to $C(\tau,\lambda)$ of the
retraction $\rho$ is distance non-increasing
as well. Since the image of $C(\tau,\lambda)$ under
$\rho$ is just the cubical complex $C(\eta,\lambda)$,
the subcomplex $C(\eta,\lambda)\subset
C(\tau,\lambda)$ is convex. The same argument also
shows that $d(C(\eta,\lambda),\Phi(\tau))=
d(\Phi(\eta),\Phi(\tau))$.

To show that $C(\tau,\eta)\subset C(\tau,\lambda)$
is convex for every complete train track
$\eta\in E(\tau,\lambda)$ we argue in the same way.
Namely, by Lemma 5.1 and its proof, for each
$\eta\in E(\tau,\lambda)$ the space $C(\tau,\eta)$ is a complete
${\rm Cat}(0)$-space. Now if
$X_1\subset X_2\subset \dots$ is a nested sequence
of complete locally compact
${\rm Cat}(0)$-spaces with complete
locally compact ${\rm Cat}(0)$-union
$\cup_iX_i=X$ and if for each $i$ the space $X_i$ is a
convex subspace of $X_{i+1}$ then for each $i$, the space
$X_i$ is convex in $X$ as well. Thus as above,
it is enough to show that
$C(\tau,\eta)\subset C(\tau,\zeta)$ is convex
whenever $\zeta$ can be obtained from $\eta$
by a single split at a large branch $e$.

Assume
without loss of generality that the number of $e$ in $\eta$ with
respect to the numbering of the branches of $\tau$ defining
our standard isometry $\Phi$
equals one. Then using our above notation we have
$\alpha_1(\zeta)=\alpha_1(\eta)+1,\alpha_i(\zeta)=\alpha_i(\eta)$ for
$i\geq 2$ and therefore the distance between 
$C(\tau,\eta)$ and $\Phi(\zeta)$ with respect to
the restriction of the euclidean metric on
$\mathbb{R}^q$ equals one.
Moreover, by construction of the map
$\Phi$, for every point $(z_1,z_2,\dots,z_q)\in C(\tau,\zeta)$
with $z_1>\alpha_1(\eta)$ the point
$(\alpha_1(\eta),z_2,\dots,z_q)$ is contained
in $C(\tau,\eta)$. As above, this implies that
$C(\tau,\eta)\subset C(\tau,\zeta)$ is convex and completes
the proof of the first part of the lemma.

For the second part of the lemma it is
enough to show that for every geodesic $c:[0,b]\to
C(\tau,\lambda)$ issuing from $c(0)=\Phi(\tau)$ and
every $i\geq 0$ the function $t\to \alpha_i(c(t))$ is
non-decreasing. Namely, in this case there
is an edge-path $\rho:[0,r]\to C(\tau,\lambda)$
in the cubical complex $C(\tau,\lambda)$
whose Hausdorff distance to $c[0,b]$ is uniformly
bounded and which has the same property. However
by construction, the successive vertices met by
such an edge-path are the image under $\Phi$
of a splitting sequence in $E(\tau,\lambda)$.
However, by our above consideration, for every
$i\geq 1$ and every $s\in \mathbb{R}$ the
set $\{\alpha_i\geq s\}\cap C(\tau,\lambda)$ is convex
in $C(\tau,\lambda)$ and hence the function
$t\to \alpha_i(c(t))$ is necessarily non-decreasing.
The lemma follows.
\end{proof}

Next we justify the notion ``flat strip'' for the sets
$E(\tau,\lambda)$. Namely, recall that the natural inclusion
$(E(\tau,\lambda),d_\lambda)\to {\cal T\cal T}$ is a
quasi-isometric embedding and that $C(\tau,\lambda)$ is
quasi-isometric to its one-skeleton $E(\tau,\lambda)$. The next
lemma shows that the path-metric on $C(\tau,\lambda)$
is quasi-isometric to the restriction of the
euclidean metric.

\begin{lemma}\label{flatstripembedding} \begin{enumerate}
\item There is a number $c>0$
such that for every flat strip $E(\tau,\lambda)$
the inclusion $C(\tau,\lambda)\to \mathbb{R}^q$ is a
$c$-quasi-isometric embedding.
\item If $\lambda\in {\cal C\cal L}$ is carried by $\tau$
and if $\zeta,\eta\in E(\tau,\lambda)$ then there is a unique
train track $\Theta(\zeta,\eta)\in E(\tau,\lambda)$ such that
$\zeta,\eta\in E(\tau,\Theta(\zeta,\eta))$ and that
$\Theta(\zeta,\eta)\in E(\tau,\xi)$ for every train track $\xi\in
{\cal V}({\cal T\cal T})$ which can be obtained from both
$\zeta,\eta$ by a splitting sequence.
\end{enumerate}
\end{lemma}

\begin{proof} Let $E(\tau,\lambda)$ be a flat strip and let
$\Phi:E(\tau,\lambda)\to \mathbb{R}^q$ be a standard isometric
embedding. Since the inclusion $\iota:\Phi(E(\tau,\lambda))\to
\mathbb{R}^q$ is a one-Lipschitz map, for the first part of the
lemma it is enough to show the
existence of a universal constant $c>0$ such that for all
$\sigma,\eta\in E(\tau,\lambda)$ we have
$d_\lambda(\sigma,\eta)\leq c \Vert \Phi(\sigma)-\Phi(\eta)\Vert$
where $\Vert\,\Vert$ is the euclidean norm on $\mathbb{R}^q$ and
$d_\lambda$ is the intrinsic path metric on $E(\tau,\lambda)$.

For this let $\sigma,\eta\in E(\tau,\lambda)$. Using the notations
from Lemma \ref{projection}, write $\zeta=\Pi^1_{E(\tau,\sigma)}(\eta)=
\Pi^1_{E(\tau,\eta)}(\sigma)$. By construction, $\zeta$ is
splittable to both $\sigma,\eta$ and this is not the case for any
train track which can be obtained from $\zeta$ by a single split.
Via replacing $\Phi$ by the composition of $\Phi$ with a
translation by a vector in $\mathbb{Z}^q$ we may assume that
$\Phi(\zeta)=0$.

We claim that up to a permutation of the standard basis of
$\mathbb{R}^q$, there is a number $\ell \geq 1$ such that for the
standard direct orthogonal decomposition
$\mathbb{R}^q=\mathbb{R}^\ell\oplus \mathbb{R}^{q-\ell}$ we have
$\Phi(\sigma)\in \mathbb{R}^\ell$ and $\Phi(\eta)\in
\mathbb{R}^{q-\ell}$. Namely, by the choice of the train track
$\zeta$ and the fact that $\sigma,\eta$ both carry the complete
geodesic lamination $\lambda$, the set of large branches ${\cal
E}(\zeta)$ of $\zeta$ can be partitioned into disjoint subsets
${\cal E}^+,{\cal E}^-$ such that a splitting sequence connecting
$\zeta$ to $\sigma$ does not contain any split at a large branch
branch $e\in {\cal E}^+$ and a splitting sequence connecting
$\zeta$ to $\eta$ does not contain any split at a large branch
$e\in {\cal E}^-$.

Following \cite{PH92}, we call a trainpath $\rho:[0,m]\to \zeta$
\emph{one-sided large} if for every $i<m$ the half-branch
$\rho[i,i+1/2]$ is large and if $\rho[m-1,m]$ is a large branch. A
one-sided large trainpath $\rho:[0,m]\to \zeta$ is embedded
\cite{PH92}, and for every $i\in \{1,\dots,m-1\}$ the branch
$\rho[i-1,i]\subset \zeta$ is mixed. For every large half-branch $\hat
b$ of $\zeta$ there is a unique one-sided large trainpath issuing
from $\hat b$. Define ${\cal A}_0^+,{\cal A}_0^-$ to be the set of
all branches of $\zeta$ contained in a one-sided large trainpath
ending at a branch in ${\cal E}^{+},{\cal E}^-$. Then the sets
${\cal A}_0^+,{\cal A}_0^-$ are disjoint, and a branch of $\zeta$
is \emph{not} contained in ${\cal A}_0^+\cup {\cal A}_0^-$ if and
only if it is small. Each endpoint of a small branch is a starting
point of a one-sided large trainpath. Define ${\cal A}^{\pm}$ to
be the union of ${\cal A}_0^{\pm}$ with all small branches $b$ of
$\zeta$ with the property that both large half-branches incident
on the endpoints of $b$ are contained in ${\cal A}_0^{\pm}$. If
$b\not\in{\cal A}^+ \cup {\cal A}^-$ then $b$ is a small branch
incident on two distinct switches, and one of these switches is
the starting point of a one-sided large trainpath in ${\cal
A}_0^+$,the other is the starting point of a one-sided large
trainpath in ${\cal A}_0^-$.

The map $\Phi$ is determined
by a numbering of the branches of $\zeta$ (compare
Lemma 5.1 of \cite{H06a}).
We may assume that this numbering is such that for the cardinality
$\ell$ of ${\cal A}^-$, the set ${\cal A}^-$ consists of
the branches with numbers $1,\dots,\ell$.
A splitting sequence connecting $\zeta$ to
$\sigma$ does not contain any split at a large branch
$e\in {\cal E}^+$ by assumption.
Therefore, such a splitting sequence only
contains splits at the branches in ${\cal A}^-$.
By the choice
of our numbering, the image of any such splitting
sequence under the map
$\Phi$ is contained in the linear subspace spanned
by the first $\ell$ vectors of the standard basis of
$\mathbb{R}^q$.
Similarly, the image under $\Phi$ of a
splitting sequence connecting $\zeta$ to $\eta$ is contained
in the subspace
$\mathbb{R}^{q-\ell}\subset \mathbb{R}^q$ spanned
by the last $q-\ell$ vectors of the standard basis. This shows our claim.

The image under $\Phi$ of any edge-path in $E(\tau,\lambda)$ defined
by a splitting sequence
is an edge-path in the standard cubical subgraph
${\cal G}$ of $\mathbb{R}^q$ whose vertices are the points
$\mathbb{Z}^q$ in $\mathbb{R}^q$ with integral coordinates and
whose edges are the integral translates of the line
segments connecting $0$ to the standard basis vectors. Such a path
is without backtracking, i.e. if $\alpha_1,\dots,\alpha_q$ is the
basis of $(\mathbb{R}^q)^*$ dual to the standard basis of
$\mathbb{R}^q$ then for each $i$ the restriction of the function
$\alpha_i$ to such a path is non-decreasing. As a consequence,
such a path is a geodesic in the graph
${\cal G}$ equipped with the intrinsic
path metric and hence a uniform
quasi-geodesic in $\mathbb{R}^q$. If $\gamma_\sigma,\gamma_\eta$
are such edge-paths connecting $0=\Phi(\zeta)$ to
$\Phi(\sigma),\Phi(\eta)$ induced by a splitting sequence then
$\gamma_\sigma\subset \mathbb{R}^\ell, \gamma_\eta\subset
\mathbb{R}^{q-\ell}$ and hence $\gamma_\eta\circ
\gamma_\sigma^{-1}$ is a uniform quasi-geodesic in $\mathbb{R}^q$
connecting $\Phi(\sigma)$ to $\Phi(\eta)$. But this just means
that the distance in $\mathbb{R}^q$ between
$\Phi(\sigma),\Phi(\eta)$ is bounded from below by a universal
multiple of the distance of $\sigma,\eta$ in $E(\tau,\lambda)$ and
shows the first part of our lemma.

To show the second part, let again $\sigma,\eta\in
E(\tau,\lambda)$ and let
$\zeta=\Pi^1_{E(\tau,\eta)}(\sigma)$. Then $\zeta$
is splittable to both $\sigma,\eta$. By our above consideration, a
splitting sequence connecting $\zeta$ to $\sigma$ commutes with a
splitting sequence connecting $\zeta$ to $\eta$. Using our above
notation, if $\Phi(\zeta)=0$ then there is a train track
$\Theta(\eta,\sigma)\in E(\tau,\lambda)$ with
$\Phi(\Theta(\eta,\sigma))= \Phi(\sigma)+\Phi(\eta)$. This train
track has the property stated in the second part of the lemma.
\end{proof}

A \emph{nonprincipal ultrafilter} is a finitely additive
probability measure $\omega$ on the natural numbers $\mathbb{N}$
such that $\omega(S)=0$ or $1$ for every $S\subset \mathbb{N}$ and
$\omega(S)=0$ for every finite subset $S\subset \mathbb{N}$. Given
a compact metric space $X$ and a sequence $(a_i)\subset X$ $(i\in
\mathbb{N})$, there is a unique element ${\omega}-\lim a_i\in X$
such that for every neighborhood $U$ of ${\omega}-\lim a_i$ we
have $\omega\{i\mid a_i\in U\}=1$. In particular, given any
bounded sequence $(a_i)\subset \mathbb{R}$, ${\omega}-\lim a_i $
is a point selected by $\omega$.

Let $(X,d)$ be any metric space and let $(z_i)\subset X$. Write
$X_\infty= \{(x_i)\in \prod_{i\in \mathbb{N}} X \mid d(x_i,z_i)/i$
is bounded$\}$. For $x=(x_i),y=(y_i)\in X_\infty$ the sequence
$d(x_i,y_i)/i$ is bounded and hence we can define $\tilde
d_\omega(x,y)=\omega-\lim d(x_i,y_i)/i$. Then $\tilde d_\omega$ is
a pseudodistance on $X_\infty$, and the quotient metric space
$X_\omega$ equipped with the projection $d_\omega$ of the
pseudodistance $\tilde d_\omega$ is called the \emph{asymptotic
cone} of $X$ with respect to the non-principal ultrafilter
$\omega$ and with basepoint defined by the sequence $(z_i)$. If
$z_i=x_0$ for all $i$ and some fixed $x_0\in X$ then we denote
this basepoint by $*$. Note that
neither the asymptotic cone defined by $X$ and the constant
sequence $(x_0)$ nor the basepoint $*$ depend on the
choice of $x_0\in X$. In the sequel we always assume
that the basepoint in the construction of an asymptotic
cone of a metric space $X$ is defined by a constant
sequence unless explicitly stated otherwise.
The cone $(X_\omega,*)$ with basepoint $*$
may depend on the choice of $\omega$. If the isometry
group of $X$ acts cocompactly then an asymptotic cone with respect
to the ultrafilter $\omega$ admits a transitive group of
isometries whose elements can be represented by sequences in ${\rm
Iso}(X)$. The asymptotic cone of a ${\rm CAT}(0)$-space is a ${\rm
CAT}(0)$-space. We refer to \cite{K99} for a careful discussion of
asymptotic cones of ${\rm CAT}(0)$-spaces.

Our next goal is to determine the asymptotic cone
$C(\tau,\lambda)_\omega$ with basepoint the constant sequence
$(\tau)$ of the maximal extension $C(\tau,\lambda)$ of a flat
strip $E(\tau,\lambda)\subset {\cal T\cal T}$ where $\lambda$ is a
complete geodesic lamination carried by a
train track $\tau\in {\cal V}({\cal T\cal T})$. For
this define a \emph{metric cone} over a metric space $(\partial
Y,\angle)$ to be a metric space $(Y,d)$ of the form
$Y=[0,\infty)\times \partial Y/\sim$ where the equivalence
relation $\sim$ identifies the set $\{0\}\times \partial Y$ with a
single point. The metric $d$ on $Y$ is given by
$d((a,\xi),(b,\eta))= \sqrt{a^2+b^2-2ab\cos \angle(\xi,\eta)}$.
The space $(Y,d)$ is a ${\rm Cat}(0)$-space if and only if
$(\partial Y,\angle)$ is ${\rm Cat}(1)$ \cite{BH99}. The
metric cone is a proper ${\rm Cat}(0)$-space
if $(\partial Y,\angle)$ is a compact
${\rm Cat}(1)$-space. We have.

\begin{lemma}\label{propercat}The asymptotic cone $C(\tau,\lambda)_\omega$
with respect to a non-principal ultrafilter $\omega$ of the
maximal extension $C(\tau,\lambda)$ of a flat strip
$E(\tau,\lambda)\subset {\cal T\cal T}$ is a proper
${\rm Cat}(0)$-metric cone
and does not depend on $\omega$.
\end{lemma}

\begin{proof}
Let for the moment $Y$ be an arbitrary proper
complete ${\rm Cat}(0)$-space
with basepoint $y_0$ and distance function $d$.
Then the \emph{ideal boundary}
$\partial Y$ of $Y$ is defined; equipped with the
\emph{cone topology}, $\partial Y$ is compact.
The boundary $\partial Y$ can also be equipped with the
\emph{angular metric} $\angle$; however, the topology
defined by this metric need not coincide with the
cone topology. The metric space
$(\partial Y,\angle)$ is a complete ${\rm CAT}(1)$-space
(Theorem II.9.13 of \cite{BH99}) which
may consist of uncountably many distinct
connected components; we call it
the \emph{angular boundary} of $Y$.
If $\xi_0\not=\xi_1\in
\partial Y$ are such that $\angle(\xi_0,\xi_1)<\pi$
then there is a geodesic in $\partial Y$ connecting
$\xi_0$ to $\xi_1$
(Proposition II.9.21 of \cite{BH99}).

Assume that for some sequence $\{i(j)\}$ going to infinity the
pointed ${\rm Cat}(0)$-spaces $(Y,y_0,d/i(j))$ converge as
$j\to\infty$ in the \emph{pointed Gromov Hausdorff topology} to a
locally compact pointed metric space $(Y_\infty,y_0,d_\infty)$.
Then $(Y_\infty,d_\infty)$ is a complete ${\rm Cat}(0)$-space. By
the discussion on p.38 of \cite{B95}, $(Y_\infty,d_\infty)$ is the
quotient $[0,\infty)\times (\partial Y,\angle)/\sim$ where
$\{0\}\times \partial Y$ is identified with a single point (the
basepoint $y_0$) and where the metric $d_\infty$ is defined by
$d_\infty((a,\xi),(b,\eta))=
\sqrt{a^2+b^2-2ab\cos\angle(\xi,\eta)}$. In other words,
$(Y_\infty,d_\infty)$ equals the metric cone defined by the ${\rm
Cat}(1)$-space $(\partial Y,\angle)$. Since $(Y_\infty,d_\infty)$
is locally compact, its ideal boundary equipped with the cone
topology is compact and hence the metric space $(\partial
Y,\angle)$ is compact. In particular, it consists of only finitely
many connected components. The limit space $(Y_\infty,d_\infty)$
is independent of the sequence $\{i(j)\}$ used to define it, and
it is uniquely determined up to isometry by a closed metric ball
of positive radius about the basepoint $y_0$ in $Y_\infty$.

Now let $\lambda\in {\cal C\cal L}$ be a complete geodesic
lamination and let $\tau$ be a train track which carries
$\lambda$. Let $\omega$ be a non-principal ultrafilter. Let
$\Phi:E(\tau,\lambda)\to \mathbb{R}^q$ be a standard isometric
embedding which maps $\tau$ to $\Phi(\tau)=0$ and determines the
maximal extension $C(\tau,\lambda)$ of $E(\tau,\lambda)$. Let
$(X_\omega,*)$ be the asymptotic cone of $C(\tau,\lambda)$ defined
by the non-principal ultrafilter $\omega$ whose basepoint $*$ is
given by the constant sequence $(\Phi(\tau))$. By Lemma 5.1,
$C(\tau,\lambda)$ is a complete ${\rm Cat}(0)$-space and hence the
same is true for $X_\omega$. By Lemma 5.3 the inclusion
$C(\tau,\lambda)\to \mathbb{R}^q$ is a quasi-isometric embedding
and therefore there is a natural bilipschitz embedding of
$X_\omega$ into $\mathbb{R}^q$, the asymptotic cone of
$\mathbb{R}^q$ (see e.g. \cite{KL97}). Since $X_\omega$ is
complete, the image of this embedding is a closed subset of
$\mathbb{R}^q$ and hence $X_\omega$ is proper. In particular,
$X_\omega$ is the limit of a sequence of scaled pointed metric
spaces $(C(\tau,\lambda),\Phi(\tau),\frac{1}{i(j)})$ in the
pointed Gromov-Hausdorff topology where $\{i(j)\}\subset
\mathbb{N}$ is a sequence with $\omega\{i(j)\mid j\}=1$ (see
\cite{K99}). Thus by our above observation, the asymptotic cone
$X_\omega$ is just the euclidean cone over the ideal boundary
$\partial C(\tau,\lambda)$ of $C(\tau,\lambda)$ equipped with the
angular metric $\angle$, and it does not depend on the sequence
$\{i(j)\}$. This shows the lemma.
\end{proof}

Let again $\lambda$ be a complete geodesic lamination.
Then $\lambda$ consists of a finite number
$\lambda_1,\dots,\lambda_k$ $(1\leq k\leq 3g-3+m)$
of minimal components which are connected by
a finite number of isolated leaves. The components
$\lambda_i$ are either simple closed curves or
minimal arational laminations.
If $\lambda_i$ is a
minimal arational component then
$\lambda_i$ \emph{fills} a unique bordered
connected subsurface $S_i\subset S$
of $S$. This means that $\lambda_i$ is contained in
$S_i$, and every essential simple closed curve $c$ on $S$ which has an
essential intersection with $S_i$, i.e. which is not
freely homotopic to a curve contained in $S-S_i$,
has an essential intersection with $\lambda_i$ as
well. Up to homotopy,
the subsurfaces $S_i$ of $S$ are pairwise disjoint.
We call $S_i$ the \emph{characteristic subsurface} of $S$ for
$\lambda_i$.

If we replace each boundary component of $S_i$ by a puncture then
we obtain a surface of finite type, again denoted by $S_i$, and of
negative Euler characteristic which we call the
\emph{characteristic surface} of $\lambda_i$ (recall that we
assumed that $\lambda_i$ is minimal arational).
Sometimes we do not distinguish between the
characteristic surface of $\lambda_i$ and the
characteristic subsurface of $S$ for $\lambda_i$.
The surface $S_i$
may be a four times punctured sphere or a once punctured torus.
The lamination $\lambda_i$ can be viewed as a geodesic lamination
on the surface $S_i$ which is minimal and fills $S_i$, i.e. every
complementary component of
$\lambda_i\subset S_i$ either is a topological
disc or a once punctured topological disc. Hence every train track
$\zeta$ on $S_i$ which carries $\lambda_i$ defines a flat strip
$E(\zeta,\lambda_i)$ with maximal extension $C(\zeta,\lambda_i)$;
note that this also makes sense if $S_i$ is a forth times
punctured sphere or a once punctured torus, see \cite{PH92}. For
simplicity we denote the asymptotic cone of $C(\zeta,\lambda_i)$
with respect to the non-principal ultrafilter $\omega$ and
basepoint the constant sequence $(\zeta)$ by $A(\lambda_i)$. Note
however that $A(\lambda_i)$ may depend on $\zeta$. If $\lambda_i$
is a simple closed curve then we define $A(\lambda_i)$ to be a
single ray $[0,\infty)$.

Call a complete geodesic lamination \emph{spread-out} if it
contains precisely $3g-3+m$ minimal components. Examples of
spread out geodesic laminations are complete geodesic laminations
whose minimal components form a pants decomposition of $S$. There
are also other types of spread out geodesic laminations. For
example, let $g\geq 1$ and let $P$ be a pants decomposition of $S$
containing a separating pants curve $c$ such that the surface
obtained from $S$ by cutting along $c$ is the union of a bordered
torus $S_0$ and a surface $S_1$ of genus $g-1$ with one boundary
component and $m$ punctures. The surface $S_0$ contains a pants
curve $c_0$ from the decomposition $P$ in its interior. There is a
spread-out geodesic lamination $\lambda$ on $S$ which contains the
components of the simple geodesic multi-curve
$P-c_0$ as minimal components and whose intersection with
$S_0$ is the union of a minimal arational geodesic lamination
$\lambda_0$ and two isolated leaves which connect $\lambda_0$ to the
boundary circle of $S_0$.

Given two euclidean cones $Y_1,Y_2$ with
basepoints $y_1,y_2$, the product $Y_1\times Y_2$
can be equipped with a product metric in such a way
that the resulting metric space is an euclidean
cone with basepoint $(y_1,y_2)$. With respect to
this metric, the cones $Y_1\times \{y_2\}$ and
$\{y_1\}\times Y_2$ are convex subspaces of $Y_1\times Y_2$.
Any two non-constant geodesic rays
$\gamma_1:[0,\infty)\to Y_1\times \{y_{2}\}$
and $\gamma_2:[0,\infty)\to \{y_1\}\times Y_2$
issuing from
the basepoint $\gamma_1(0)=\gamma_2(0)=
(y_1,y_2)$ bound a flat convex subspace in $Y_1\times Y_2$
which is isometric
to the closed quadrant $\{(x_1,x_2)\in \mathbb{R}^2\mid
x_i\geq 0\}$. We call $Y_1\times Y_2$
equipped with this metric the \emph{conical product}
of $Y_1$ and $Y_2$.
The angular boundary of $Y_1\times Y_2$
equals the \emph{spherical join} $Y_1*Y_2$ of $Y_1$ and $Y_2$
(see \cite{BH99} Chapter I.5).

Call the cone $Z=\{(x_1,\dots,x_{3g-3+m})\in
\mathbb{R}^{3g-3+m}\mid x_i\geq 0\}$ the \emph{standard partition
cone} of dimension $3g-3+m$; it equals the iterated conical
product of $3g-3+m$ single rays, viewed as euclidean cones over
single points. We have.

\begin{lemma}\label{partitioncone} For a complete
geodesic lamination $\lambda$ on $S$ with minimal components
$\lambda_1,\dots,\lambda_k$ and
a complete train track $\tau$ which carries
$\lambda$, the asymptotic
cone $C(\tau,\lambda)_\omega$ with basepoint $(\tau)$
equals the conical product
of $k$ metric cones which are bilipschitz
equivalent to the cones
$A(\lambda_i)$. If $\lambda$ is spread out then
$C(\tau,\lambda)_\omega$ is isometric to
a standard partition cone of
dimension $3g-3+m$.
\end{lemma}

\begin{proof} Let $\lambda$ be a complete geodesic lamination and
let $\lambda_1,\dots,\lambda_k$ be the minimal components of
$\lambda$.
Let $s\leq k$ be such that (after reordering) the components
$\lambda_1,\dots,\lambda_s$ of $\lambda$ are minimal arational and
that the components $\lambda_{s+1},\dots,\lambda_k$ are simple
closed curves.
We say that a train track $\eta$ which carries
$\lambda$ \emph{separates} $\lambda$ if $\eta$ contains disjoint
subtracks $\zeta_1,\dots,\zeta_k$ with the following property. For
each $i$, the train track $\zeta_i$ carries $\lambda_i$. If
$i\leq s$ then $\zeta_i$ is
contained in the characteristic subsurface $S_i$ of $S$
for $\lambda_i$, and complementary
components of $\zeta_i$ on $S_i$ are in one-to-one correspondence
with the complementary components of $\lambda_i$.
If $i\geq s+1$ then $\zeta_i$ is a simple closed curve.
Moreover, every large
branch of $\eta$ is a subbranch of $\cup_i\zeta_i$. We claim that
for every train track $\tau\in {\cal V}({\cal T\cal T})$ which
carries $\lambda$ there is a finite splitting sequence
$\{\tau(i)\}_{0\leq i\leq m} \subset E(\tau,\lambda)$ issuing from
$\tau(0)=\tau$ such that $\tau(m)$ separates $\lambda$.

For this we use the results of \cite{PH92}. Recall that a
collision of a train track $\eta$ at a large branch $e$ is a split
of $\eta$ at $e$ followed by the removal of the diagonal of the
split. A collision strictly decreases the number of branches of
our train track $\eta$. Moreover, the train track obtained from
$\eta$ by a collision at $e$ is a subtrack of a train track
obtained from $\eta$ by a split at $e$. A \emph{degenerate
splitting sequence} is a sequence $\{\eta(i)\}$ of train tracks
such that for every $i$ the train track $\eta(i+1)$ can be
obtained from $\eta(i)$ by a split or a collision.

For $i\leq s$ let $S_i$ be the characteristic
surface of $\lambda_i$. As in Section 2, call a train track
$\xi$ on a surface $\tilde S$ \emph{large} if the complementary
components of $\xi$ are all topological discs and
once punctured topological discs.
For each $i\leq s$ choose a large train
track $\sigma_i$ on $S_i$ which carries $\lambda_i$ and such that
there is a one-to-one correspondence between the complementary
components of $\sigma_i$ and the complementary components of
$\lambda_i$ on $S_i$. For $i>s$ let $\sigma_i$ be the train track
which is just the simple closed curve $\lambda_i$ together with
the choice of one switch. We may assume that the train tracks
$\sigma_i$ are in fact train tracks on $S$ which are pairwise
disjoint. Then $\sigma=\cup_i\sigma_i$ is a train track on $S$
which carries $\cup_i\lambda_i$. Moreover, if $\nu$ is a train
track obtained from $\sigma$ by a splitting sequence, if $\nu$
carries $\cup_i\lambda_i$ and if $\nu$ is a subtrack of a complete
train track $\eta$ which carries $\lambda$ then $\eta$ separates
$\lambda$ provided that $\eta$ does not contain any large branch
in $\eta-\nu$.

Choose a transverse measure $\mu$ on $\cup_i\lambda_i$ with full
support. By Theorem 2.3.1 of \cite{PH92} there is a degenerate
splitting sequence $\{\tilde \tau(i)\}_{0\leq i\leq \ell}$ issuing
from $\tilde \tau(0)=\tau$ with the following properties. The
train track $\tilde \tau(\ell)$ carries $\cup_i\lambda_i$ and can
be obtained from $\sigma$ by a splitting sequence. Moreover, for
each $i$, the measure $\mu$ induces a \emph{positive} transverse
measure on $\tilde \tau(i)$. In particular, the
train tracks $\tilde\tau(i)$ are recurrent.
The sequence $\{\tilde\tau(i)\}$
contains a uniformly bounded number of collisions. Say that there
is a sequence $0\leq i_1<\dots<i_p <\ell$ (where $p$ is bounded
from above by the number of branches of a complete train track on
$S$) such that for each $j\leq p$ the train track $\tilde
\tau(i_{j}+1)$ is obtained from $\tilde \tau(i_j)$ by a collision
at a large branch $e_j$ and that for $j\not\in \{i_1,\dots,i_p\}$
the train track $\tilde \tau(j+1)$ is obtained from $\tilde
\tau(j)$ by a split.

Since $\tau$ carries the complete geodesic lamination
$\lambda\supset \cup_i\lambda_i$ by assumption and since $\mu$
define a positive transverse measure on $\tilde \tau(i)$ for each
$i$, we may assume that $\tilde\tau(j)$ carries $\lambda$ for
every $j\leq i_1$ (compare the discussion in the proof of Lemma
4.3 of \cite{H06b}). The train track $\tilde \tau(i_1+1)$ is
obtained from $\tilde \tau(i_1)$ by a collision at the large
branch $e_1$. Let $\tau(i_1+1)$ be the unique train track which
carries $\lambda$ and which can be obtained from $\tilde
\tau(i_1)$ by a \emph{split} at $e_1$. Then $\tau(i_1+1)$ carries
$\lambda$ and contains $\tilde \tau(i_1+1)$ as a subtrack. As a
consequence, the splitting sequence connecting $\tilde
\tau(i_1+1)$ to $\tilde \tau(i_2)$ induces as in Section 4 a
splitting sequence connecting $\tau(i_1+1)$ to a train track
$\tau(q)$ which carries $\lambda$ and contains $\tilde
\tau(i_2)$ as a subtrack. Inductively in finitely many steps we
obtain in this way a splitting sequence $\{\tau(i)\}_{0\leq i\leq
m}\subset E(\tau,\lambda)$ with the property that $\tau(m)$
contains $\tilde \tau(\ell)$ as a subtrack. In particular,
$\tau(m)$ separates $\lambda$ provided that $\tau(m)-\tilde
\tau(\ell)$ does not contain any large branch.

Now if $e$ is a large branch of $\tau(m)$ which is not contained
in $\tilde \tau(\ell)$ then $e$ is not incident on a switch
contained in $\tilde \tau(\ell)$ and the preimage of $e$ under a
carrying map $\lambda\to \tau(m)$ does not intersect a minimal
component of $\lambda$. Therefore this preimage consists of
\emph{finitely many} arcs. Let $\hat \tau$ be the train track
obtained from $\tau(m)$ by a single $\lambda$-split at $e$.
The branch $\hat e$ of $\hat\tau$
corresponding to the branch $e$ in $\tau$ is small. Since
$\lambda$ is complete by assumption, a carrying map
$\lambda\to\hat \tau$ is surjective and hence the number of
components of the preimage in $\lambda$ of the
branch $\hat e$ of $\hat\tau$ under such a carrying
map is strictly smaller than the number of components of
the preimage in $\lambda$
of the branch $e$ of $\tau(m)$. As a consequence, after possibly
replacing $\tau(m)$ by a train track which can be obtained from
$\tau(m)$ by a finite splitting sequence we may assume that
$\tau(m)$ separates $\lambda$.

By Lemma 5.2, $C(\tau(m),\lambda)$ is a \emph{convex} subspace of
the ${\rm CAT}(0)$-space $C(\tau,\lambda)$ whose $m$-neighborhood
in $C(\tau,\lambda)$ is all of $C(\tau,\lambda)$ (compare the
discussion in the proof of Lemma 5.2). Therefore the asymptotic
cone $C(\tau(m),\lambda)_\omega$ with basepoint $(\tau(m))$ is
isometric to the asymptotic cone $C(\tau,\lambda)_\omega$ with
basepoint $(\tau)$. Thus for the purpose of our lemma we may
assume without loss of generality that $\tau$ separates $\lambda$.
In particular, for a transverse measure $\mu$ on $\cup_i\lambda_i$
with full support, the subtrack $\sigma$ of $\tau$ of all branches
of $\tau$ with \emph{positive} $\mu$-weight decomposes into $k$
connected components $\sigma_1,\dots,\sigma_k$ where $\sigma_i$
carries $\lambda_i$ for each $i$. If $i\leq k$ is such that the
component $\lambda_i$ is a simple closed curve then $\sigma_i$ is
an embedded simple closed curve in $\tau$, and if $i$ is such that
$\lambda_i$ is minimal arational then $\sigma_i$ is a large train
track on the characteristic subsurface $S_i$ of $S$ for
$\lambda_i$.

Let again $s\leq k$ be the such that for $i\leq s$ the component
$\lambda_i$ is minimal arational and that for $i>s$ the component
$\lambda_i$ is a simple closed curve. By the discussion in Section
4, for every $i\leq s$, every splitting sequence issuing from
$\sigma_i\subset S_i$ which consists of train tracks carrying
$\lambda_i$ induces a splitting sequence issuing from $\tau$ which
is contained in the flat strip $E(\tau,\lambda)$. Moreover, for
$i\not=j$ a splitting sequence in $E(\tau,\lambda)$
induced by a sequence of
$\lambda_i$-splits issuing from $\sigma_i$ commutes with a
splitting sequence in $E(\tau,\lambda)$
induced by a sequence of $\lambda_j$-splits
issuing from $\sigma_j$. Since $\tau$ separates $\lambda$ by
assumption, up to reordering and composing the isometric embedding
$\Phi:E(\tau,\lambda)\to \mathbb{R}^q$ with a translation by an
element in $\mathbb{Z}^q$, the maximal extension $C(\tau,\lambda)$
of the flat strip $E(\tau,\lambda)$ is of the form
$C(\sigma_1,\lambda_1)\times\dots \times
C(\sigma_k,\lambda_k)\subset \mathbb{R}^{n_1}\times \dots\times
\mathbb{R}^{n_k}\times \mathbb{R}^u=\mathbb{R}^q$
where for each $i\leq k$,
$C(\sigma_i,\lambda_i)$ is the convex intersection of
$C(\tau,\lambda)$ with the euclidean subspace $\mathbb{R}^{n_i}$
of $\mathbb{R}^q$ spanned by all standard basis vectors which
correspond to subbranches of $\sigma_i$ in $\tau$ and where
$\mathbb{R}^u$ is spanned by all standard basis vectors which
correspond to branches of $\tau$ not contained in any of the
subtracks $´\sigma_i$. The convex subspace
$C(\sigma_i,\lambda)$ of $C(\tau,\lambda)$
just equals the maximal extension of the subgraph of
$E(\tau,\lambda)$ of all train tracks which carry $\lambda$ and
can be obtained from $\tau$ by a splitting sequence induced by a
sequence of $\lambda_i$-splits of $\sigma_i$. In particular, the
space $C(\sigma_i,\lambda_i)$ is bilipschitz equivalent to the
maximal extension of the flat strip $E(\sigma_i,\lambda_i)$ (in
general, however, it is not isometric to this extension).
As a consequence,
the asymptotic cone $C(\tau,\lambda)_\omega$ is a
product cone of the form stated in the lemma.

We are left with showing that for a spread out
complete geodesic lamination
$\lambda$ the asymptotic cone $C(\tau,\lambda)_\omega$ is a
standard partition cone of dimension $3g-3+m$.
Thus let $\lambda\in {\cal C\cal L}$ be a complete
geodesic lamination which
contains $3g-3+m$ minimal components. We claim that
$\lambda$ contains a sublamination $Q$ which is a union of simple
closed curves dividing $S$ into pairs of pants, borderd tori with
one boundary circle and \emph{$X$-pieces}, i.e. bordered punctured
spheres of Euler characteristic $-2$. Namely, if $\lambda$ does
not contain any minimal arational component then the minimal
components of $\lambda$ consist of a collection of $3g-3+m$ simple
closed curves. In other words, these components form a pants
decomposition for $S$ and our claim is immediate. Otherwise let
$\lambda_0$ be a minimal arational component of $\lambda$ with
characteristic subsurface $S_0$ of $S$.
Then every essential simple closed
curve on $S$
which has an essential intersection with $S_0$ (i.e. which
can not be freely homotoped to a curve contained in $S-S_0$)
intersects $\lambda_0$ transversely. A boundary component of $S_0$
has vanishing intersection number with $\lambda$ and hence since
the number of minimal components contained in $\lambda$ equals
$3g-3+m$, the boundary circles of $S_0$ are necessarily minimal
components of $\lambda$. Moreover, since $\lambda_0$ is the only
minimal component of $\lambda$ which intersects $S_0$, either
$S_0$ is a bordered torus with one boundary component or an
$X$-piece as claimed above.

Let $\tau\in {\cal V}({\cal T\cal T})$ be
a complete train track which carries $\lambda$.
By our above consideration, for the identification
of the asymptotic cone of $C(\tau,\lambda)$
we may assume without loss of
generality that $\tau$ separates $\lambda$. In particular, the
images of the minimal components
$\lambda_1,\dots,\lambda_{3g-3+m}$ of $\lambda$ under a carrying
map $\lambda\to \tau$ are disjoint subtracks $\sigma_i$ of $\tau$.
If $\lambda_i$ is a minimal arational component then $\sigma_i$ is
a train track contained in the interior of a bordered subsurface
$S_i$ of $S$ which either is a one-holed torus of Euler
characteristic $-1$ or a four holed sphere of Euler
characteristic $-2$ (where some of the holes may be punctures) and
whose boundary consists of simple closed embedded curves in
$\tau$. As a consequence, $\sigma_i$ consists of at most six
branches and four switches (Corollary 1.1.3 of \cite{PH92}), and
it contains a single large branch since
otherwise $S_i$ contains two disjoint simple closed
not mutually freely homotopic essential curves.
The mapping class group of the
surface $S_i$ contains the free group with two generators as a
subgroup of finite index, and an infinite splitting sequence of a
complete train track on $S_i$ corresponds to choosing an
infinite word in these generators. As a consequence, for
\emph{every} $i\in \{1,\dots,3g-3+m\}$ a sequence of
$\lambda_i$-splits issuing from $\sigma_i$
is unique, and the asymptotic
cone of $C(\sigma_i,\lambda_i)$ is just the single ray
$[0,\infty)$. A sequence of $\lambda$-splits issuing from $\tau$
then consists in choosing in each step one of the subtracks
$\tilde\sigma_i$ which
are filled by the laminations $\lambda_i$ and performing
either a $\tilde\sigma_i$-split at a proper large
subbranch of $\tilde\sigma_i$ or a split
which is induced by a $\lambda_i$-split of $\tilde\sigma_i$. Together
with the above, this shows that the asymptotic cone
$C(\tau,\lambda)_\omega$ is isometric to
the standard
$3g-3+m$-dimensional partition cone. This
completes the proof of the lemma.
\end{proof}

Call a proper complete ${\rm CAT}(0)$-metric
cone $Y$ \emph{standard} if its defining
${\rm Cat}(1)$-space $\partial Y$ is of diameter strictly smaller
than $\pi$. Since $Y$ is a ${\rm Cat}(0)$-space, this then implies
that the angular boundary $(\partial Y,\angle)$ of $Y$ is arcwise
connected \cite{BH99}. The following lemma gives additional
information on the asymptotic cones of all maximal extensions of
flat strips in ${\cal T\cal T}$. We always equip
the boundary $\partial C(\tau,\lambda)$ of
$C(\tau,\lambda)$ with the angular metric.

\begin{lemma}\label{diameter} Let $\tau\in {\cal V}({\cal T\cal T})$,
let $\lambda$ be a complete geodesic lamination which is carried
by $\tau$ and let $\omega$
be non-principal ultrafilter. Then the
asymptotic cone $C(\tau,\lambda)_\omega$ of the ${\rm
Cat}(0)$-space $C(\tau,\lambda)$ is a standard proper
${\rm CAT}(0)$ cone
with boundary $\partial C(\tau,\lambda)$
of diameter not bigger than $\pi/2$. There is a number $b\in (0,1)$
and an embedding of $\partial C(\tau,\lambda)$
onto a compact arcwise connected subset of
a spherical shell $\{x=(x_1,\dots,x_q)\in \mathbb{R}^q\mid
0\leq x_i\leq 1,b\leq \Vert x\Vert\leq 1\}$.
\end{lemma}

\begin{proof} Let $\tau\in {\cal V}({\cal T\cal T})$ and
let $\lambda$ be a complete geodesic lamination
carried by $\tau$. We show that the
diameter of the angular boundary $(\partial C(\tau,\lambda),\angle)$
of $C(\tau,\lambda)$ is at most $\pi/2$.

If $\partial C(\tau,\lambda)$ consists of a single point
then there is nothing to show, so assume that
$\partial C(\tau,\lambda)$ contains at least two points.
Let $\xi\not=\xi^\prime\in \partial C(\tau,\lambda)$ be points
such that the angle $\angle(\xi,\xi^\prime)$ between
$\xi,\xi^\prime$ is maximal. Such points
exist since by Lemma 5.4, the space
$(\partial C(\tau,\lambda),\angle)$
is compact. Let $\gamma\not=
\gamma^\prime:[0,\infty)\to C(\tau,\lambda)$ be geodesic
rays issuing from $\gamma(0)=\gamma^\prime(0)=\tau$ which define
the points $\xi,\xi^\prime$
in $\partial C(\tau,\lambda)$.

Assume that $C(\tau,\lambda)$ is
defind by a standard isometric embedding
$\Phi:E(\tau,\lambda)\to \mathbb{R}^q$. Let
$\alpha^1,\dots,\alpha^q\subset (\mathbb{R}^q)^*$ be the dual
basis of the standard basis of $\mathbb{R}^q$, i.e. the functions
$\alpha^i$ are the standard coordinate functions on
$\mathbb{R}^q$.
By Lemma 5.2, the restriction of each of the
euclidean coordinate functions $\alpha^i$ to any geodesic arc
in $C(\tau,\lambda)$ issuing from $\Phi(\tau)$ is
non-decreasing. More precisely, there is a number $p>0$ and there
is a splitting sequence
$\{\tau(i)\}$ in $E(\tau,\lambda)$ issuing from $\tau$ such
that the Hausdorff distance between $\{\Phi(\tau(i)\}$ and
$\gamma$ does not
exceed $p$.
Similarly, there is a
splitting sequence $\{\eta(i)\}\subset E(\tau,\lambda)$ such that
the Hausdorff distance between $\gamma^\prime[0,\infty)$ and
$\Phi(\{\eta(i)\})$ does not exceed $p$.

For $k>0$ let $\ell(k)\geq k,\ell^\prime(k)\geq k$
be such that the distance between $\gamma(k)$ and
$\Phi(\tau(\ell(k)))$ and the distance between
$\gamma^\prime(k)$ and $\Phi(\eta(\ell^\prime(k)))$ is bounded
from above by $p$.
Using the notations from Lemma \ref{projection}, for $k\geq
0$ define
$\zeta(k)=\Pi^1_{E(\tau,\tau(\ell(k)))}(\eta(\ell^\prime(k))).$
Then the train
track $\zeta(k)$ is splittable to both $\tau(\ell(k))$ and
$\eta(\ell^\prime(k))$
but this is not true for any train track
in $E(\zeta(k),\lambda)-\zeta(k)$.

Denote by $d$ the ${\rm Cat}(0)$-metric on $C(\tau,\lambda)$. We
claim that there is a number $\alpha\in
(0,\angle(\xi,\xi^\prime))$ such that
$\min\{d(\Phi(\zeta(k)),\Phi(\tau(\ell(k)))),
d(\Phi(\zeta(k)),\Phi(\eta(\ell^\prime(k))))\}\geq \alpha k$ for
every sufficiently large $k>0$. To show this claim, consider the
triangle $\Delta$ in $C(\tau,\lambda)$ with vertices
$\Phi(\tau),\Phi(\zeta(k)),\Phi(\tau(\ell(k)))$ and the triangle
$\Delta^\prime$ with vertices $\Phi(\tau),\Phi(\zeta(k)),
\eta(\ell^\prime((k)))$. The triangles $\Delta,\Delta^\prime$ have
a common side which consists of the geodesic arc connecting
$\Phi(\tau)$ to $\Phi(\zeta(k))$. Let $\Delta_0,\Delta_0^\prime$
be comparison triangles in the Euclidean plane; we may assume that
$\Delta_0,\Delta_0^\prime$ have a common side with vertices $A,C$
correponding to the points $\Phi(\tau),\Phi(\zeta(k))$. Let
$c:[0,b]\to C(\tau,\lambda)$ be the geodesic arc connecting
$c(0)=\Phi(\zeta(k))$ to $c(b)=\Phi(\tau(\ell(k)))$. By Lemma 5.2,
the geodesic arc $c$ is contained in the convex
subset $C(\zeta(k),\lambda)$ of $C(\tau,\lambda)$
whose distance to $\Phi(\tau)$ equals
$d(\Phi(\tau),\Phi(\zeta(k)))$. Thus by
convexity of the distance function on a
${\rm CAT}(0)$-space,
the distance between $\Phi(\tau)$ and $c(s)$ is non-decreasing with
$s$. By comparison, this implies that the angles of the triangles
$\Delta_0,\Delta_0^\prime$ at the vertex $C$ are not smaller than
$\pi/2$.

Since $C(\tau,\lambda)$ is a ${\rm Cat}(0)$-space,
there is a number $a\in
(0,\angle(\xi, \xi^\prime))$ and there is a number $t(a)
>0$ such that $d(\gamma(t),\gamma^\prime(t))\geq
at+2p$ for all $t\geq t(a)$ (compare \cite{BH99}). Thus if for some
$\epsilon >0$ and large enough $k$ the distance between
$\Phi(\zeta(k))$ and $\Phi(\tau(\ell(k)))$ is smaller than
$\epsilon ak$ then the distance between $\Phi(\zeta(k))$ and
$\Phi(\tau)$ is not smaller than $(1-a\epsilon)k$, and the
distance between $\Phi(\zeta(k))$ and $\Phi(\eta(\ell^\prime(k)))$
is at least $(1-\epsilon)ak$. Since the angle at $C$ of the
triangle $\Delta_0^\prime$ is not smaller than $\pi/2$, comparison
shows that the distance between $\Phi(\eta(\ell^\prime(k)))$ and
$\Phi(\tau)$ is not smaller than the length of the side opposite
to the right angle of an euclidean right-angled triangle
whose sides adjacent to the right angle
have length not smaller than $(1-a\epsilon)k,(a-a\epsilon)k$. Therefore
this distance is not smaller than
$k\sqrt{1+a^2-2a\epsilon-2a^2\epsilon+2a^2\epsilon^2}$ which is
strictly bigger than $k+2p$ provided that $k>0$ is sufficiently
large and $\epsilon >0$ is sufficiently small compared to $a$. But
the distance between $\Phi(\tau)$ and
$\Phi(\eta(\ell^\prime(k)))$ is at
most $k+p$ by the choice of $\eta(\ell^\prime(k))$ which is a
contradiction. This shows the existence of a number $\alpha
>0$ as claimed above.

As in the proof of Lemma 5.3, observe that the set ${\cal E}$ of
large branches of $\zeta(k)$ can be partitioned into two disjoint
subsets ${\cal E}={\cal E}^+\cup {\cal E}^-$ so that a splitting
sequence connecting $\zeta(k)$ to
$\tau(\ell(k)),\eta(\ell^\prime(k))$ does not contain any split at
a branch $e\in {\cal E}^+,e^\prime\in{\cal E}^-$. Using the
notations from the proof of Lemma 5.3, let ${\cal A}_0^+,{\cal
A}_0^-$ be the set of all branches of $\zeta(k)$ contained in a
one-sided large trainpath on $\zeta(k)$ terminating at a large
branch in ${\cal E}^+, {\cal E}^-$ and let ${\cal A}^{\pm}$
be the union of ${\cal A}_0^{\pm}$ with those small branches whose
endpoints are both starting points of a one-sided large trainpath
in ${\cal A}_0^{\pm}$. If we denote by ${\cal A}^0$ the collection
of all small branches not contained in ${\cal A}^+\cup {\cal A}^-$
then we obtain a partition of the set ${\cal A}$ of all branches
of $\zeta(k)$ into the disjoint sets ${\cal A}^+,{\cal A}^-,{\cal
A}^0$ (compare the proof of Lemma 5.3). Normalize the map $\Phi$
by a composition with a translation in such a way that
$\Phi(\zeta(k))=0$. After possibly a permutation of the standard
basis of $\mathbb{R}^q$, the partition of the branches of
$\zeta(k)$ into the disjoint sets ${\cal A}^+,{\cal A}^-,{\cal
A}^0$ determines a direct decomposition
$\mathbb{R}^q=\mathbb{R}^{q_1}\times \mathbb{R}^{q_2}
\times\mathbb{R}^{q_3}$ with $q_1>0,q_2>0$ and $q_3\geq 0$ such
that the image under $\Phi$ of a splitting sequence $\{\beta(i)\}$
connecting $\zeta(k)$ to $\tau(\ell(k))$ is contained in
$\mathbb{R}^{q_1}$, and the image under $\Phi$ of a splitting
sequence $\{\xi(i)\}$ connecting $\zeta(k)$ to
$\eta(\ell^\prime(k))$ is contained in $\mathbb{R}^{q_2}$ (with
the obvious interpretation as linear subspaces of $\mathbb{R}^q$).
Since splits at large branches in ${\cal A}^+,{\cal A}^-$ commute
and since both train tracks $\tau(\ell(k)),\eta(\ell^\prime(k))$
are contained in the flat strip $E(\tau,\lambda)$, if we denote by
$C^+$ and $C^-$ the maximal extensions of the flat strips
$E(\zeta(k),\tau(\ell(k)))$, $E(\zeta(k),\eta(\ell^\prime(k)))$,
viewed as convex subsets of $C(\tau,\lambda)\subset \mathbb{R}^q$
(see Lemma 5.2), then for every $x\in C^+$ and $y\in C^-$ we have
$x+y\in C(\tau,\lambda)$.

Let again $c:[0,b]\to
C(\tau,\lambda)$ be the geodesic connecting $c(0)=\Phi(\zeta(k))$ to
$c(b)=\Phi(\tau(\ell(k))))$ and let $c^\prime:[0,b^\prime]\to
C(\tau,\lambda)$ be the geodesic connecting
$c^\prime(0)=\Phi(\zeta(k))$
to $c^\prime(b^\prime)=\Phi(\eta(\ell^\prime(k))))$. Then
$c,c^\prime$ are curves in $\mathbb{R}^q$ which are
parametrized by arc length.
Let $a_1:[0,\infty)\to \mathbb{R}^2$ be two rays in the
euclidean plane parametrized by arc length and issuing
from $a_1(0)=a_2(0)=0$ which enclose a right angle at $0$.
We may assume that the tangent vectors of $a_i$ at $0$
are the standard basis vectors $e_1,e_2$.
Then for $s\in [0,b]$ and $t\in [0,b^\prime]$ the
line segment $\ell(s,t)$ in $\mathbb{R}^2$
connecting $a_1(s)$ to $a_2(t)$ and parametrized
proportional to arc length on $[0,1]$
can uniquely be represented
in the form $\ell(s,t)(u)=a_1(\rho_1(s,t)(u))+a_2(\rho_2(s,t)(u))$
for functions $\rho_1(s,t), \rho_2(s,t)$ on $[0,1]$
with values in
$[0,s],[0,t]$ and depending
continuously on $s,t$. By our above consideration, for
all $s\in
[0,b],t\in [0,b^\prime]$ and all $u\in [0,1]$ the point
$c(\rho_1(s,t)(u))+c^\prime(\rho_2(s,t)(u))$ is
contained in $C(\tau,\lambda)$. The
curve $u\to c(\rho_1(s,t)(u))+c^\prime(\rho_2(s,t)(u))$
connects $c(s)$ to $c^\prime(t)$, and its
length coincides with the length of the curve
$\ell(s,t)$. By comparison, this implies that
the triangle in $C(\tau,\lambda)$ with
vertices $\Phi(\zeta(k)),\Phi(\tau(\ell(k))),
\Phi(\eta(\ell^\prime(k))))$
is flat, and its angles at the vertices
$\Phi(\tau(\ell(k))), \Phi(\eta(\ell^\prime(k)))$ sum up to
$\pi/2$. By comparison and our above discussion, the distance
between $\Phi(\tau)$ and $\Phi(\tau(\ell(k))),
\Phi(\eta(\ell^\prime(k)))$ is not smaller than the distance
between $\Phi(\zeta(k))$ and $\Phi(\tau(\ell(k))),
\Phi(\eta(\ell^\prime(k)))$. Therefore if we denote by $\tilde
\Delta(k)$ a comparison triangle in the euclidean plane for the
triangle $\Delta(k)$ in $C(\tau,\lambda)$ with vertices
$\Phi(\tau),\Phi(\tau(\ell(k))),\Phi(\eta(\ell^\prime(k)))$, then
by comparison, the angle at the point corresponding to $\tau$ in
$\tilde \Delta(k)$ is not bigger than $\pi/2$. Since $k>0$ was
arbitrary and the distance between $\gamma(k),\gamma^\prime(k)$
and $\Phi(\tau(\ell(k))),\Phi(\eta(\ell^\prime(k)))$ is uniformly
bounded, by the definition of the angle between $\xi,\xi^\prime$
and the results in Chapter II.9 of \cite{BH99}, this means that
$\angle(\xi,\xi^\prime)\leq \pi/2$. However, we chose
$\xi,\xi^\prime$ in such a way that their angular distance is
maximal among all distances in the angular boundary of
$C(\tau,\lambda)$ and therefore the diameter of $\partial
C(\tau,\lambda)$ with respect to the angular metric is at most
$\pi/2$. This completes the proof of the first part of
our lemma.

To show the second part, let again
$\alpha^1,\dots,\alpha^q$
be the basis of $\mathbb{R}^q$ which is dual to the
standard basis of $\mathbb{R}^q$.
For every $z\in \partial C(\tau,\lambda)$ there
is a unique geodesic ray $\gamma_z:[0,\infty)\to
C(\tau,\lambda)$ issuing from $\gamma_z(0)=\Phi(\tau)$
which is asymptotic to $z$. For
$j\leq q$ let $\alpha_\omega^j(z)=\omega-\lim_{k\to \infty}
\alpha^j(\gamma_z(k))/k$. Since $\alpha^j(\gamma_z(k))\leq k$ for
all $k$, this limit exists. By our explicit construction,
the point
$\rho(z)=(\alpha_\omega^1(z),\dots,\alpha_\omega^q(z))$
has non-negative entries, has norm bounded in
$[b,1]$ for universal constant $b>0$ and depends
continuously on $z$. Moreover
the map $z\to \rho(z)$ is injective and hence
the assignment $\rho:z\to \rho(z)$ defines an
embedding of the boundary of $C(\tau,\lambda)$ onto
a compact path-connected subset of
the spherical shell $\{x=(x_1,\dots,x_q)\in \mathbb{R}^q\mid
x_i\in [0,1],b\leq \Vert x\Vert \leq 1\}$.
\end{proof}

Finally we are able to
estimate from above the
topological dimensions of the asymptotic cones
$C(\tau,\lambda)_\omega$.

\begin{lemma}\label{topologicaldimension}
The topological dimension of the
cones $C(\tau,\lambda)_\omega$ is bounded from
above by $3g-3+m$.
\end{lemma}

\begin{proof} Let $\tau\in {\cal V}({\cal T\cal T})$ and let
$\lambda$ be a complete geodesic lamination
carried by $\tau$.
Choose a non-principal
ultrafilter $\omega$.
We have to show that the
topological dimension of $C(\tau,\lambda)_\omega$
does not exceed $3g-3+m$.
For this note first that by Lemma 5.5, this holds true
for spread-out complete geodesic laminations.
In particular, it holds true for an exceptional
surface $S$, i.e. a one-punctured torus or a forth
punctured sphere with the obvious
interpretation of flat strips for these exceptional
surfaces (see the proof of Lemma 5.5).
By induction, we therefore may assume that
our dimension estimate is valid
for all proper subsurfaces of $S$. By Lemma 5.5, it
then also holds for every geodesic lamination which
does not contain
a minimal component which fills up $S$.

Thus let $\lambda$ be a complete geodesic lamination
which contains a minimal component $\lambda_0$
which fills up $S$.
Let $\tau\in {\cal V}({\cal T\cal T})$
by a train track which carries $\lambda$.
Let $\Phi:E(\tau,\lambda)\to \mathbb{R}^q$ be
a standard isometric embedding with
$\Phi(\tau)=0$ which defines
the maximal extension $C(\tau,\lambda)$.
By Lemma 5.5, the asmyptotic cone $C(\tau,\lambda)_\omega$ is
the metric cone over the angular boundary
$(\partial C(\tau,\lambda),\angle)$,
and by Lemma 5.6, this boundary is a compact
${\rm CAT}(1)$ geodesic metric space of diameter
at most $\pi/2$.  By Lemma 5.6, there is a number
$b>0$ such that the
ultrafilter $\omega$ defines an embedding $\rho$
of the boundary $\partial C(\tau,\lambda)$
onto a compact
connected subset $C$ of a spherical shell
$\{z=(z_1,\dots,z_q)\in \mathbb{R}^q\mid
z_i\geq 0,b\leq  \Vert z\Vert\leq 1\}$ for some
$b>0$.

We have to show that the topological dimension of the
set $C$ is at most $3g-4+m$. For this
let again $\alpha^1,\dots,\alpha^q$ be the
basis of $(\mathbb{R}^q)^*$ which is dual to the
standard basis. Let $\gamma_1,\dots,\gamma_n$
be geodesic rays in $C(\tau,\lambda)$ issuing from
$\Phi(\tau)$ with corresponding points
$z_1,\dots,z_n$ in the compact set $C$.
We choose the rays $\gamma_i$ in such a
way that for each of the points $z_j\in C$
there is a linear function $\alpha^{i_j}$ which
assumes a maximum at $z_j$. We may also assume that
there is some
$\epsilon >0$ with the property that
for $j\not=k$ the value
of $\alpha^{i_j}$ on $z_k$
is smaller than $\alpha^{i_j}(z_j)/(1+2\epsilon)$.
After reordering of the standard basis
vectors, we may assume that
$i_j=j$ for all $j\leq n$.

By Lemma 5.2 and its proof, for
each $i$ there is a splitting sequence
$\{\tau_i(j)\}$ whose image under the map
$\Phi$ is of
Hausdorff distance to $\gamma_i[0,\infty)$
bounded from above by a universal constant
$p>0$. For $k>0$ let $\ell_i(k)$
be the such that the distance between
$\Phi(\tau_i(\ell_i(k)))$ and $\gamma_i(k)$ is
at most $p$.
We may assume that for $\omega$-all $k$ we have
$\alpha^1(\Phi(\tau_1(\ell_1(k)))/
\alpha^1(\Phi(\tau_i(\ell_i(k)))\geq 1+\epsilon$ for
all $i\geq 2$.

Using the notation from Lemma 4.8,
for $i\geq 2$ let $\eta(i)=
\Pi^1_{E(\tau,\tau_1(\ell_1(k)))}\tau_i(\ell_i(k))\in
E(\tau,\tau_1(\ell_1(k)))$. As in Lemma 4.8, there
is a train track $\eta\in E(\tau,\tau_1(\ell_1(k)))$
such that for each $i$, $\eta(i)$ is splittable
to $\eta$ and that moreover if
$\zeta\in E(\tau,\tau_1(\ell_1(k)))$
is such that $\eta(i)$ is splittable to $\zeta$
for each $i$ then $\eta$ is splittable to $\zeta$.
The coordinate functions of $\Phi(\eta)$
satisfy
$\alpha^j(\Phi(\eta))=\max\{\alpha^j(\Phi(\eta(i)))\mid i\geq 1\}$.
Since $\alpha^1(\tau_1(\ell_1(k))>(1+\epsilon)
\max_{i\geq 2}\alpha^1(\tau_i(\ell_i(k)))$
for $i\geq 2$ there is a partition of the
set ${\cal E}$ of large branches of $\eta$
into disjoint sets ${\cal E}={\cal E}_1\cup {\cal E}_2$
with the property that a splitting sequence
connecting $\eta$ to $\tau_1(\ell_1(k))$ does not contain
any split at a large branch $e\in {\cal E}_2$ and that
moreover the following holds.
For $i\geq 2$ and using the
notation from Lemma 5.2, let
$\zeta_i=\Theta(\eta,\tau_i(\ell_i(k)))\in E(\tau,\lambda)$.
Then both $\eta$ and
$\tau_i(\ell_i(k)))$ are
splittable to $\zeta_i$, and
a splitting sequence
connecting $\eta$ to $\zeta_i$ does not contain
any split at a large branch in ${\cal E}_1$.
Define moreover inductively $\nu_2=\zeta_2$ and
$\nu_i=\Theta(\zeta_i,\nu_{i-1})$
for $i\geq 3$ and write $\zeta=\nu_k$.
Then each of the train tracks
$\zeta_i$ is splittable
to $\zeta$, and every train track with this
property can be obtained from $\zeta$ by a splitting
sequence.

By construction, we have $\alpha^1(\tau_1(\ell_1(k)))/
\alpha^1(\tau_i(\ell_i(k)))\geq 1+\epsilon$ for all
$i$. Now by monotonicity of the
coordinate functions on geodesic rays issuing from
$\Phi(\tau)$ we necessarily have
$\alpha^1(\tau_1(\ell(k))\to \infty$ $(k\to \infty)$
and therefore for sufficiently large $k$
the union ${\cal A}(\eta)$ of the set ${\cal A}_0(\eta)$
of all branches which
are contained in a one-sided large trainpath
on $\eta$ terminating at a branch in ${\cal E}_1$
with the set of all small branches whose endpoints
are both contained in ${\cal A}_0(\eta)$ is a subgraph of
$\eta$ which contains a simple closed curve.
The train tracks $\tau_i(\ell_i(k))$ $(i\geq 2)$ are contained
in a flat strip $E(\tau,\zeta)$ where
$\zeta$ can be obtained from $\eta$ by a splitting
sequence which does not contain a split
at any of the branches in ${\cal A}(\eta)$.
Hence a splitting sequence
connecting $\tau$ to $\zeta$ is induced from a splitting sequence
of a subtrack of $\eta$ contained in
a proper subsurface of $S$ of strictly
bigger Euler characteristic. Thus by induction hypothesis,
the images under the map $\rho$ of the rays
$\gamma_2,\dots,\gamma_n$
are contained in a compact subset of $C$ of dimension
at most $3g-5+m$. Then the dimension of the convex
hull in $\partial C(\tau,\lambda)$
with respect to the angular distance of the points
$z_1,\dots,z_n$ is at most
$3g-4+m$. Since the points of $C$ were arbitrarily chosen
with the above properties
this shows the lemma.
\end{proof}

\section{A quasi-convex bicombing of the train track complex}

A \emph{bicombing} of a metric space $(X,d)$
assigns to every pair of points
$x,y\in X$ a curve
$c_{x,y}:[0,1]\to X$ connecting $x=c_{x,y}(0)$ to
$y=c_{x,y}(1)$. The curve $c_{x,y}$ is called the
\emph{combing line} connecting $x$ to $y$.
We call the bicombing
\emph{symmetric} if $c_{x,y}(t)=c_{y,x}(1-t)$ for
all $x,y$ and all $t\in [0,1]$,
\emph{reflexive} if $c_{x,x}(t)=x$ for all $x\in X$ and
all $t\in [0,1]$ and
\emph{$L$-Lipschitz} for some $L\geq 1$
if for all $x,y\in X$ the curve
$t\to c_{x,y}(t/d(x,y))$ $(t\in [0,d(x,y)])$ is
$L$-Lipschitz. Call moreover the
bicombing \emph{$L$-quasi-convex} for some
$L>0$ if for all $x,y,x^\prime,y^\prime\in X$
and all $t>0$ we have $d(c_{x,y}(t),c_{x^\prime,y^\prime}(t))
\leq L(d(x,x^\prime)+d(y,y^\prime))+L$.
As an example, if $X$ is a ${\rm Cat}(0)$-space then any two
points can be connected by a unique geodesic
parametrized proportional to arc length, and
these geodesics define a
reflexive symmetric $1$-Lipschitz $1$-quasi-convex bicombing of $X$
which we call the \emph{geodesic bicombing}.

The purpose of this section is to construct
a reflexive symmetric $L$-Lipschitz $L$-quasi-convex
bicombing for the train track complex ${\cal T\cal T}$.
For this fix a framing $F$ of $S$ and let
$X$ be the set of \emph{all} complete
train tracks which can be obtained from a train track
in standard form for $F$
by a splitting sequence. Then $X$ is
$r$-dense in ${\cal T\cal T}$ for some $r>0$.
As a consequence, it is sufficient to construct
such a reflexive symmetric $L$-Lipschitz
$L$-quasi-convex bicombing for the set $X$ equipped
with the restriction of the metric on ${\cal T\cal T}$.

We begin with constructing for a train track $\tau$ in standard
form for $F$ and for a complete geodesic lamination $\lambda$
carried by $\tau$ a bicombing for the flat strip
$E(\tau,\lambda)$. The combing path connecting $\tau$ to a train
track $\eta\in E(\tau,\lambda)$ is obtained from a particular
splitting sequence connecting $\tau$ to $\eta$. First we
establish some suitable notations. Namely, let $\rho:[0,m]\to
\tau$ be any trainpath on $\tau$. Then for every $i\in
\{1,\dots,m-1\}$ there is a single branch of $\tau$ which is
incident on $\rho(i)$ and not contained in $\rho$. We call such a
branch a \emph{neighbor} of $\rho$ at $\rho(i)$. The switch is
called a \emph{left} switch (or a \emph{right} switch) if the
neighbor of $\rho$ at $\rho(i)$ is to the left (or to the right) of
$\rho$ with respect to the orientation of $\rho$ and the
orientation of $S$. Define a \emph{special trainpath} on a train
track $\sigma$ to be a trainpath $\rho:[0,2k-1]\to \sigma$ of
length $2k-1$ for some $k\geq 1$ with the following properties.
\begin{enumerate}
\item
$\rho[0,2k-1]$ is embedded in $\sigma$.
\item
For each $j\leq k-1$ the branch $\rho[2j,2j+1]$ is large and the
branch $\rho[2j+1,2j+2]$ is small.
\item
With respect to the orientation of $S$ and the orientation of
$\rho$, right and left switches in $\rho[1,2k-2]$ alternate.
\end{enumerate}
The left part of Figure G shows a special trainpath of length 5.
\begin{figure}[hb]
\includegraphics{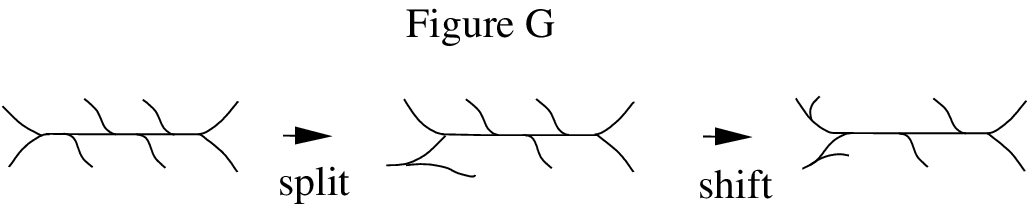}
\end{figure}

A \emph{special circle} in a train track $\sigma$ is
a trainpath $\rho:[0,2k-1]\to \sigma$ with
$\rho[0,1]=\rho[2k-2,2k-1]$ (as oriented arcs) for some
$k\geq 2$ such that $\rho[0,2k-2)$ is embedded in $\sigma$
and which satisfies the requirements 2), 3) in the
definition of a special trainpath.

If $\tau^\prime$ is any train track containing a special
trainpath or a special circle $\rho^\prime$ of length $2k-1$ and if
$\tau$ is shift equivalent to $\tau^\prime$, then
there is a natural bijection of the branches of
$\tau^\prime$ onto the branches of $\tau$, and this
bijection preserves the type of branches. In particular,
there is a trainpath $\rho$ on $\tau$ of
length at least $2k-1$ which contains the same large
and small branches as $\rho^\prime$ under our identification
of branches. We say that $\rho$ \emph{corresponds}
to $\rho^\prime$. We call a trainpath $\rho:[0,m]\to \tau$ on a
train track $\tau$ \emph{symmetric large}
if $\rho$ corresponds to a special trainpath
in this way, and we call $\rho$ a \emph{symmetric circle}
if it corresponds to a special circle. Finally, a
\emph{symmetric trainpath} is a trainpath which
either is a symmetric large trainpath or a symmetric circle.

If $\rho:[0,m]\to \tau$ is a special
trainpath on $\tau$ of length $m\geq
2$ and if $i\geq 0$ is such that $\rho[i-1,i]$ is a large branch
then there is a unique choice of a right or left split of 
$\tau$ at $\rho$ for which 
the branches $\rho[i-2,i-1]$ and $\rho[i,i+1]$ are
winners. We call such a split the \emph{$\rho$-split} of
$\tau$ at $\rho[i-1,i]$ (note that one of the two branches
$\rho[i-2,i-1]$ or $\rho[i,i+1]$ may be empty).
If the length $m$
of $\rho$ equals one, then every split of $\tau$ at $\rho[0,1]$ is
a $\rho$-split by definition. If $\tau^\prime$ is obtained from
$\tau$ by a $\rho$-split at a large branch $\rho[i-1,i]$ then
there is a natural bijection $\phi(\tau,\tau^\prime)$ from the
branches of $\tau$ onto the branches of $\tau^\prime$ (compare
the discussion in \cite{H06a}). The image of $\rho$ under the map
$\phi(\tau,\tau^\prime)$ is a trainpath $\rho^\prime$ on
$\tau^\prime$ of the same length as $\rho$. Thus for every
symmetric trainpath $\rho$ on $\tau$ it makes sense to talk about
a splitting sequence consisting of a single $\rho$-split at
\emph{every} branch of $\rho$. We first observe.

\begin{lemma}\label{combinatorics}
\begin{enumerate}
\item
Let $\tau\in {\cal V}({\cal T\cal T})$ and let $\rho:[0,m]\to
\tau$ be a symmetric large trainpath. Then there is a unique train
track $\tau^\prime$ which can be obtained from $\tau$ by a
splitting sequence of length $m$ consisting of a single
$\rho$-split at each branch $\rho[i-1,i]$ $(i\leq m)$. The train
track $\tau^\prime$ contains a symmetric large trainpath
$\rho^\prime$ of maximal length $n\in \{0,\dots, m-2\}$ which is
contained in the image of $\rho$ under $\phi(\tau,\tau^\prime)$.
If $\rho$ corresponds to a special trainpath of length one then
$\rho^\prime$ is trivial.
\item Let $\rho:[0,m]\to \tau$
be a symmetric circle. Then there is a unique train track
$\tau^\prime$ which can be obtained from $\tau$ by a splitting
sequence of length $m-1$ consisting of a single $\rho$-split at
each branch $\rho[i-1,i]$ $(i\leq m-1)$. The image of $\rho$ in
$\tau^\prime$ under the map $\phi(\tau,\tau^\prime)$ is a
symmetric circle of the same length. Moreover, there is a number
$k<m/2$ such that the train track obtained from $\tau$ by $k$ such
modifications is the image of $\tau$ under a simple Dehn twist
along the circle $\rho[0,m-1]$.
\end{enumerate}
\end{lemma}

\begin{proof}
We show the first part of the
lemma by induction of the length $m$ of a
symmetric large trainpath $\rho$ on
a complete train track $\tau$.
The case that this length equals one is trivial, so
assume that we showed the claim for all
symmetric large
trainpaths of length at most $m-1$ for some $m\geq 2$.

Let $\rho:[0,m]\to \tau$ be a symmetric large trainpath 
of length $m$ on a
complete train track $\tau$. If $\rho[0,1]$ is a
large branch then by assumption, there is a unique choice of a
right or left split of $\tau$ at $\rho[0,1]$ so that the branch
$\rho[1,2]$ is a winner of the split. Let $\tilde \tau$ be the
train track obtained from $\tau$ by this choice of a split. Then
the image of $\rho$ under the map $\phi(\tau,\tilde\tau)$ is a
trainpath in $\tilde\tau$ which begins with a small branch, and
its subpath $\phi(\tau,\tilde\tau)\rho[1,m]$ is a symmetric large
trainpath $\tilde\rho$ of length $m-1$. Clearly $\tilde\rho$ is
the maximal symmetric large trainpath on $\tau^\prime$ which is
contained in the image of $\rho$ under the map $\phi(\tau,\tilde
\tau)$. We then can apply our induction hypothesis to
$\rho^\prime$ and deduce in this way the statement of the lemma.

If the branch $\rho[0,1]$ is \emph{not} large then it is mixed and
the half-branch $\rho[0,1/2]$ is large. Thus $\rho(0)$ is the
starting point of a unique one-sided large trainpath
$\zeta:[0,k]\to \tau$ which terminates at a large branch
$\zeta[k-1,k]$. Note that we have $2\leq k\leq m$, moreover
necessarily $\zeta[0,k]=\rho[0,k]$. There is a unique choice of a
split of $\tau$ at $\rho[k-1,k]$ with the branch $\rho[k-2,k-1]$
as a winner. Let $\tau_1$ be the split track. The image of $\rho$
in $\tau_1$ under the map $\phi(\tau,\tau_1)$ is a symmetric large
trainpath $\rho_1:[0,m]\to \tau_1$ with the additional property
that the subpaths $\rho_1[0,k-1]$ and $\rho_1[k,m]$ are both
symmetric large, however the trainpath $\rho_1[k,m]$ may be
trivial. Since the lengths of these trainpaths are strictly
smaller than $m$, we can apply our induction hypothesis and obtain
the statement of the lemma for symmetric large trainpaths.

If $\rho:[0,m]\to \tau$ is a symmetric circle then we may assume
that $\rho[0,1]$ is a large branch. Let $\hat\tau$ be the train
track obtained from $\tau$ by a single $\rho$-split at
$\rho[0,1]$. The image of $\rho$ under the map
$\phi(\tau,\hat\tau)$ is a special circle $\hat\rho:[0,m]\to
\hat\tau$ of the same length as $\rho$. Moreover,
$\hat\rho[1,m-1]$ is a symmetric large trainpath. By our above
consideration, there is a splitting sequence consisting of a
single split at each branch of $\hat\rho$. The image of $\hat\rho$
under this sequence is a trainpath $\rho^\prime$ on a train track
$\tau^\prime$ beginning and ending with a small half-branch which
combines with $\phi(\tau,\tau^\prime)\hat\rho[0,1]$ to a symmetric
circle. Moreover, the map $\phi(\tau,\tau^\prime)$ of $\rho$ onto
$\rho^\prime$ preserves the type of the branches. By our explicit
construction, if $2k+1$ is the length of the special circle which
is shift equivalent to $\rho$ then after $k$ such steps we obtain
a train track $\tilde\tau$ which is the image of $\tau$ under a
single Dehn twist along the circle $\rho[0,m-1]$, with the twist
direction determined by $\rho$. This shows the lemma.
\end{proof}

We say that the train track $\tau^\prime$
as in Lemma \ref{combinatorics}
is obtained from the symmetric large trainpath
$\rho:[0,m]\to \tau$ on
$\tau$ by a \emph{level-one $\rho$-multi-split}.

Now let $\tau\in {\cal V}({\cal T\cal T})$ be
a complete train track which is splittable
to a complete train track $\sigma$.
Define a \emph{level-one splittable $\sigma$-configuration}
to be a symmetric large trainpath $\rho:[0,m]\to \tau$
of \emph{maximal length}
with the property that a level-one $\rho$-multi-split
of $\tau$ is splittable to $\sigma$.
A \emph{level-one non-splittable $\sigma$-configuration}
consists of a single large branch $e$ of
$\tau$ so that no train track obtained from
$\tau$ by a split at $e$ is splittable to $\sigma$.
We have.

\begin{lemma}\label{unique}
Let $\tau$ be a train track which is splittable to $\sigma$.
Then every large branch $e$ of $\tau$ is contained in a unique
level-one $\sigma$-configuration, and two such
level-one $\sigma$-configurations either
coincide or are disjoint.
\end{lemma}

\begin{proof}
By definition of a level-one $\sigma$-configuration
and by uniqueness of splitting sequences,
a large branch $e$ of $\tau$ such that
no split of $\tau$ at $e$ is splittable to $\sigma$
is contained in a unique level-one
$\sigma$-configuration, and
this configuration is
non-splittable.

Now assume that there is a train track $\tilde \tau$ obtained from
$\tau$ by a split at $e$ which
is splittable to $\sigma$.
Let $v$ be a switch of $\tau$
on which $e$ is incident and let $a$ be the branch
which is incident and small at $v$ and which is a
winner of the
split connecting $\tau$ to $\tilde \tau$.
Let $v^\prime$ be the second vertex on which $a$ is incident.
Assume first that $a$ is a small branch.
Then there is a unique trainpath $\rho:[0,m]\to \tau$
with $\rho[0,1]=e,\rho[1,2]=a$ and such that $\rho[2,m]$ is the
one-sided large trainpath on $\tau$ starting at $v^\prime=
\rho(2)$. Write $a^\prime=\rho[m-1,m]$.

If $a^\prime=a$ and if $\rho(1)=\rho(m-1)$ then $\rho[0,m]$ is a
loop which has a cusp at $\rho(1)$ and therefore the branch $a$ is
\emph{not} contained in the level-one $\sigma$-configuration
containing $e$. In other words, in this case the switch $v$ is
contained in the boundary of any level-one $\sigma$-configuration
containing $e$. If $a^\prime=a$ and if $\rho(1)=\rho(m)$ then
$\rho[0,m-1]$ is a $C^1$-embedded circle in $\tau$. By
construction, this circle $c$ is clearly symmetric. Thus $c$
equals the level-one $\sigma$-configuration containing $e$ if and
only if the train track obtained from $\tau$ by a $c$-multi-split
is splittable to $\sigma$. If this is not the case then the switch
$v$ is contained in the boundary of any level-one
$\sigma$-configuration containing $e$.

If $a^\prime\not=a$ then the trainpath $\rho[0,m]$ on $\tau$ is
symmetric large, and the branch $a$ is contained in a level-one
$\sigma$-configuration containing $e$ if and only if the level-one
$\rho$-multi-split of $\tau$ is splittable to $\sigma$.

Now if the branch $a$ is mixed then
the branch $\tilde a$ corresponding to $a$ in
the train track $\tilde\tau$ obtained from
$\tau$ by a single $\sigma$-split at $e$
is large. The branch $a$ belongs to a
level-one $\sigma$-configuration if and
only if the train track obtained from
$\tilde\tau$ by a $\sigma$-split at $\tilde a$
with the small branch corresponding
to $e$ in $\tilde\tau$ as a winner
is splittable to $\sigma$. Moreover,
this condition chooses uniquely one
of the two neighbors of $a$ which can possibly
be contained in a level-one
$\sigma$-configuration containing
$e$.

In finitely many steps of this form
determined by $\sigma$ we extend in this
way our trainpath starting at $e$ beyond the switch $v$ and
passing through $a$ until no further such
extension is possible. If the resulting
trainpath is not closed then we repeat this
construction with the second switch $w$ on
which $e$ is incident.
By uniqueness
of our procedure, in finitely
many steps we construct in this way a
maximal symmetric trainpath $\rho$
so that the train track obtained from $\tau$ by a
level-one $\rho$-multi-split is
splittable to $\sigma$.
\end{proof}

Now let $\tau\in {\cal V}({\cal T\cal T})$ and
let $\rho:[0,m]\to \tau$ be a symmetric large trainpath.
Let $\tau_1$ be the train track obtained from
$\tau$ by the level-one multi-split along $\rho$.
Then $\tau_1$ contains a symmetric large
trainpath $\rho_1:[0,n]\to \tau_1$ of length
$n\leq m-2$ which is contained in the image of $\rho$
under the map $\phi(\tau,\tau_1)$.
Define the \emph{level-two $\rho$-multi-split} to
be the train track $\tau_2$ obtained from $\tau_1$ by
a level-one $\rho_1$-multi-split.
Then $\tau_2$ contains a symmetric large trainpath
$\rho_2$ of length at most $m-4$ which is contained
in the image of $\rho_1$ under the natural bijection
$\phi(\tau_1,\tau_2)$. Inductively in at most $m/2$ steps
we repeat this construction until the length of
our trainpath vanishes. The train track obtained
from $\tau$ by a \emph{$\rho$-multi-split} is by
definition the train track defined inductively in this
way.

Similarly, if $\rho:[0,m]\to \tau$ is a symmetric
circle then we define the \emph{$\rho$-multi-split}
to be the train track obtained from $\tau$ by a
sequence of $\rho$-splits and which is the image
of $\tau$ under a simple Dehn twist along the circle
$\rho[0,m-1]$ as described in Lemma \ref{combinatorics}.

Now assume that $\tau$ is splittable to $\sigma$ and let
$\rho$ be a level-one $\sigma$-configuration as
defined in Lemma \ref{unique}. We define the
train track obtained from $\tau$ by a \emph{$\sigma$-move
at $\rho$} to be the unique train track $\tau^\prime$ with the
following two properties.
\begin{enumerate}
\item $\tau^\prime$ is splittable to a train track obtained
from $\tau$ by a $\rho$-multi-split.
\item If $\eta$ is splittabe to both $\sigma$ and the
train track obtained from $\tau$ by a $\rho$-multi-split
then $\eta$ is splittable to $\tau^\prime$.
\end{enumerate}

For a train track $\tau$ which is splittable to a train track
$\sigma$ we define a \emph{$\sigma$-move} to be the following
modification of $\tau$. Let $\rho_1,\dots,\rho_k$ be the
splittable level-one $\sigma$-configurations of $\tau$. By Lemma
\ref{unique}, these are uniquely defined pairwise disjoint
symmetric trainpaths on $\tau$. We define the train track
$\tau^\prime$ obtained from $\tau$ by a $\sigma$-move to be the
train track resulting from $\sigma$-moves at each of the symmetric
trainpaths $\rho_i$.

For each complete train track $\tau$ which is
splittable to a complete train track $\sigma$
define now inductively a sequence $\{\tau(i)\}_{0\leq i\leq m}\subset
E(\tau,\sigma)$ beginning at $\tau$ and ending at $\sigma$
by requiring that for each $i$ the train track
$\tau(i+1)$ is obtained from $\tau(i)$ by a
$\sigma$-move. We call the sequence the \emph{tight
multi-sequence} connecting $\tau$ to $\sigma$, and
we denote it by $\gamma(\tau,\sigma)$.
Note that a tight multi-sequence is
uniquely determined by $\tau$ and $\sigma$.
Moreover, since $\tau(i+1)$ can be obtained from $\tau(i)$ by
a non-trivial splitting sequence of uniformly bounded length,
there is a universal constant
$\kappa >0$ such that every tight multi-sequence
defines a $\kappa$-quasi-geodesic in ${\cal T\cal T}$.

We call two curves $c_1:[0,a_1]\to {\cal T\cal T},
c_2:[0,a_2]\to {\cal T\cal T}$ for some 
$0\leq a_1\leq a_2<\infty$
\emph{weight-$L$ fellow travellers} if
the following holds.
\begin{enumerate}
\item 
$d(c_1(t),c_2(t))\leq L(d(c_1(0),c_2(0))+d(c_1(a_1),c_2(a_2)))$
for every $t\in [0,a_1]$.
\item $d(c_1(a_1),c_2(t))\leq Ld(c_1(a_1),c_2(a_2))$ for
all $t\in [a_1,a_2]$.
\end{enumerate}
If $c_1,c_2$ are weight-$L$ fellow travellers then
the Hausdorff distance in ${\cal T\cal T}$ 
between the images 
$c_1[0,a_1]$ and $c_2[0,a_2]$ is bounded from above by
$L(d(c_1(0),c_2(0))+d(c_1(a_1),c_2(a_2)))$.

\begin{lemma}\label{flatconfig}
There is a number $L>0$ with the following
property.
Let $\lambda$ be a complete geodesic
lamination carried by a complete train track $\tau$
and let $\sigma,\eta\in E(\tau,\lambda)$. Then the
tight multi-sequences
$\gamma(\tau,\sigma)$ and $\gamma(\tau,\eta)$
are weight-$L$ fellow travellers.
\end{lemma}

\begin{proof}
By Corollary \ref{qiembed}, it suffices to
show the existence of a number $a>0$ with the
following property. Let $d_\lambda$ be the intrinsic
path-metric on $E(\tau,\lambda)$. Then
for $\sigma,\eta\in E(\tau,\lambda)$ the curves 
$\gamma(\tau,\sigma)$ and
$\gamma(\tau,\eta)$ are weight-$L$ fellow
travellers. Moreover,
by the explicit description of the intrinsic
distance function on $E(\tau,\lambda)$,
for our purpose it is in fact enough to
show this property for train tracks
$\eta,\sigma\in E(\tau,\lambda)$ such that
$\sigma$ can be obtained from $\eta$ by a
single split at a large branch $e$.

Thus let $\{\tau(i)\}_{0\leq i\leq k}$ be the tight multi-sequence
connecting $\tau$ to $\sigma$ and let $\{\zeta(i)\}_{0\leq i\leq
\ell}$ be the tight multi-sequence connecting $\tau$ to $\eta$;
then there is a largest number $i\leq k$ such that
$\zeta(i-1)=\tau(i-1)$. If $i-1=\ell$ then we have
$\tau(\ell)=\eta$. Since $\sigma$ can be obtained from $\eta$ by a
single split at a large branch $e$ we obtain $k=\ell+1$, and the
distance between corresponding points on
the tight multi-sequences connecting
$\tau$ to $\sigma,\eta$ 
is at most one which shows our claim.

Now consider the case that $i-1<\ell$.
Since $\zeta(i)\not=\tau(i)$
the train track $\tau(i-1)$ contains a
level-one $\sigma$-configuration
$\rho:[0,m]\to \tau(i-1)$
so that the train track $\tau_1$ obtained
from the $\sigma$-move at $\rho$ is not
splittable to $\eta$.

Assume first that the level-one $\rho$-multi-split
is not splittable to $\eta$.
Then there is some $j\leq m$ such that a splitting
sequence connecting $\tau(i-1)$ to $\eta$ does not
contain a split at the branch $\rho[j-1,j]$.
We distinguish three cases.

{\sl Case 1:} $\rho[j-1,j]$ is a large branch.

By the fact that $\sigma$ can be obtained from $\eta$ by a single
split at a large branch $e$ and uniqueness of splitting sequences,
the branch $\rho[j-1,j]$ coincides with $e$ via the map
$\phi(\tau(i-1),\eta)$. Now the branch $\rho[j,j+1]$ is incident
and small at $\rho(j)$ and hence no splitting sequence issuing
from $\tau(i-1)$ which does not contain a split at $\rho[j-1,j]$
contains a split at $\rho[j,j+1]$. This implies that a splitting
sequence connecting $\tau(i-1)$ to $\eta$ does not contain a split
at $\rho[j,j+1]$, and the same is true for a splitting sequence
connecting $\tau(i-1)$ to $\sigma$. By our assumption on $\eta,\sigma$
and by the definition of a level-one $\sigma$-configuration, we
conclude that $\rho$ is not a symmetric circle and that $m=j$. The
same argument also shows that $j=1$ and hence $\rho$ consists of a
single large branch $e$. As a consequence, the train track
$\tau(i)$ can be obtained from $\eta(i)$ by a single split at $e$.
Inductively we conclude that for every $u\in \{i,\dots,k\}$ the
train track $\tau(u)$ can be obtained from $\zeta(u)$ by a single
split at $e$. In other words, the distance between 
cooresponding points on the
tight multi-sequences $\gamma(\tau,\sigma)$ and
$\gamma(\tau,\eta)$ is at most one.

{\sl Case 2:} $\rho[j-1,j]$ is a mixed branch.

Assume without loss of generality that $\rho[j-1,j]$ is large at
$\rho(j-1)$, i.e. that the one-sided large trainpath issuing from
$\rho[j-1,j]$ is the path $\rho[j-1,q]$ for some $q\leq m$. If
$j\not=1$ then the branch $\rho[j-2,j-1]$ is small at $\phi(j-1)$
and a splitting sequence which does not contain a split at
$\rho[j-1,j]$ can not contain a split at $\rho[j-2,j-1]$. It now
follows as in Case 1 above from the definition of a level-one
$\sigma$-configuration that we necessarily have $j=1$. Then the
trainpath $\rho[1,m]$ is symmetric large and defines a level-one
$\eta$-configuration. By construction, this implies that for every
$u\in \{i,\dots,k\}$ the train track $\tau(u)$ can be obtained
from $\zeta(u)$ by a single split at the branch $e$ and hence the
distance between corresponding points on the tight multi-sequences
$\gamma(\tau,\sigma)$ and $\gamma(\tau,\eta)$ is at most one.

{\sl Case 3:} $\rho[j-1,j]$ is a small branch.

Assume first that $\rho$ is a symmetric large trainpath.
Then by definition,
we necessarily have $2\leq j\leq m-1$ and the
trainpaths $\rho[0,j-1]$ and $\rho[j,m]$ are
level-one $\eta$-configurations. It follows from
our explicit construction that for
every $u\in \{i,\dots,k\}$ the train track
$\tau(u)$ can be obtained from $\eta(u)$ by a single
split at $e$.

If $\rho$ is a symmetric circle then 
the trainpath $\zeta:[0,m-2]\to \tau(i-1)$ defined
by $\zeta[k,k+1]=\rho[k+j-1,k+j]$ (indices are taken modulo
$m-1$) is symmetric large and an $\eta$ configuration.
By our definition, the $\zeta$-multi-split is splittable
to the $\rho$-multi-split and therefore as before, 
for every $u\in \{i,\dots,k\}$ the 
the train track $\tau(u)$ can be obtained from 
$\eta(u)$ by a single split at $e$.

In the case that the level-one $\sigma$-configuration $\rho$
is also a level-one $\eta$-configuration
we can apply the above consideration
to the train tracks obtained from $\tau(i-1)$
and $\sigma(i-1)$ by the level-one $\rho$-multi-split.
Our control on the distance between corresponding points on
$\gamma(\tau,\eta)$ and $\gamma(\tau,\sigma)$ follows.
This shows the lemma.
\end{proof}

In the following proposition, we denote by
$E(F,\lambda)$ the flat strip defined by a
train track in standard form for $F$ which
carries the complete geodesic lamination $\lambda$.

\begin{proposition}\label{bicombing}
There is a number
$L>0$ with the following property.
Let $F$ be any framing of $S$ and let 
$X\subset {\cal V}({\cal T\cal T})$ to be the
set of all train tracks which can be 
obtained from a train track in standard form for
$F$ by a splitting sequence.
Then there is a reflexive symmetric
$L$-Lipschitz 
$L$-quasi-convex bicombing of $X$
equipped with the restriction of
the metric on ${\cal T\cal T}$.
If $x\in E(F,\lambda),y\in E(F,\nu)$ for some
$\lambda,\nu\in {\cal C\cal L}$ then the
combing line connecting $x$ to $y$ is contained
in the $L$-neighborhood of
$E(F,\lambda)\cup E(F,\nu)$.
\end{proposition}

\begin{proof} Let $F$ be a framing for $S$
and let $\tau$ be a complete train
track in standard form for $F$. Let $\lambda$ be
a complete geodesic lamination carried by $\tau$.
We construct first a reflexive symmetric 
$L$-Lipschitz and $L$-quasi-convex bicombing
of the flat strip $E(\tau,\lambda)$ as follows.

Let $\eta,\sigma\in E(\tau,\lambda)$.
Using the notations from Section 4, let
$\zeta=\Pi^1_{E(\tau,\eta)}(\sigma)$.
Define $\gamma$ to be the composition of the inverse
of the tight multi-sequence connecting $\zeta$ to
$\eta$ with the tight multi-sequence connecting
$\zeta$ to $\sigma$. Define
$c_{\eta,\sigma}$ to be the constant speed
reparametrization of $\gamma$ on $[0,1]$. Note that
$c_{\sigma,\eta}$ is just the inverse of
$c_{\eta,\sigma}$ and hence this defines
a symmetric reflexive $L$-Lipschitz bicombing
of $E(\tau,\lambda)$. By the results of Section 4
and Lemma \ref{flatconfig}, this bicombing is
moreover $L$-quasi-convex for a universal constant
$L>1$. If $\tau$ is a train track in standard form for
the framing $F$ then we also write $c_{F,\sigma}$ instead
of $c_{\tau,\sigma}$.

Now let $\eta\in E(F,\lambda),\beta
\in E(F,\mu)$ and write
$\zeta=\Pi_{E(F,\eta)}(\beta),\tilde \zeta=\Pi_{E(F,\beta)}(\eta)$.
We claim that 
the distance between corresponding points on the 
curves $c_{F,\zeta}$ and
$c_{F,\tilde \zeta}$ is uniformly bounded.

For this let $\sigma$ be a subtrack of a
train track $\tau$ and let $\sigma^\prime$ be
a train track obtained from $\sigma$ by
a single split at a large branch $e$. 
Using the terminology from Section 4, let $\tau^\prime$
be the train track which contains $\sigma^\prime$
as a subtrack and is obtained from $\tau$ 
as follows. Modify $\tau$ to a train track
$\tilde \tau$ obtained from $\tau$ by a $\sigma$-move
at $e$ and let $\tau^\prime$ be the train track
obtained from $\tilde \tau$ by a single split at the
tight branch $e$ and which contains $\sigma^\prime$
as a subtrack. If 
the $\sigma$-complexity $\chi(\tau,\sigma)$
of $\tau$ coincides with the 
$\sigma^\prime$-complexity $\chi(\tau^\prime,\sigma^\prime)$
of $\tau^\prime$ then the
large branch $e$ of $\sigma$ defines an embedded
trainpath in $\tau$ which is just the $\tau^\prime$-configuration
of $\tau$ as defined above. Moreover, 
$\tau^\prime$ is obtained from $\tau$ by a $\tau^\prime$-move.

Together this shows the following.
Let $\sigma$ be any birecurrent generic train track
on $S$. Assume that $\sigma$
is splittable to a train track $\sigma^\prime$.
Then we can define as above a tight multi-sequence
$\{\sigma(i)\}_{0\leq i\leq p}$ connecting 
$\sigma=\sigma(0)$ to a train track $\sigma^\prime=\sigma(p)$.
Let $\tau$ be a complete train 
track which contains $\sigma$ as a subtrack and let
$\tau^\prime$ be obtained from $\tau$ by a splitting
sequence induced from a splitting sequence connecting
$\sigma$ to $\sigma^\prime$. Assume that $\chi(\tau,\sigma)=
\chi(\tau^\prime,\sigma^\prime)$ and let $\{\tau(j)\}$ be the
tight multi-sequence connecting $\tau$ to $\tau^\prime$. Then for
every $i\leq p$, the train track $\tau(i)$ contains $\sigma(i)$ 
as a subtrack. In particular, if $\eta$ is another 
complete train track containing $\sigma$ as a subtrack,
if $\eta^\prime$ is obtained from $\eta$
by a splitting sequence induced from a splitting
sequence connecting $\sigma$ to $\sigma^\prime$ and
if $\{\eta(j)\}$ is the tight multi-sequence
connecting $\eta=\eta(0)$ to $\eta^\prime=\eta(p)$ then
the distance between corresponding points on 
$\{\tau(i)\},\{\eta(j)\}$ is bounded from above
by $Ld(\tau,\eta)+L$ for a 
universal constant $L>0$.

Now for every splitting sequence $\{\tau(i)\}\subset
{\cal V}({\cal T\cal T})$ induced
by a splitting sequence $\sigma(j(i))$ 
of subtracks $\sigma(j(i))<\tau(i)$   
the number of splits $\tau(i)\to \tau(i+1)$ which
reduce the complexity, i.e. such that
$\chi(\tau(i),\sigma(j(i)))>\chi(\tau(i+1),\sigma(j(i+1)))$,
is uniformly bounded. Together with Lemma \ref{flatconfig}
and using the explicit construction of the maps
$\Pi_{E(\tau,\eta)}$ and $\Pi_{E(\tau,\beta)}$ we deduce that
the distance between corresponding points on 
the tight splitting sequences connecting $\tau$ to 
$\beta,\eta$ is uniformly bounded.

Now define $\gamma$ to be the composition of the
inverse of the tight
multi-sequence in $E(\tau,\beta)$ connecting
$\tilde \zeta$ to $\beta$ with the tight multi-sequence
connecting $\zeta$ to $\eta$. The curve $\gamma$ is not
continuous but can be made continuous by inserting
an arc of uniformly bounded length parametrized
on $[0,1]$ which connects $\tilde \zeta$ to $\zeta$. Let
$c_{\sigma,\eta}$ be the constant speed reparametrization
of $\gamma$ on $[0,1]$. By the considerations in Section 4,
this defines indeed a reflexive symmetric
$L$-Lipschitz $L$-quasi-convex bicombing of $X$.
\end{proof}

\section{The geometric rank}

This section is devoted to the proof of Theorem B
from the introduction. We use an argument which
is motivated by the work of Kleiner and Leeb
\cite{KL97}.

Choose again a non-principal ultrafilter $\omega$ and
consider the asymptotic cone
${\cal T\cal T}_\omega$ of ${\cal T\cal T}$ with
respect to $\omega$ and basepoint
the constant sequence $(\tau_0)$ where
$\tau_0$ is a train track in standard form for
a framing $F$ of $S$.
Then ${\cal T\cal T}_\omega$ is
a complete geodesic metric space
(Lemma 2.5.2 of \cite{KL97}).
Since the mapping class group
${\cal M}(S)$ acts properly and cocompactly
as a group of isometries
on ${\cal T\cal T}$,
the asymptotic cone
${\cal T\cal T}_\omega$ is independent of the
point $\tau_0$ and admits
a transitive group of isometries (Proposition 2.5.6
of \cite{KL97}). If we denote by
$X\subset{\cal V}({\cal T\cal T})$ the set of all
complete train tracks which can be obtained from
a train track in standard form for $F$ by a splitting
sequence, equipped with the restriction of the
metric on ${\cal T\cal T}$, then
the cone ${\cal T\cal T}_\omega$ is bilipschitz equivalent to
the asymptotic cone $X_\omega$
of $X$ with respect to $\omega$ and the basepoint
$(\tau_0)$. By Corollary 4.10, for every complete
geodesic lamination $\lambda$ carried by
a train track $\tau\in {\cal V}({\cal T\cal T})$
in standard form for $F$ the inclusion
$E(\tau,\lambda)\to {\cal T\cal T}$ is a quasi-isometric
embedding and hence the
cone ${\cal T\cal T}_\omega$ contains a uniform bilipschitz
image $C(\lambda)$ of the asymptotic cone
of the flat
strip $E(\tau,\lambda)$ with basepoint the constant
sequence $(\tau)$.
Since $E(\tau,\lambda)$ is uniformly quasi-isometric
to its maximal extension $C(\tau,\lambda)$,
the set $C(\lambda)$ is uniformly bilipschitz equivalent
to the asymptotic cone $C(\tau,\lambda)_\omega$
of $C(\tau,\lambda)$. Therefore
by Lemma 5.7, $C(\lambda)$ is locally compact and its
topological dimension is bounded from above
by $3g-3+m$. We call the image of a cone $C(\lambda)$
of this form
under an isometry of ${\cal T\cal T}_\omega$ a \emph{cone}.
Note that each cone in ${\cal T\cal T}_\omega$ is
an ultralimit of a sequence of flat strips
in ${\cal T\cal T}$.

For $k\geq 0$ let $\Delta^k=\{(x_1,\dots,x_{k+1})\in
\mathbb{R}^{k+1}\mid x_i\geq 0,\sum_ix_i= 1\}$
be the standard $k$-simplex in
$\mathbb{R}^{k+1}$. For
$i\geq 0$ let $\Delta^i$ be the subsimplex of $\Delta^k$ which is
the standard face of dimension $i$ obtained by intersecting
$\Delta^k$ with
$\mathbb{R}^{i+1}\subset \mathbb{R}^{k+1}$; we have
$\Delta^0\subset \Delta^1\subset \dots \subset \Delta^k$.
A singular $k$-simplex in a topological space
$Y$ is a continuous map $\sigma:\Delta^k\to
Y$. Denote by $C_*(Y)$ the chain complex of
singular chains in $Y$. For a subset $V$ of $Y$ and a number
$k>0$ let $C_k(Y,V)$ be the set of all
singular $k$-chains
whose boundaries are contained in $V$. We use the results
from Section 6 to show.

\begin{lemma}\label{straightening}
Let $\emptyset\not= V\subset U\subset {\cal T\cal T}_\omega$
be open sets, let $k\geq 1$
and let $\sigma \in C_k(U,V)$ be a singular $k$-chain with
boundary in $V$.
Then there is a singular $k$-chain
$Str(\sigma)\in
C_k(U,V)$ with the following properties.
\begin{enumerate}
\item $Str(\sigma)$ and $\sigma$ define the same class
in $H_k(U,V)$.
\item $Str(\sigma)$ is contained in a finite union of cones.
\end{enumerate}
\end{lemma}

\begin{proof} As before, denote by $X\subset {\cal T\cal T}$
the set of
all train tracks which can be obtained from a train track
in standard form for some fixed framing $F$ of $S$ by
a splitting sequence. We equip $X$ with the restriction of the
metric on ${\cal T\cal T}$. Let $\tau_0$ be a train track
in standard form for $F$. Since $X$ is $r$-dense in
${\cal T\cal T}$ for some $r>0$, for every
non-principal ultrafilter $\omega$ the $\omega$-asymptotic cone
$(X_\omega,d_\omega)$ of $X$ with basepoint
the constant sequence $(\tau_0)$ is bilipschitz
equivalent and hence homeomorphic to the
$\omega$-asymptotic cone of ${\cal T\cal T}$
with basepoint $(\tau_0)$.
In Proposition \ref{bicombing} we constructed for some $L>1$
a reflexive symmetric $L$-Lipschitz $L$-quasi-convex
bicombing of $X$. Taking the $\omega$-limits of the combing
lines defines a reflexive symmetric $L$-Lipschitz
bicombing of $X_\omega$ and hence
${\cal T\cal T}_\omega$ for a possibly different constant
$L>0$. If we denote for $x,y\in X_\omega$
by $c_{x,y}$ the
combing line connecting $x$ to $y$, then this bicombing
is moreover $L$-convex in the
following sense. For every quadruple $x,y,x^\prime,y^\prime$
of points in $X_\omega$ and all $t>0$, we have
$d_\omega(c_{x,y}(t),c_{x^\prime,y^\prime}(t))\leq
L(d_\omega(x,x^\prime)+d_\omega(y,y^\prime))$.
In particular, the combing lines depend continuously on their
endpoints.

To every singular simplex
$\sigma:\Delta^k\to X_\omega$ we associate
a \emph{straightened simplex}
$Str(\sigma)$ inductively as follows. First let
$S^0(\Delta^k)$ be the $0$-skeleton of $\Delta^k$
consisting of $k+1$ vertices and
define $Str(\sigma)(S^0(\Delta^k))=
\sigma(S^0(\Delta^k))$. Assume by
induction that the restriction of $Str(\sigma)$ to
the subsimplex $\Delta^j$ has been defined
for some $j\in \{0,\dots,k-1\}$. Let
$x$ be the vertex of $\Delta^{j+1}$ which is
not contained in $\Delta^{j}$ and extend $Str(\sigma)$ to
$\Delta^{j+1}$ by connecting $Str(\sigma)(x)$ to each point
in $Str(\sigma)(\Delta^j)$ by the combing line with
the same endpoints. By the above observation, this defines
a continuous map $Str(\sigma):\Delta^k\to X_\omega$
which coincides with $\sigma$ on the vertices
of $\Delta^k$.
Since our bicombing is symmetric,
straightening commutes with the boundary
maps. In particular,
the boundary of the straightening of $\sigma$ is the
straightening of the boundary of $\sigma$.
This means that $Str$ defines a chain map
of the chain complex $C_*(Y)$.
Moreover, there is a number
$L(k)>L$ such that if the diameter of
the vertex set
of a singular simplex
$\sigma$ is smaller than some $r>0$ then
the diameter of $Str(\sigma)$ is smaller than
$L(k)r$.

Now let $\emptyset\not=
V\subset U$ be open sets in $X_\omega$ and
let $\sigma\in C_k(U,V)$ be a singular chain.
Since the image of $\sigma$ is compact,
there is a positive lower bound $\delta>0$ for
the distance betwen the image of $\sigma$
and $X_\omega-U$ and for the distance between
the boundary of $\sigma$ and $X_\omega-V$.
After a sufficiently fine barycentric subdivision we may
assume that the diameter
of each simplex in our chain is at most
$\delta/4L^2(k)$.
Then the diameter of each singular simplex in the
straightened chain $Str(\sigma)$ is at most
$\delta/4L(k)$.
Since the $0$-skeleton of $\sigma$
and $Str(\sigma)$ coincide,
the distance between a point $z\in\sigma$
and the corresponding point in the straightening
$Str(\sigma)$ is bounded from above by the sum
of the diameter of $\sigma$ and $Str(\sigma)$ and
hence this distance is
at most $\delta/2L(k)$.
In particular, the
singular chain $Str(\sigma)$
is contained in $C_k(U,V)$.

We claim that $Str(\sigma)$ and $\sigma$ define the
same class in $H_k(U,V)$. Namely, connect each point
in $\sigma$ to the corresponding point in its straightening
by the combing line connecting
these two points. Since the length of a combing
line is bounded from above by $L$ times the distance
between its endpoints, these combing lines are entirely
contained in $U$, and the combing lines which connect
a boundary point of $\sigma$ to a boundary point of $Str(\sigma)$
are entirely contained in $V$.
Thus the collection of these combing lines
define a $k+1$-chain in $C_{k+1}(U,V)$ with boundary
$\sigma- Str(\sigma)$. In other words,
the relative cocycles $\sigma$,
$Str(\sigma)$ are homologous.

Now by construction, each straightened
simplex of the chain $Str(\sigma)$ is contained
in finitely many cones. More precisely, a vertex
$v$ of a singular simplex $\sigma$ can be represented
by a sequence $(x_i)\subset X$. If $w$ is another
vertex which is represented by the sequence
$(y_i)$ then for each $i$ there is a combing
line $c_i$ connecting $x_i$ to $y_i$, and the
combing line $c_{v,w}$ is the $\omega$-limit
of the sequence $(c_i)$. In particular, this
line is contained in the union of the two cones
containing $(x_i),(y_i)$.
As a consequence, the straightened chain $Str(\sigma)$
is contained in finitely many cones as well.
This completes the proof of the lemma.
\end{proof}

The following proposition completes the proof of Theorem B
and of the corollary from the
introduction.

\begin{proposition}\label{geometricrank}
If $k>3g-3+m$ then $H_k(U,V)=0$ for
all pairs of open sets $V\subset U$ in ${\cal T\cal T}_\omega$.
\end{proposition}

\begin{proof} Let $\emptyset\not=
V\subset U\subset {\cal T\cal T}_\omega$ and assume
that $H_k(U,V)\not=0$ for some $k\geq 1$.
Then there is a singular
$k$-chain $c=\sum_i a_i c_i$
for some $a_i\in \mathbb{Z}$ and for
continuous maps $c_i:\Delta^k\to
{\cal T\cal T}_\omega$ whose boundary is contained in $V$
and such that this chain is not homologous to a chain in $V$.
By Lemma \ref{straightening} we may assume without loss of generality
that the chain $c$ is straightened.
This means
in particular that $c$ is contained in a finite
union of cones. These cones are embedded in
${\cal T\cal T}_\omega$ and are glued
along closed subsets.
By Lemma 5.7 the topological dimension
of a cone in $X_\omega$
is not bigger than $3g-3+m$.
In other words, $\sigma$ defines a nontrivial
relative homology class in $H_k(U^\prime,V^\prime)$ where
$U^\prime$ is an open subset in a topological space
obtained from glueing a finite disjoint
collection of standard proper cones of dimension
at most $3g-3+m$ along a finite collection of
closed subsets. But this just means that the topological
dimension of the set $U^\prime$ is at most $3g-3+m$ and hence
we necessarily have $k\leq 3g-3+m$.
This completes the proof of the proposition.
\end{proof}

Now let $k\geq 1$, let $c>1$ and let $\eta:\mathbb{R}^k\to
{\cal T\cal T}$ be a $c$-quasi-isometric embedding with
$\eta(0)=\tau_0$ for our basepoint $\tau_0$.
Then the $\omega$-asymptotic cone
of $\mathbb{R}^k$ admits
a bi-Lipschitz embedding into the asymptotic
cone ${\cal T\cal T}_\omega$ of ${\cal T\cal T}$. Thus
there is a bilipschitz embedding
of $\mathbb{R}^k$ into ${\cal T\cal T}_\omega$.
Since $\mathbb{R}^k$ is an absolute retract,
there are open subsets
$U\supset V$ in ${\cal T\cal T}_\omega$ such that
the relative homology group $H_k(U,V)$ is non-trivial.
By Proposition \ref{geometricrank} this means that $k\leq 3g-3+m$
which shows the corollary from the
introduction.

\begin{corollary}
The geometric rank of
${\cal M}(S)$ equals $3g-3+m$.
\end{corollary}

\section{Quasi-isometric rigidity}

This section is devoted to the proof of Theorem A from
the introduction. We call a
finitely generated group $\Gamma$ \emph{quasi-isometrically
rigid} \cite{M03b} if for every finitely generated group $H$
which is quasi-isometric to $\Gamma$ there is a
finite index subgroup $H^\prime$ of $H$ and a
homomorphism of $H^\prime$ with finite kernel onto a subgroup
of $\Gamma$ of finite index.
Our goal is to show that ${\cal M}(S)$ is quasi-isometrically
rigid; in the case of once-punctured surfaces (i.e. if $m=1$)
this result is due to Mosher and Whyte (see \cite{M03b}).

The \emph{curve graph} ${\cal C}(S)$ of $S$
is the locally infinite metric graph whose
vertices are the free homotopy classes of essential simple
closed curves on $S$, i.e. simple closed curves
which are neither contractible nor freely homotopic
into a puncture, and where two such vertices $c_1,c_2$ are
connected by an edge if and only if the curves $c_1,c_2$
can be realized disjointly. There
is a natural homomorphism $\rho$
from the \emph{extended mapping class
group} of \emph{all} isotopy classes
of homeomorphisms of $S$ into the
group ${\rm Aut}({\cal C}(S))$
of simplicial automorphisms of ${\cal C}(S)$.
By a result of Ivanov (see \cite{I02}) and Luo \cite{L00},
if $S$ is different from the closed surface of genus
$2$ and the twice punctured torus, then
$\rho$ is an isomorphism. If $S$ is a closed surface
of genus $2$ then $\rho$ is surjective, with
kernel the group $\mathbb{Z}/2\mathbb{Z}$
generated by the hyperelliptic involution.
If $S$ is the twice punctured torus, then the kernel
of $\rho$ is again the subgroup $\mathbb{Z}/2\mathbb{Z}$
generated by the hyperelliptic involution, and
the image of $\rho$ is a subgroup of index 5 in
${\rm Aut}({\cal C}(S))$.
As a consequence, for the purpose of our
theorem it suffices to construct
for every finitely generated group
$\Gamma$ which is quasi-isometric to ${\cal M}(S)$
a homomorphism $\rho:\Gamma\to {\rm Aut}({\cal C}(S))$ with
finite kernel and finite index image.

We begin with constructing a
\emph{Tits boundary} ${\cal T\cal B}$ for
${\cal M}(S)$. For this
let $\lambda$ be a complete geodesic lamination
with $k\geq 1$ minimal components $\lambda_1,\dots,\lambda_k$.
After reordering we may assume that for some
$s\leq k$ the laminations $\lambda_1,\dots,\lambda_s$
are minimal arational and the laminations
$\lambda_{s+1},\dots,\lambda_k$ are simple closed geodesics.
For $i\leq k$ the complete lamination
$\lambda$ determines a sign
${\rm sgn}_\lambda(\lambda_i)\in \{+,-\}$ for $\lambda_i$
as follows. If $\lambda_i$ is minimal arational
then we define the sign to be positive. If $\lambda_i$
is a simple closed curve then for a given
orientation of $\lambda_i$, the orientation
of $S$ determines the right and the left side
of $\lambda_i$ in a tubular annulus about $\lambda_i$.
The complete lamination $\lambda$ contains at least one
leaf which spirals about $\lambda_i$ from the left
side. If this spiraling leaf approaches $\lambda_i$
in the direction given by the orientation of $\lambda_i$
then we choose the sign to be positive, otherwise
the sign is chosen to be negative.
Since $\lambda$ is complete by assumption, this choice
of sign does not depend on the
orientation of $\lambda_i$ used to define it (see the discussion
in \cite{H06a}).
We simply write ${\rm sgn}_\lambda$
for this collection
of signs
or also ${\rm sgn}$ if no confusion is possible.

Using the notations from Section 5, for $i\leq s$
denote by $A(\lambda_i)$
the asymptotic cone of a flat strip $C(\zeta,\lambda_i)$
where $\zeta$ is a complete train track on the characteristic
surface $S_i$ for $\lambda_i$ which carries $\lambda_i$.
Recall that $A(\lambda_i)$ is a proper
${\rm CAT}(0)$ cone
defined by a compact ${\rm CAT}(1)$-space
$\partial A(\lambda_i)$, and it does not
depend on $\zeta$ up to uniform bilipschitz identification.
Let again $\Delta^{k-1}=\{(x_1,\dots,x_k)\in \mathbb{R}^k
\mid x_i\geq 0,\sum_ix_i=1\}$ be the standard $k-1$-dimensional
simplex in $\mathbb{R}^k$ and
write $\Delta(\lambda)=\{\sum_i {\rm sgn}_\lambda(\lambda_i) x_i
\mu_i\mid \mu_i\in \partial A(\lambda_i),
(x_1,\dots,x_k)\in\Delta^{k-1}\}$.
The space $\Delta(\lambda)$
is uniquely determined by $\lambda$
up to permutations of the minimal components of $\lambda$
and therefore the topology on
$\Delta(\lambda)$ induced from the topology
on $\Delta^{k-1}$ and the compact ${\rm Cat}(1)$-spaces
$\partial A(\lambda_i)$
is independent of any choices. More
precisely, with this topology
the space $\Delta(\lambda)$ is
homeomorphic to the spherical join
$\partial A(\lambda_1)*\dots *\partial A(\lambda_k)$ of
the spaces $\partial A(\lambda_i)$ $(i\leq k)$
and hence $\Delta(\lambda)$
is homeomorphic to the
boundary $\partial C(\tau,\lambda)_\omega$
of the asymptotic cone of $C(\tau,\lambda)$.
We equip the collection
$\tilde {\cal T\cal B}=\{\Delta(\lambda)\mid
\lambda\in {\cal C\cal L}\}$ with the
topology as a disjoint union of the
spaces $\Delta(\lambda)$, i.e.
$\tilde{\cal T\cal B}$ is naturally a disjoint
union of compact ${\rm Cat}(1)$-spaces of dimension
at most $3g-4+m$ (compare Section 5).

Define an equivalence relation $\sim$ on $\tilde{\cal T\cal B}$
as follows. Let again $\lambda_1,\dots,\lambda_k$
be the minimal components of the complete
geodesic lamination $\lambda$ and let
$\mu_1,\dots,\mu_m$ be the minimal components of
the complete geodesic lamination $\mu$.
After reordering
there is some $\ell\leq \min\{m,k\}$
such that the following holds.
\begin{enumerate}
\item $\mu_i=\lambda_i$ for all
$i\leq \ell$.
\item If $\lambda_i$ is a simple closed curve for some $i\leq \ell$
then ${\rm sgn}_\lambda(\lambda_i)={\rm sgn}_\mu(\lambda_i)$.
\item If there is some $i>\ell$ and a component
$\lambda_i$ of $\lambda$
which coincides with a component
$\mu_j$ of $\mu$ for $j>\ell$ then
${\rm sgn}_\lambda(\lambda_i)\not=
{\rm sgn}_\mu(\mu_j)$.
\end{enumerate}
Then both $\Delta(\lambda)$ and $\Delta(\mu)$ contain
the (signed) spherical join
$\partial A(\lambda_1)*\dots *\partial A(\lambda_\ell)$
as a closed subspace,
and we define
$x\in \Delta(\lambda)$ to be equivalent to
$y\in \Delta(\mu)$ if and only if
$x,y$ are both contained in this signed
spherical join and define the same point there.

By the definition of our topology
on $\tilde {\cal T\cal B}$,
the equivalence relation $\sim$ is closed.
It identifies the spaces $\Delta(\lambda),
\Delta(\mu)$ for all complete geodesic
laminations $\lambda,\mu$ which contain
the same unique minimal component $\lambda_0$
which fills up $S$, i.e.
which is such that every simple closed curve on $S$
intersects $\lambda_0$ transversely.
Note that the number of such spaces
which are identified in this way with a fixed
space $\Delta(\lambda)$ is bounded
from above by a universal constant.
The quotient space ${\cal T\cal B}=\tilde {\cal T\cal B}/\sim$
has naturally the structure of a (locally infinite)
complex of dimension
$3g-4+m$; it is closely related to
the curve graph on $S$ (compare \cite{MM99}
for a discussion of the curve complex). Note that we do not
claim that ${\cal T\cal B}$ is a cell complex in the usual sense.
Moreover, ${\cal T\cal B}$ has infinitely
many connected components. Namely, every minimal
geodesic lamination $\nu$ which fills up $S$ defines
a connected component of ${\cal T\cal B}$ which is
homeomorphic to a single compact ${\rm Cat}(1)$-space.
However, since
the curve graph is connected, the
components of $\tilde {\cal T\cal B}$
defined by geodesic laminations with
a minimal
component which does not fill up $S$
all map to the same connected component
${\cal T\cal B}_0$ of ${\cal T\cal B}$.

We call the image in ${\cal T\cal B}$ of a set
$\Delta(\lambda)\subset \tilde {\cal T\cal B}$
a \emph{cell}.
Every spread out geodesic lamination on $S$ defines a
cell which is naturally homeomorphic
to a standard $3g-4+m$-dimensional
simplex, and we call these cells \emph{chambers}.
The vertices of a chamber in ${\cal T\cal B}$
either correspond to signed simple closed geodesics
or to minimal arational geodesic laminations which fill a
bordered subtorus or an $X$-piece
of $S$, i.e. a bordered punctured sphere
of Euler characteristic $-2$.
Denote by ${\cal T\cal B}_1$ the subcomplex
of ${\cal T\cal B}$ which is the closure of the
chambers. The complex ${\cal T\cal B}_1$ is
a connected simplicial complex in the usual sense.

If $S$ is not a closed surface
of genus $2$ or a twice punctured torus then we define
${\cal M}_0(S)$ to be the extended mapping class
group of all isotopy classes of
\emph{any} homeomorphism of $S$ including the
orientation reversing ones.
For a closed surface $S$ of genus
$2$ or a twice punctured torus
we define ${\cal M}_0(S)$ to be the
quotient of the extended mapping class group under
the hyperelliptic involution. Then ${\cal M}_0(S)$
naturally acts on ${\cal T\cal B}$ as
a group homeomorphisms preserving the cell structure and
the subcomplex ${\cal T\cal B}_1$.

As before,
denote by ${\rm Aut}({\cal C}(S))$ the
group of simplicial automorphisms of
the curve graph ${\cal C}(S)$ of $S$.
The automorphism group ${\rm Aut}({\cal C}(S))$ naturally
contains the group ${\cal M}_0(S)$ as a subgroup.
In fact, equality holds if $S$ is not a twice
punctured torus. Using this fact, the next
lemma gives a description of the group of
isotopy classes of homeomorphisms of ${\cal T\cal B}$.

\begin{lemma}\label{boundarymap}
\begin{enumerate}
\item Every homeomorphism of ${\cal T\cal B}$
preserves the subcomplex ${\cal T\cal B}_1$.
\item
There is an injective homomorphism of
the group of isotopy classes
of homeomorphisms
of ${\cal T\cal B}_1$ into ${\rm Aut}({\cal C}(S))$
whose restriction to ${\cal M}_0(S)$ is the identity.
\end{enumerate}
\end{lemma}

\begin{proof} Let $\phi$ by an arbitrary
homeomorphism of ${\cal T\cal B}$. Since ${\cal T\cal B}$
is the disjoint
union of the locally infinite connected
subcomplex ${\cal T\cal B}_0$ and infinitely
many compact components, $\phi$ preserves ${\cal T\cal B}_0$.
We claim that $\phi$ also preserves ${\cal T\cal B}_1$.

For this let $\lambda$ be
a complete geodesic lamination with minimal
components $\lambda_1,\dots,\lambda_k$ which
defines a cell $\Delta(\lambda)$ in ${\cal T\cal B}$
of dimension $3g-4+m$ (by
abuse of notation we use now the same
symbol for a the space $\Delta(\lambda)\subset \tilde {\cal T\cal B}$
and its image cell in ${\cal T\cal B}$). After reordering,
there is a number $s\geq 0$ such that the
components $\lambda_1,\dots,\lambda_s$ are precisely
those simple closed curve components
of $\lambda$ which are contained in the boundary
of a characteristic subsurface $S_j$ of
a minimal arational component $\lambda_j$ of $\lambda$.
Then for every $i\leq s$,
\emph{every} complete
geodesic lamination $\mu$ which defines
a cell $\Delta(\mu)$ of dimension $3g-4+m$
and which contains each of the minimal components
$\lambda_j$ for $j\not=i$ also contains $\lambda_i$.

Recall that a point in $\Delta(\lambda)$
defines a tuple $(s_1,\dots,s_k)\in \Delta^{k-1}$
and a tuple $(x_1,\dots,x_k)\in
\partial A(\lambda_1)\times\dots\times \partial A(\lambda_k)$.
By the above consideration,
for each $i\leq s$ the set of all points in $\Delta(\lambda)$
corresponding to a tuple $(s_1,\dots,s_k)\in \Delta^{k-1}$
with $s_i=0$ and $s_j>0$ for $j\not=i$ is contained
in the boundary of precisely two cells of maximal dimension,
namely one cell
for each choice of a sign for $\lambda_i$.
If $x$ has a neighborhood in $\Delta(\lambda)$
which is homeomorphic to a closed half-space in
$\mathbb{R}^{3g-4+m}$ containing $x$ in its boundary
(which is always the case if $\Delta(\lambda)$ is a chamber),
then $x$ has a neighborhood in ${\cal T\cal B}_0$
which is homeomorphic to
$\mathbb{R}^{3g-4+m}$. Call the union
of all cells in ${\cal T}{\cal B}_0$ of maximal dimension
whose intersection with $\Delta(\lambda)$
is of this form the \emph{star} of $\Delta(\lambda)$.
Define moreover
inductively the \emph{multi-cell} $R(\lambda)$
containing $\Delta(\lambda)$
to be the smallest subcomplex of ${\cal T\cal B}_0$
which contains $\Delta(\lambda)$ and which contains
with each cell
its star. Note that a multi-cell consists of a uniformly bounded
number of cells.
Namely, if $\mu\in {\cal C\cal L}$ defines a cell contained
in the multi-cell $R(\lambda)$ then the minimal components of
$\mu$ coincide with the minimal components of $\lambda$.
If the star of a cell coincides with
the cell itself then we also
call this single cell a multi-cell.

We claim that a point $x\in {\cal T\cal B}_0$
which is not contained
in the interior of a multi-cell
does not have a neighborhood in ${\cal T\cal B}_0$ which
is homeomorphic to an open subset of $\mathbb{R}^{3g-4+m}$.
Namely, let $x\in {\cal T\cal B}_0$ be any
point admitting a neighborhood in ${\cal T\cal B}_0$
which is homeomorphic to an open subset of $\mathbb{R}^{3g-4+m}$.
If $x$ is \emph{not} an interior point
of a multi-cell in ${\cal T\cal B}_0$, i.e. if
no neighborhood of $x$ in ${\cal T\cal L}_0$
is entirely contained in a multi-cell,
then $x$ is necessarily contained
in the boundary of a cell
of maximal dimension.
Thus assume that $x$ is contained in the boundary of the cell
$\Delta(\lambda)$ for some $\lambda\in {\cal C\cal L}$.
Let $\lambda_1,\dots,\lambda_k$ be the
minimal components of $\lambda$.
Then up to reordering,
the point $x$ defines a tuple $(s_1,\dots,s_k)\in \Delta^{k-1}$
and a tuple $(x_1,\dots,x_k)\in
\partial A(\lambda_1)\times\dots\times \partial A(\lambda_k)$
where $s_1=0$.

Next assume that $\lambda_1$ is a minimal
arational geodesic lamination with characteristic
surface $S_1$.
Since the dimension of $\Delta(\lambda)$ equals
$3g-4+m$ by assumption, there are infinitely
many pants decompositions of $S_1$
whose union with $\cup_{j\geq 2} \lambda_j$
determine a cell in ${\cal T\cal B}_0$ of maximal
dimension.
Then $x$ is contained in the boundary of
infinitely many such cells and hence
a neighborhood of $x$ in ${\cal T\cal B}_0$ is not
locally compact.
As a consequence, necessarily $\lambda_1$ is a simple
closed geodesic.

If  $\lambda_1$ is a simple closed geodesic
which is not contained in the boundary of
the characteristic surface of a minimal
arational component $\lambda_i\not=\lambda_1$ of
$\lambda$ then there is a bordered subsurface $S_1$ of
$S$ of negative Euler characteristic
containing $\lambda_1$ which
does not have an essential intersection with $\cup_{i\geq 2}\lambda_i$.
Thus there
are infinitely many pairwise distinct simple closed
geodesics which do not intersect
$\cup_{i\geq 2}\lambda_i$ and which
define together with $\lambda_2,\dots,\lambda_k$ a
cell of maximal dimension. As a consequence,
$x$ is contained in infinitely
many distinct such cells
which contradicts our assumption that there is a
neighborhood of $x$ in ${\cal T\cal B}_0$
which is homeomorphic to
a ball in $\mathbb{R}^{3g-4+m}$.
This shows that indeed a point $x\in {\cal T\cal B}_0$
has a neighborhood in
${\cal T\cal B}_0$ which is homeomorphic to
$\mathbb{R}^{3g-4+m}$ only if $x$ is
an interior point of a multi-cell.
In particular, every homeomorphism of ${\cal T\cal B}_0$
preserves the union of all interior points of all
such multi-cells in ${\cal T\cal B}_0$.

Our next goal is to give a topological characterization
of multi-cells defined by spread out geodesic laminations.
Thus let $\lambda$ be a spread out geodesic
lamination with minimal components
$\lambda_1,\dots,\lambda_{3g-3+m}$
which defines the multi-cell
$R(\lambda)$. Then $R(\lambda)$ has a canonical structure
of a finite simplicial complex. The boundary
of $R(\lambda)$ is partitioned into finitely many
sides, i.e. simplices of codimension one.
Such a side either is defined by a
minimal arational component of $\lambda$
or by a simple closed geodesic which is not
a boundary component of the characteristic surface
of a minimal arational minimal component of $\lambda$. Thus
a neighborhood in ${\cal T\cal B}$ of an interior point $x$
of such a side contains a set which is
homeomorphic to a countable collection of closed
half-spaces in $\mathbb{R}^{3g-4+m}$
glued along their boundary, and $x$ is contained
in the interior of the boundary of these
half-spaces. In particular, such a neighborhood
is not homeomorphic to the neighborhood of
a point in a side of a chamber of codimension at least 2.
In other words, there is an open dense subset $V$ of the
boundary of the multi-cell $R(\lambda)$
which consists of $p>0$ open disjoint sides,
whose closure is
the whole boundary of $R(\lambda)$ and which
admits a purely topological
characterization. On the other hand, if $R(\lambda)$ is
a multi-cell and
if $\lambda$ is \emph{not} spread out,
then $R(\lambda)$ contains points which admit
a compact neighborhood in $R(\lambda)$ not
containing any neighborhood which is homeomorphic
to an open subset of $\mathbb{R}^{3g-3+m}$.
But this just means that
multi-cells defined by chambers in ${\cal T\cal B}$ admit
a topological characterization and hence
the image of a multi-cell defined by
a chamber under a homeomorphism
$\phi$ of ${\cal T\cal B}$ is again a
multi-cell defined by a chamber.
As a consequence, every homeomorphism
$\phi$ of ${\cal T\cal B}$
restricts to a homeomorphism of
the simplicial complex ${\cal T\cal B}_1$
which is homotopic to a simplicial map.

Call a vertex $v$ of ${\cal T\cal B}$ \emph{simple}
if $v$ is defined by a signed simple closed curve.
For a vertex $a$
of ${\cal T\cal B}_1$ let ${\cal F}(a)$ be the collection
of all multi-cells containing $a$; we claim that $a$ is
simple if and only if for every $k\geq 2$ there
are multi-cells $F_1,\dots,F_k\in {\cal F}(a)$
with $F_i\cap F_j=\{a\}$ for all $i\not=j$.
Namely, if $a$ is simple then for every $k\geq 1$ we
can find $k$ pants
decompositions $P_1,\dots,P_k$ of $S$ containing $a$ as one
of their pants curves and such that every pants
curve $c\in P_i-\{a\}$ intersects at least one of the
curves in the collection $P_j-\{a\}$ transversely.
As a consequence, if $C_i,C_j$ are two chambers
in ${\cal T\cal B}_1$ defined by the pants decompositions
$P_i,P_j$ and a system of signs which coincide
for the vertex $a$ then we have $C_i\cap C_j=\{a\}$ by
the definition of our complex ${\cal T\cal B}$.

On the other hand, if $a$ is a vertex
of the simplicial complex ${\cal T\cal B}_1$ which is not
simple then $a$ is defined by
a minimal arational geodesic lamination which
fills a non-trivial subsurface $S_0$ of
$S$. More precisely, there is a unique bordered
subsurface $S_0$ of $S$ with boundary
$\partial S_0\not=\emptyset$
and with the additional property that every
simple closed curve on $S$ which has an essential
intersection with $\partial S_0$ also has an essential
intersection with $a$. Since $a$ is a vertex
of a chamber, the surface $S_0$ either is a once
punctured torus or a forth punctured sphere.
Therefore every multi-cell of ${\cal T\cal B}_1$ which
contains $a$ as a vertex also contains each
component of the boundary of
$S_0$ equipped with a choice of a sign. Since the boundary
of $S_0$ consists of at most $4$ components,
if $F_1,\dots,F_{2^{4}+1}$
are $2^{4}+1$ multi-cells containing
the vertex $a$ then at least two of them, say the
multi-cells $C_1,C_2$, contain the boundary
$\partial S_0$ of $S_0$ equipped with the same
collection of signs and hence the intersection
$C_1\cap C_2$ contains an edge. As a consequence,
simple vertices can be distinguished from non-simple
ones by the multi-cells they are contained in and therefore
a simplicial automorphism
of ${\cal T\cal B}_1$ has to preserve the set
of simple vertices.

Recall that every simple vertex of ${\cal T\cal B}_1$
consists of a simple closed curve on $S$ together with
a sign. In other words, a simple closed curve
$a$ gives rise to two distinct vertices $(a,+),(a,-)$ which
differ by their sign. We claim that a homeomorphism
of ${\cal T\cal B}_1$ preserves the set of pairs
$((a,+),(a,-))$ of such vertices. For this we argue
as before. Namely, let $a$ be a simple
closed geodesic and let $(a,+)$ be any simple vertex
of ${\cal T\cal B}$ defined by $a$
(here we mean that $+$ is the sign of our vertex).
Let $C\subset{\cal T\cal B}_1$ be
any chamber containing $(a,+)$. Then
the vertices of the side $F$ of $C$ opposite
to $a$ correspond to
a geodesic lamination with $3g-4+m$ minimal components.
The side $F$ is also a side of a chamber which contains
the simple vertex $(a,-)$. In other words,
for \emph{every} chamber $C_+$ containing $(a,+)$ as a vertex,
the point $(a,-)$ is a vertex of a chamber $C_-$
with the same opposite side.

Now let $C$ be a chamber
which contains a vertex $x$ defined by a minimal
arational component filling up
a subsurface which contains the simple closed curve
$a$ in its boundary; note that
such a chamber always exists.
Then $a$ is contained in the boundary of a multi-cell
which is not a chamber, and this multi-cell contains
both vertices $(a,+),(a,-)$. The multi-cell
is distinguished from a chamber by the
number of its sides, and the intersection of all
these multi-cells consists precisely of the
two points $(a,+),(a,-)$.
This shows that the
pair of vertices $(a,+),(a,-)$ is characterized
by its intersection pattern with the chambers
of ${\cal T\cal B}_1$. Therefore
every simplicial automorphism of ${\cal T\cal B}_1$
preserves the pairs of simple vertices
determined by a simple closed geodesic on $S$. Thus
such an automorphism $\phi$ induces a bijection of the
$0$-skeleton of the
curve graph ${\cal C}(S)$ of $S$.

This bijection preserves
disjointness for simple closed curves. Namely,
two simple closed curves on $S$ can be
realized disjointly if and only if they
define two distinct vertices of a common chamber.
In other words, a homeomorphism $\phi$ of
${\cal T\cal B}_1$ defines a simplicial automorphism
of ${\cal C}(S)$. This completes the proof of our
lemma.
\end{proof}

Let ${\cal T\cal C}$ be the (simplicial) cone over the
Tits boundary ${\cal T\cal B}$ (compare
\cite{KL97}). Define a
\emph{chamber} in ${\cal T\cal C}$ to
be the cone over a chamber in ${\cal T\cal B}$.
Note that a chamber is a
$3g-3+m$-dimensional standard proper cone.
Define moreover an
\emph{apartment} in ${\cal T\cal C}$ to
be a union of chambers which is homeomorphic to
$\mathbb{R}^{3g-3+m}$. Every pants decomposition $P$ of $S$
defines $2^{3g-3+m}$
chambers, one for each choice of sign combination
for our pants curves, and the union of all
these chambers defines an apartment.
The cone over a cell in ${\cal T\cal B}$
is contained in an apartment if and
only if it is a chamber.
Namely, we saw in the proof of Lemma \ref{boundarymap} that
a cell in ${\cal T\cal B}$
of maximal dimension which is not a chamber
contains points which are not contained in
subsets of ${\cal T\cal B}$ homeomorphic to
$\mathbb{R}^{3g-4+m}$ and hence
such a cell can not be contained in
an apartment. On the other hand, if
$\lambda$ is a spread out complete geodesic lamination
with minimal components $\lambda_1,\dots,\lambda_{3g-3+m}$
and minimal arational components
$\lambda_1,\dots,\lambda_s$ then for each
$i\leq s$ we can find a simple closed
curve $c_i$ which does not intersect
$\lambda_j$ for $j\not=i$ and such that the
curves $c_1,\dots,c_s$ are disjoint.
In particular, for every subset
$\{\lambda_{i_1},\dots,\lambda_{i_\ell}\}$ of
the set $\{\lambda_1,\dots,\lambda_s\}$ for some $\ell\leq s$
we obtain
a new chamber by replacing each of the
components $\lambda_{i_j}$ of $\lambda$
by the signed simple
closed curve $(c_{i_j},+)$ $(j\leq \ell)$.
A signed simple closed curve component $\lambda_i$ of
$\lambda$ $(i\geq s+1)$
can be replaced by the same curve with the opposite sign.
The union of the chambers defined by $\lambda$
with all those chambers obtained from all
possible combinations of such replacements
defines an apartment by construction.

Choose again a non-principal ultrafilter $\omega$.
Fix a framing $F$ for $S$ and let $X$ be the
set of all train tracks which can be obtained from
a train track in standard form for $F$ by a splitting
sequence. Choose a train track
$\tau_0$ in standard form for $F$ and
use the constant sequence $(\tau_0)$
as a basepoint for the
asymptotic cone ${\cal T\cal T}_\omega$.
Recall that for every complete geodesic
lamination $\lambda\in {\cal C\cal L}$
the flat strip $E(F,\lambda)$
is quasi-isometric to its
maximal extension $C(F,\lambda)$ which
is a ${\rm Cat}(0)$-space, moreover $E(F,\lambda)$
is quasi-isometrically embedded in
${\cal T\cal T}$. As a consequence, the $\omega$-asymptotic
cone $C(F,\lambda)_\omega$
of $C(F,\lambda)$ with basepoint a constant sequence is
topologically embedded in the $\omega$-asymptotic cone
${\cal T\cal T}_\omega$ of ${\cal T\cal T}$.
The cone $C(F,\lambda)_\omega$ is homeomorphic to the
euclidean cone over the cell $\Delta(\lambda)$.
If $\lambda_1,\dots,\lambda_k$ are the minimal
components of $\lambda$ then $C(F,\lambda)_\omega=
\partial A(\lambda_1)*\dots*\partial A(\lambda_k)$
where for each $i$, $\partial A(\lambda_i)$ is a
compact connected ${\rm Cat}(1)$-space.

If $\Pi_{E(F,\lambda)}:X\to E(F,\lambda)$
is the projection as in Proposition \ref{shortestdistance} then
for every complete geodesic lamination $\nu$ the
projection $\Pi_{E(F,\lambda)}
E(F,\nu)$ of $E(F,\nu)$ is a combinatorially
convex subset of $E(F,\lambda)$
and hence its maximal extension is
a complete ${\rm Cat}(0)$-space. The
Hausdorff distance in ${\cal T\cal T}$ between the graph
$\Pi_{E(F,\nu)}E(F,\lambda)$ and the graph
$\Pi_{E(F,\lambda)}E(F,\nu)$ is
uniformly bounded. The asymptotic cone
$E(F,\lambda)_\omega$ of $E(F,\lambda)$ contains
the asymptotic cone of the projection
$\Pi_{E(F,\lambda)}E(F,\nu)$. Let
$\lambda_1,\dots,\lambda_s$ be those minimal
components of $\lambda$
which are minimal components of $\nu$ as well
and with the additional property that for
every minimal component $\lambda_i$ which
is a simple closed curve, the sign ${\rm sgn}_\lambda(\lambda_i)$
defined as above by $\lambda$ for $\lambda_i$ coincides
with the sign ${\rm sgn}_\nu(\lambda_i)$
defined by $\nu$ for $\lambda_i$.
By Lemma \ref{topologicaldimension},
the asymptotic cone
of the projection $\Pi_{E(F,\lambda)}E(F,\nu)$
is the subset of $E(F,\lambda)_\omega$ which
is the cone over $\partial A(\lambda_1)*\dots *\partial A(\lambda_s)$.
As a consequence,
our Tits cone ${\cal T\cal C}$ is topologically
embedded in the asymptotic cone ${\cal T\cal T}_\omega$
of ${\cal T\cal T}$.

As in Section 6.3 of \cite{KL97},
for a point $z\in Z$ we say that two subsets
$S_1,S_2$ of $Z$ \emph{have the same germ at $z$} if
$S_1\cap U=S_2\cap U$ for some neighborhood $U$ of $z$.
The equivalence classes of subsets with the same germ
at $z$ will be denoted ${\rm Germ}_zZ$.
Write $k=3g-3+m$ for short.
For $x\in {\cal T\cal T}_\omega$ consider the collection
${\cal S}_1(x)$ of all germs of topological embeddings
of $\mathbb{R}^k$ into ${\cal T\cal T}_\omega$
passing through $x$. Note that each
such germ determines a local homology
class of degree $k$
whose support contains $x$. Let ${\cal S}_2(x)$ be the
lattice of germs generated by ${\cal S}_1(x)$ under
finite intersection and union. The following lemma is
the analog of Lemma 6.3.1 of \cite{KL97}.

\begin{lemma}\label{lattice}
The lattice ${\cal S}_2(x)$ admits a natural
embedding into the lattice ${\cal K}{\cal T\cal C}$
of the Tits cone ${\cal T\cal C}$
for ${\cal T\cal T}_\omega$ generated
by the cells of maximal dimension under finite
union. The image of this embedding
contains the sublattice of the Tits cone ${\cal T\cal C}$
generated by the chambers under finite
intersection and union.
\end{lemma}

\begin{proof}
Following Kleiner and Leeb (Section 6 of \cite{KL97}),
for a subset $Y$ of a topological space
$Z$ and for $[c]\in H_k(Z,Y)$
define ${\rm Supp}(Z,Y,[c])\subset Z-Y$
to be the set of points $z\in Z-Y$ such that the image
of $[c]$ in the local homology group
$H_k(Z,Z-\{z\})$ is nonzero.
Then ${\rm Supp}(Z,Y,[c])$ is a closed
subset of $Z-Y$ contained in the image of
any chain $c$ representing the relative class $[c]$.

Since the isometry group of ${\cal T\cal T}_\omega$
acts transitively, it is sufficient to show the
claim of the lemma for the lattice ${\cal S}_2(*)$ defined
by the basepoint $*$ of ${\cal T\cal T}_\omega$.

By Lemma 6.2.1 of \cite{KL97}, every germ
of a topological embedding of $\mathbb{R}^k$ through
$*$ defines a nontrivial class in $H_k({\cal T\cal T}_\omega,
{\cal T\cal T}_\omega)-\{*\}$, i.e. there is an
open subset $U$ of ${\cal T\cal T}_\omega-\{*\}$
and a singular chain $c\in C_k({\cal T\cal T}_\omega,U)$
whose support contains $*$.

Recall from Section 7 the construction
of straightening which associates to a singular chain $c$
representing $[c]$ the straightened chain $Str(c)$.
Via passing to a
sufficiently small barycentric subdivision
we may assume that the boundary of the
straightening ${Str}(c)$ of $c$ is contained in $U$.
By Lemma \ref{straightening}, the chain
$Str(c)$ is contained in a finite union
${\cal P}=\cup_{i=1}^\ell P_i$
of cones $P_i$ which intersect along their
boundaries. Now the support of the class
$[c]\in H_k({\cal T\cal T}_\omega,U)$ is contained
in the image of $Str(c)$, on the other hand
it is defined by the germ of our embedding of
$\mathbb{R}^k$ into ${\cal T\cal T}_\omega$. Since
the dimension of each cone is at most $k$, elements
of ${\cal S}_1(*)$ define finite unions of
cells of maximal dimension in ${\cal T\cal C}$.

On the other hand, by the discussion in
the beginning of this section,
each chamber is a finite intersection of
apartments and hence of
elements of ${\cal S}_1(x)$. Intersections of chambers
yield sides of the Tits cone, so we have a well defined
inclusion of lattices $\Theta:{\cal S}_2(x)\to {\cal K}{\cal T\cal C}$
containing the sublattice generated by
the apartments under finite union and intersection.
\end{proof}

If $\phi:{\cal T\cal T}_\omega\to {\cal T\cal T}_\omega$
is any homeomorphism which
fixes the basepoint $*=(\tau_0)$ then $\phi$
induces a homeomorphism of lattices ${\cal S}_2(*)\to
{\cal S}_2(*)$ and therefore by
Lemma \ref{boundarymap} and Lemma \ref{lattice},
this homeomorphism
defines an element of ${\rm Aut}({\cal C}(S))$. We state this
fact as a corollary.

\begin{corollary}\label{homomorphism1} There is a homomorphism
from the group of homeomorphisms
of ${\cal T\cal T}_\omega$ which fix the basepoint $*$
into the group
${\rm Aut}({\cal C}(S))$.
\end{corollary}

\begin{proof} If $\phi$ is a homeomorphism of
${\cal T\cal T}_\omega$ fixing $*$ then
by our above discussion, $\phi$ determines
a homeomorphism $\tilde \phi$ of subcone of the
Tits cone ${\cal T\cal C}$ containing
the simplicial cone over the complex
${\cal T\cal B}_1$. We then
associate to $\phi$ the element $\rho(\phi)\in
{\rm Aut}({\cal C}(S))$ whose restriction to the cone over
${\cal T\cal B}_1$
coincides with the restriction of $\tilde\phi$
up to an isotopy preserving
the simplicial structure of ${\cal T\cal B}_1$.
\end{proof}

Now we are ready to complete the proof of Theorem A from
the introduction.
Namely, let $\Gamma$ be a finitely generated group
with a word norm $\vert\, \vert$ defined
by a finite symmetric set of generators.
The norm $\vert \,\vert $ defines a distance
function $d$ on $\Gamma$ via $d(g,h)=\vert g^{-1}h\vert$
which is invariant under left translation.
With respect to this distance, the group $\Gamma$
acts on itself as a group of isometries by
left translation. Assume that $\Gamma$ is
quasi-isometric to ${\cal M}(S)$, i.e. that
there is a quasi-isometry $\Theta_0:\Gamma\to {\cal M}(S)$.
Since ${\cal M}(S)$ is equivariantly
quasi-isometric to ${\cal T\cal T}$ \cite{H06a}
and hence it is quasi-isometric to $X$, there is
a quasi-isometry $\Theta:\Gamma\to X\subset
{\cal V}({\cal T\cal T})$ with
$\Theta(e)=\tau_0$ (here $e$ denotes the unit) and
inverse $\Lambda:X\to \Gamma$ where
$\tau_0$ is a train track in standard form
for the framing $F$ defining $X$. Via this
quasi-isometry, the group $\Gamma$ induces a
quasi-action as a group of uniform quasi-isometries
on $X$ as follows. The quasi-isometry
determined by $g\in \Gamma$ is the
map $\phi(g)$ defined by
$\phi(g)(\eta)=\Theta\circ g\circ \Lambda(\eta)$ where $g$
acts on $\Gamma$ by left translation. By construction,
there is a universal constant $L>0$ such that
$d(\phi(g)\phi(h)(\eta),\phi(gh)(x))\leq L$
for all $g,h\in \Gamma$ and all $x\in X$
(compare the discussion in \cite{M03b}).
The quasi-action of $\Gamma$ on $X$
then induces an action of $\Gamma$
as a group of uniformly bilipschitz
homeomorphisms on the asymptotic cone
${\cal T\cal T}_\omega=X_\omega$. Since
the basepoint $*=(\tau_0)$ is the
ultralimit of both the constant sequence $(\tau_0)$ and
the constant sequence $(\Theta\circ g\circ \Lambda(\tau_0))$,
this action preserves the basepoint $*=(\tau_0)$.
By Corollary \ref{homomorphism1} there is a homomorphism
of the group of homeomorphisms
of ${\cal T\cal T}_{\omega}$ preserving $*$
into the group ${\rm Aut}({\cal C}(S))$
and therefore we obtain a homomorphism $\Gamma\to
{\rm Aut}({\cal C}(S))$. We summarize our discussion as follows.

\begin{lemma}\label{homomorphism2}
Let $\Gamma$ be quasi-isometric
to ${\cal M}(S)$; then there is a homomorphism
$\rho:\Gamma\to {\rm Aut}({\cal C}(S))$.
\end{lemma}

For the proof of the theorem
in the introduction we are left with showing that
the kernel of our homomorphism is finite and that its image
is of finite index.

We follow again \cite{KL97}. Namely,
the following result is the analog
of Proposition 7.1.1 of \cite{KL97}. For its formulation,
for a constant $D>0$ define a \emph{$D$-Hausdorff envelope}
of a set $A\subset X$ to be a set $B$ containing
$A$ in its $D$-neighborhood.
Define a \emph{maximal
quasi-flat} in ${\cal T\cal T}$ to be
the image under an element of ${\cal M}(S)$
of a finite union of flat strips
which is uniformly quasi-isometric to $\mathbb{R}^{3g-3+m}$.
Then the asymptotic
cone of such a maximal quasi-flat with respect to a basepoint defined
by a constant sequence
is an apartment in the Tits cone ${\cal T\cal C}$.
Define more generally an \emph{apartment} (or
a \emph{chamber}) in the
asymptotic cone ${\cal T\cal T}_\omega$ to be
the image of an apartment (or a chamber)
in the Tits cone ${\cal T\cal C}$
under an isometry of ${\cal T\cal T}_\omega$.
Then an apartment consists
of finitely many chambers. A chamber
is bilipschitz equivalent to a standard partition
cone of dimension $3g-3+m$. This
cone is determined by its boundary which admits
a natural identification with
a compact convex subset in the standard unit sphere
$S^{3g-4+m}$ with dense interior.
For a number $\delta >0$
define a \emph{$\delta$-truncated chamber}
to be the closed subcone of a chamber $A$ defined by the compact
convex subset of the boundary $\partial A$ of
$A$ consisting of all points whose
spherical distance to the boundary of $\partial A$ in 
$S^{3g-4+m}$ is at least $\delta$.
A \emph{$\delta$-truncated apartment} is obtained from an
apartment by replacing a chamber by the
$\delta$-truncated chamber contained in its interior.
We define moreover a \emph{$\delta$-truncated
maximal quasi-flat} in ${\cal T\cal T}$ to be a subset
of a maximal quasi-flat $Y$ whose asymptotic cone is
the $\delta$-truncated apartment contained in the
asymptotic cone of $Y$. We have (compare \cite{KL97}).

\begin{proposition}\label{approximation} Let $\cal Q$ be a family
of subsets of ${\cal T\cal T}$ which are uniformly
quasi-isometric to $\mathbb{R}^{3g-3+m}$.
Then for every $\delta >0$ there is
a constant $D=D(\delta)>0$ so that any set $Q\in {\cal Q}$
is a $D$-Hausdorff envelope of a $\delta$-truncated
maximal quasi-flat.
\end{proposition}

Before we show the proposition we use it to derive
the theorem from the introduction.

\begin{corollary}\label{final}
If $\Gamma$ is quasi-isometric
to ${\cal M}(S)$ then there is a homomorphism
$\Gamma\to {\rm Aut}({\cal C}(S))$ with finite kernel
and finite index image.
\end{corollary}

\begin{proof} We observed above that a quasi-isometry
$\Theta_0:\Gamma\to {\cal M}(S)$ of a
finitely generated group $\Gamma$ into ${\cal M}(S)$
induces a quasi-isometry $\Theta:\Gamma\to X$ with
inverse $\Lambda:X\to \Gamma$. We
may assume that $\Theta$ maps the identity $e$ in $\Gamma$
to a train track $\tau_0$ in standard form for the framing $F$
as above.
We showed above that there is a homomorphism
$\rho:\Gamma\to {\rm Aut}({\cal C}(S))$ obtained from the fact
that $\Gamma$ acts as a group of homeomorphisms on the
asymptotic cone ${\cal T\cal T}_\omega$ preserving the
basepoint $(\tau_0)$.
The homeomorphism induced by an
element $h\in \Gamma$ is the ultralimit of the map
$\phi(h)=\Theta\circ h\circ \Lambda$.

We have to show that the kernel of this homomorphism is finite.
For this assume to the contrary that the
kernel is an infinite subgroup $H$ of $\Gamma$. Let
$d$ be the distance in $\Gamma$ induced by a word norm; then
there is a sequence of elements $h_i\in H$
$(i\geq 0)$ with
$d(e,h_i)\to\infty$ $(i\to \infty)$.
Let ${\cal Q}_0\subset {\cal T\cal T}$ be a finite
collection of maximal quasi-flats in ${\cal T\cal T}$
defined by a finite collection ${\cal P}=\{P_1,\dots,P_s\}$
of pants decompositions of $S$.
By this we mean that for each $Q\in {\cal Q}_0$
there is some $j\leq s$ such that $Q$ is a union
of flat strips $E(\tau_j^\ell,\lambda_j^\ell)$
where $\ell\leq 2^{3g-3+m}$,
where $\tau_j^\ell$ is a complete train track
in standard form for $P_j$,
where the complete geodesic laminations $\lambda_j^\ell$
have $P_j$ as the union of their minimal components and where
the signs defined by $\lambda_j^\ell$ for the pants curves
of $P_j$ run through all possible sign combinations.
We require
that for every $i$ and every pants curve
$\gamma\in P_i$ there is some $j\not=i$ and some
$\tilde \gamma\in P_j$ which intersects $\gamma$ transversely.
We may assume that for every $D>0$ the intersection of
the $D$-neighborhoods of these quasi-flats is a compact
neighborhood of $\tau_0$ in ${\cal T\cal T}$.

Let ${\cal Q}=\cup \phi(h_i){\cal Q}_0$; since the maps
$\phi(h_i)=\Theta\circ h_i\circ \Lambda$
are uniform quasi-isometries of ${\cal T\cal T}$, the family
${\cal Q}$ satisfies the assumptions in
Proposition \ref{approximation}. Thus
by Proposition \ref{approximation}, for a small number $\delta>0$
there is a constant $D>0$ such that
each of the sets $\phi(h_i)Q$ $(Q\in {\cal Q}_0,i>0)$ is
a $D$-Hausdorff envelope of a $\delta$-truncated
maximal quasi-flat $F(h_i,Q)$ 
in ${\cal T\cal T}$. 

Let $A(h_i,Q)$ be the asymptotic
cone of $\phi(h_i)Q$ with basepoint the constant
sequence $\phi(h_i)\tau_0=*$.
Then $A(h_i,Q)$ is a
bilipschitz-embedded euclidean space of dimension $3g-3+m$
passing through $*$ which contains the asymptotic cone of the
$\delta$-truncated maximal quasi-flat $F(h_i,Q)$ with
basepoint
$(\phi(h_i)\tau_0)=*$. As a consequence, 
$A(h_i,Q)$ equals the unique apartment which contains 
the $\delta$-truncated quasi-flat $F(h_i,Q)$. In particular, 
this apartment contains the
basepoint $*$ and hence it is an apartment in the
Tits cone ${\cal T\cal C}$.
However, by assumption the elements
$h_i$ are contained in the kernel of the homomorphism
$\rho$ and therefore $A(h_i,Q)=Q_\omega$ for all $i$ where
$Q_\omega$ is the asymptotic cone of $Q$ with basepoint
a constant sequence. As a consequence, $\phi(h_i)$ maps
each of the maximal quasi-flats $Q\in {\cal Q}_0$ to a set
containing the $\delta$-truncated maximal quasi-flat 
$F(h_i,Q==Q_\delta\subset
Q$ in its $D$-neighborhood for
a universal constant $D >0$.

After possibly increasing $D$, the $D$-neighborhoods of the
maximal truncated quasi-flats
$Q_\delta\subset Q\in {\cal Q}_0$ intersect,
and this intersection is contained
in a uniformly bounded
neighborhood of $e$. Therefore for each $i$ the $2D$-neighborhoods of
the sets $\phi(h_i)(Q)$ $(Q\in {\cal Q}_0)$ contain
intersection points in a uniformly bounded
neighborhood of $e$ which is independent of $i$.
On the other hand, the maps $\phi(h_i)$ are $L$-quasi-isometries
for a universal constant $L>0$ and hence
the intersection of the $2D$-neighborhoods
of the images $\phi(h_i)(Q)$ $(Q\in {\cal Q}_0)$ is contained
in a uniformly bounded neighborhood of the image
under $\phi(h_i)$
of the intersections of the $2LD$-neighborhoods
of the sets $Q\in {\cal Q}_0$. By the choice of our family
${\cal Q}_0$, this implies
that the distance between
$\phi(h_i)(\tau_0)$
and $\tau_0$ is bounded from above by a universal constant
not depending on $i$.
Now the map $\Theta:\Gamma\to X$
is a quasi-isometry with inverse $\Lambda$ and therefore
we conclude that the distance between $h_i=h_i(\Lambda\tau_0)$ 
and $e=\Lambda(\tau_0)$ in $\Gamma$
is uniformly bounded from above as well. This is a
contradiction and shows that
the kernel $H$ of our homomorphism $\rho$
is a finite subgroup of
$\Gamma$.

Our above argument shows
that for every $h\in \Gamma$, the
distance between $\phi(h)(\tau_0)$ and $\Theta(h)$ is
uniformly bounded.
From this we deduce that the image of the homomorphism
$\rho$ is of finite
index in ${\cal M}(S)$.
Namely, if this is not the case
then there is a sequence
$\{g_i\}\subset {\cal M}(S)$ with $d(g_i,\rho(\Gamma))\to \infty$.
Since the quasi-isometry $\Theta_0:\Gamma\to
{\cal M}(S)$ is coarsely surjective there is
a sequence $\{h_i\}\subset\Gamma$ with
$d(\Theta_0(h_i),g_i)\leq c$ for a fixed constant $c>0$.
This means that
$d(\Theta_0(h_i),\phi(h_i)(\tau_0))\to \infty$ which is
impossible by our above argument.
This completes the proof of our corollary
and hence the proof of the theorem from the introduction.
\end{proof}

We are left with the proof of Proposition \ref{approximation}.
For this we follow Section 7 of \cite{KL97}. Namely,
let ${\cal Q}$ be a family of uniformly quasi-isometrically
embedded euclidean spaces $\mathbb{R}^{3g-3+m}$ in $X$ 
as in Proposition \ref{approximation}. Consider first a single
set $Q$ from the family ${\cal Q}$ and choose a basepoint
$q\in Q$. By assumption, the ultralimit
$\omega-\lim(\frac{1}{n}Q,q)$ is an apartment $A(Q)$
in the asymptotic cone $\omega-\lim(\frac{1}{n}X,q)$
which contains the basepoint $*=(q)$. Such
an apartment consists of a finite union of chambers
which meet at the chamber walls.

Let $\overline{* x_\omega}$ be a line segment
contained in the $\delta$-truncated
apartment $A_\delta\subset A(Q)$
and issuing from the basepoint
$*$. Note that a chamber has naturally
the structure of a standard cone, so this is well defined.
Note also that by the definition of a chamber,
the line segment $\overline{*x_\omega}$ 
coincides with the combing line
connecting $*$ to $x_\omega$. There is a
sequence $(x_n)\subset Q$ of points such that
$\omega-\lim(x_n)=x_\omega$. Since $x_n\in X$ for all $n$,
for every $n$ there is a train track in standard form
for $F$ which is splittable to $x_n$. For $\omega$-almost
all $n$ this train track coincides with a fixed
track $\tau_0$, so we may assume that $\tau_0$ is splittable
to $x_n$ for all $n$. Then for each $n$ the flat
strip $E(\tau_0,x_n)$ is defined. We view $E(\tau_0,x_n)$
as a subset of its ${\rm Cat}(0)$-extension
$C(\tau_0,x_n)$.

Let $d_n\leq cn$ for all $n$ and a fixed number
$c>0$ and assume that
$d_n\to \infty$ $(n\to\infty)$. Then
the ultralimit $\omega-\lim(\frac{1}{d_n}Q,q)$ is an
apartment in ${\cal T\cal T}_\omega$ which
necessarily coincides with $A(Q)$.
If $\omega-\lim d_n/n\to 0$
then the ultralimit of the geodesic arcs in $C(\tau_0,x_n)$
connecting $\tau_0$ to $x_n$ is a
geodesic ray in the asymptotic
cone $C=\omega-\lim_{n\to\infty}\frac{1}{d_n}C(\tau_0,x_n)$
which defines a point $\zeta$ in the Tits boundary of
${\cal T\cal T}_\omega$. This point coincides with
the point in the boundary of the apartment 
$A(Q)$ defined by the unique
geodesic extension of the line segment $\overline{*x_\omega}$.

{\bf Sublemma:} There is a number
$r>0$ so that for $\omega$-all $n$ the sets
$E(\tau_0,x_n)$ are contained in the tubular $r$-neighborhood of $Q$.

Choose a point $z_n\in E(\tau_0,x_n)$
at maximal distance $d_n$ from $Q$. Note that
there is a universal constant $c>0$ such that $d_n\leq cn$.
We argue by contradiction and we assume that
$\omega-\lim d_n=\infty$.
Since $E(\tau_0,x_n)$ is connected, via possibly
changing $z_n$ we may assume that
$\omega-\lim d_n/n=0$. The asymptotic
cone $\omega-\lim(\frac{1}{d_n}{\cal T\cal T},\tau_0)$
contains the cell $E=
\omega-\lim_n\frac{1}{d_n}E(\tau_0,x_n)$ and the apartment
$\omega-\lim_n\frac{1}{d_n}Q=A(Q)$.
The point $z_\omega=(z_n)$ is contained in $E$
but not in $A(Q)$ and therefore
$E$ is not contained in $A(Q)$.
We may assume that there is a point $x_\omega^\prime\in
\partial_\infty A(Q)$ in the asymptotic
boundary of $A(Q)$ obtained from an
$\omega$-limit of  
combing lines in $Q$ connecting $*$ to $x_n$.
However, this limit of lines is a combing segment in
${\cal T\cal T}_\omega$ connecting $*$ to $x_\omega$ 
and hence is necessarily contained in
$E$ as well. As a consequence, the flat
strips $E,A(Q)$ intersect
in the ray. Since the ray is regular, it is contained
in the interior of a unique regular chamber and hence
the chamber $E$ is contained in $A(Q)$,
a contradiction which shows the sublemma.

Now since the point $x_\omega$
in the $\delta$-truncated apartment
$A_\delta$ was
arbitrary we conclude as in \cite{KL97} 
that there is a number $R>0$ such that
each $Q\in {\cal Q}$ contains the truncated maximal
quasi-flat defining the truncated
apartment $A_\delta$ in its $R$-neighborhood.
This completes the proof of the proposition.
\qed

\bigskip

\noindent
MATHEMATISCHES INSTITUT DER UNIVERSIT\"AT BONN\\
BERINGSTRASSE 1, D-53115 BONN, GERMANY


\begin{thebibliography}{CEG87}

\bibitem[B95]{B95} W.~Ballmann, {\em Lectures
on spaces of nonpositive curvature}, DMV Seminar,
Band 25, Birkh\"auser, Basel 1995.

\bibitem[BM05]{BM05} J.~Behrstock, Y.~Minsky, {\em Dimension
and rank for mapping class groups}, arXiv:math.GT/0512352.



\bibitem[BH99]{BH99} M.~Bridson, A.~Haefliger, {\sl Metric
spaces of non-positive curvature}, Springer Grund\-leh\-ren 319,
Springer, Berlin 1999.




\bibitem[CEG87]{CEG87} R.~Canary, D.~Epstein, P.~Green,
{\em Notes on notes of Thurston}, in ``Analytical and geometric
aspects of hyperbolic space'', edited by D.~Epstein, London Math.
Soc. Lecture Notes 111, Cambridge University Press, Cambridge 1987.

\bibitem[CB88]{CB88} A.~Casson with S.~Bleiler, {\sl Automorphisms
of surfaces after Nielsen and Thurston}, Cambridge University
Press, Cambridge 1988.


\bibitem[FLM01]{FLM01} B.~Farb, A.~Lubotzky, Y.~Minsky, {\em Rank one phenomena
for mapping class groups}, Duke Math. J. 106 (2001), 581-597.

\bibitem[FLP91]{FLP91} A.~Fathi, F.~Laudenbach, V.~Po\'enaru, {\sl Travaux de
Thurston sur les surfaces,} Ast\'erisque 1991.

\bibitem[G87]{G87} M.~Gromov, {\em Hyperbolic groups},
Essays in group theory (S.M.~Gersten, editor),
MSRI Publications 8, Springer 1987, 75--263.

\bibitem[G93]{G93} M.~Gromov, {\em Asymptotic invariants
of infinite groups}, in ``Geometric group theory'',
London Math. Soc. Lecture Notes 182 (1993).



\bibitem[H06a]{H06a} U.~Hamenst\"adt, {\em Geometry of the
mapping class groups I: Boundary amenability},
revised version, March 2006, arXiv:math.GR/0510116.

\bibitem[H06b]{H06b} U.~Hamenst\"adt, {\em
Geometry of the mapping class groups II: (Quasi)-geodesics},
revised version, May 2006, arXiv:math.GT/0511349.

\bibitem[Har86]{Har86} J.~Harer, {\it The virtual cohomological
dimension of the mapping class group of an oriented surface},
Invent. Math. 84 (1986), 157-176.

\bibitem[Hat91]{Hat91} A.~Hatcher, {\it On triangulations of
surfaces}, Topology Appl. 40 (1991), 189-194.


\bibitem[I02]{I02} N.~V.~Ivanov, {\sl Mapping class groups},
Chapter 12 in Handbook of Geometric Topology (Editors
R.J.~Daverman and R.B.~Sher), Elsevier Science (2002), 523-633.


\bibitem[K06]{K06} Y.~Kida, {\em Measure equivalence
rigidity for the mapping class group}, preprint 2006.

\bibitem[K99]{K99} B.~Kleiner, {\em The local
structure of length spaces with curvature bounded above},
Math. Zeit. 231 (1999), 409--456.

\bibitem[KL97]{KL97} B.~Kleiner, B.~Leeb, {\em
Rigidity of quasi-isometries for symmetric spaces and Euclidean
buildings}, Inst. Hautes Etudes Sci. Publ. Math.  No. 86 (1997),
115--197.


\bibitem[L00]{L00} F.~Luo, {\em Automorphisms
of the complex of curves}, Topology 39 (2000), 283--298.



\bibitem[MM99]{MM99} H.~Masur, Y.~Minsky, {\em Geometry of the
complex of curves I: Hyperbolicity}, Invent. Math. 138 (1999),
103-149.


\bibitem[MM00]{MM00} H.~Masur, Y.~Minsky, {\em Geometry of the complex
of curves II: Hierarchical structure}, GAFA 10 (2000), 902-974.





\bibitem[M03a]{M03} L.~Mosher, {\em Train track expansions of measured
foliations}, unpublished manuscript.

\bibitem[M03b]{M03b} L.~Mosher, {\em Homology and dynamics
in quasi-isometric rigidity of once-punctured mapping
class groups}, preprint 2003.

\bibitem[O96]{O96} J.~P.~Otal, {\sl Le Th\'{e}or\`{e}me d'hyperbolisation
pour les vari\'et\'es fibr\'ees de dimension 3}, Ast\'erisque 235,
Soc. Math. Fr. 1996.





\bibitem[PH92]{PH92} R.~Penner with J.~Harer, {\sl Combinatorics
of train tracks}, Ann. Math. Studies 125, Princeton University
Press, Princeton 1992.

\bibitem[T79]{T79} W.~Thurston, {\em Three-dimensional geometry
and topology}, unpublished manuscript.


\end{thebibliography}
\end{document}